\documentclass{article}

\usepackage{cite}
\usepackage[utf8x]{inputenc}
\usepackage{amsmath}
\usepackage{amssymb}
\usepackage{amsthm}
\usepackage{graphicx}
\usepackage[colorinlistoftodos]{todonotes}
\usepackage[colorlinks=true, allcolors=blue]{hyperref}
\usepackage{caption}
\usepackage{subcaption}
\usepackage{epstopdf} 
\usepackage{authblk}
\usepackage{bbm}

\newtheorem{Thm}{Theorem}

\newtheorem{cor}{Corollary}

\theoremstyle{definition}

\newtheorem{expl}{Example}
\newtheorem{rem}{Remark}

\newcommand{\R}{\mathbb{R}}
\newcommand{\T}{\mathbb{T}}
\newcommand{\X}{\mathbb{X}}

\newcommand{\prob}{\mathbb{P}}
\newcommand{\XT}{(\boldsymbol{X}(t))_{t\in \T}}
\newcommand{\YT}{(\boldsymbol{Y}(t))_{t\in \T}}

\newcommand{\NT}{(\boldsymbol{N}(t))_{t\in \T}}
\newcommand{\hatnt}{(\boldsymbol{\hat{N}}(t))_{t\in \T}}
\newcommand{\hatn}{\boldsymbol{\hat{N}}}

\newcommand{\bN}{\boldsymbol{N}}

\newcommand{\propensity}{f}
\newcommand{\rate}{\gamma}

\newcommand{\hG}{\mathcal{G}}

\newcommand{\hcG}{\hat{\mathcal{G}}}
\newcommand{\hcL}{\hat{\mathcal{L}}}

\newcommand{\revision}[1]{\textcolor{black}{#1}}

\title{\revision{Mathematical modeling of \\ spatio-temporal population dynamics \\ and application to epidemic spreading}}

\author[a]{Stefanie Winkelmann\footnote{winkelmann@zib.de}}
\author[a,b]{Johannes Zonker}
\author[a,b]{Christof Sch\" utte}
\author[a]{Nata\v sa Djurdjevac Conrad}

\affil[a]{Zuse Institute Berlin, 14195 Berlin, Germany}
\affil[b]{Freie Universit\" at Berlin, Institut f\"ur Mathematik und Informatik, 14195 Berlin, Germany}

\begin{document}

\maketitle

\begin{abstract}
Agent based models (ABMs) are a useful tool for modeling spatio-temporal population dynamics, where many details can be included in the model description. Their computational cost though is very high and for stochastic ABMs a lot of individual simulations are required to sample quantities of interest. 
Especially, large numbers of agents render the sampling infeasible. 
Model reduction to a metapopulation model leads to a significant gain in computational efficiency, while preserving important dynamical properties. Based on a precise mathematical description of spatio-temporal ABMs, we present two different metapopulation approaches (stochastic and piecewise deterministic) and discuss the approximation steps between the different models within this framework. Especially, we show how the stochastic metapopulation model results from a Galerkin projection of the underlying ABM onto a finite-dimensional ansatz space. Finally, we
 utilize our modeling framework to provide a conceptual model for the spreading of COVID-19 that can be scaled to real-world scenarios.

\end{abstract}


\section{Introduction}

Recently, epidemiological models have received lots of attention due to the COVID-19 pandemics. State of the art models for analyzing the epidemics spreading and doing forecasts include agent-based models (ABMs) \cite{muller2020realistic,nagelpaper,Goldenbogen2020,Wulkow2020}, meta\-population models  \cite{CoronaMetapopulationKeeling2020,chinazzi2020effect,Ajelli2010} and deterministic models based on ordinary differential equations (ODE) \cite{kermack1927contribution,hethcote2000mathematics,alvarez2020simple,roda2020difficult,shao2020dynamic,chen2020time,ivorra2020mathematical}. While ABMs have a spatial resolution on a microscale, capturing the spatial movement and interactions of each individual agent, meta\-population models  distinguish between \textit{subpopulations} as groups of agents which are homogeneously mixed. Both of these two models can take into account stochastic effects in the interaction dynamics and the spatial exchange. In contrast, ODE models are purely deterministic and assume homogeneous mixing in the overall space of motion, not involving any mobility dynamics. 
Although ABMs often reflect the reality better than metapopulation and ODE-based models, they are very costly to simulate and to re-calibrate to new parameter sets, especially for large population numbers that are considered when modeling a pandemic.  This motivates  \revision{the  consideration of} model reduction techniques which allow to keep the characteristic properties of the detailed ABM, but significantly reduce the computational effort for simulation and analysis. 
In this work we present a step-wise model reduction which decreases both spatial resolution and stochasticity of an ABM to a level appropriate for  fast and accurate simulations of real-world systems. 

In the ABM, which we will introduce as a first modeling approach, each individual agent owns a \textit{position} in some given spatial environment, as well as a \textit{status} which can, e.g., refer to the agent's state of health, opinion or level of knowledge \cite{conrad2018mathematical}. In the course of time, the agents move within the environment and interact with each other, thereby possibly adopting the status of other agents in their spatial vicinity. The propensity for such adoption events to occur is defined by a \textit{adoption rate function}.  Simulations follow the movement and interaction of each individual agent and produce a high computational burden.   

In order to approximate these detailed dynamics by the less complex \textit{stochastic metapopulation model (SMM)}, we assume that the spatial environment splits up into areas between which mixing is rare, e.g., areas of metastability where the transition rate to other areas is small.  Groups of agents that are located in the same area are termed \textit{subpopulations}. 
In application to epidemic dynamics, one can think of these subpopulations as inhabitants of a city or country, but also as smaller groups of people, e.g. pupils of the same school, age classes, working groups or households.
Within each of the subpopulations, individuals of the same status are considered as indistinguishable. An agent can anytime interact with all other agents of the same subpopulation, while an interaction with agents of a different subpopulation requires a preceding spatial transition to the associated spatial area. Both the interactions within and the transitions between the subpopulations are treated as stochastic events inducing a jump in the \textit{population state}, which is defined by the number of agents of each status in the different subpopulations. Analytically, this stochastic metapopulation model results from a Galerkin projection of the underlying ABM based on a spatial partition into metastable subsets \cite{Djurdjevac2012,sarich2010approximation,winkelmann2016spatiotemporal}. In this work, we will define the corresponding projection operator and derive explicit formulas for the matrix representations of the projected diffusion and interaction operators for both first- and second-order status-changes. We will analyze the approximation quality for a simple example. 
Simulations of the  Markov jump process described by the SMM can be done by means of the event-based stochastic simulation algorithm \cite{gillespie1977exact}. This is much more efficient than the agent-based simulations, but the effort still scales linearly with the population size. 

Given that the number of agents in each subpopulation is large, a further approximation of the dynamics is possible, replacing the purely stochastic jump dynamics by a \textit{piecewise-deterministic Markov processes}. This approach, discussed for, e.g., chemical reaction networks \cite{bookC03}, has well-established mathematical foundations, going back to Kurtz \cite{kurtz1971limit,kurtz1972relationship}, who studied under which conditions the stochastic dynamics of a large population are well approximated by deterministic processes. 
In the resulting hybrid modeling approach, the interactions within each subpopulation are described by continuous, deterministic dynamics in form of ODEs, while the comparatively rare exchange between the subpopulations are still treated as stochastic, discrete events. These piecewise-deterministic dynamics can be efficiently simulated by a combination of the stochastic simulation algorithm and an ODE solver \cite{vestergaard2015temporal,conrad2018human}, and the effort is independent of the population size. 

There exist several approaches to describe epidemic (or other types of) population dynamics by semi-deterministic processes 
    \cite{abboud2018piecewise,soubeyrand2009spatiotemporal}. 
    A piecewise-deterministic metapopulation model (PDMM), where deterministic dynamics within the subpopulations  are combined with stochastic transitions between the subpopulations, is formulated in \cite{MONTAGNON2020} for modeling a cattle trade network. We will here use a similar hybrid approach in order to model and analyze epidemic spreading kinetics.
        This efficient PDMM-approach will allow to easily simulate and compare several scenarios induced by various measures that are taken to control the spreading of the COVID-19 within  a population. We will investigate the interrelation of the kinetics within the different subpopulations and analyze the critical transition time of the virus spreading between separate subpopulations depending on several choices of containment measures. 

    \begin{figure}[h!]
    \centering
    \includegraphics[width=0.7\textwidth]{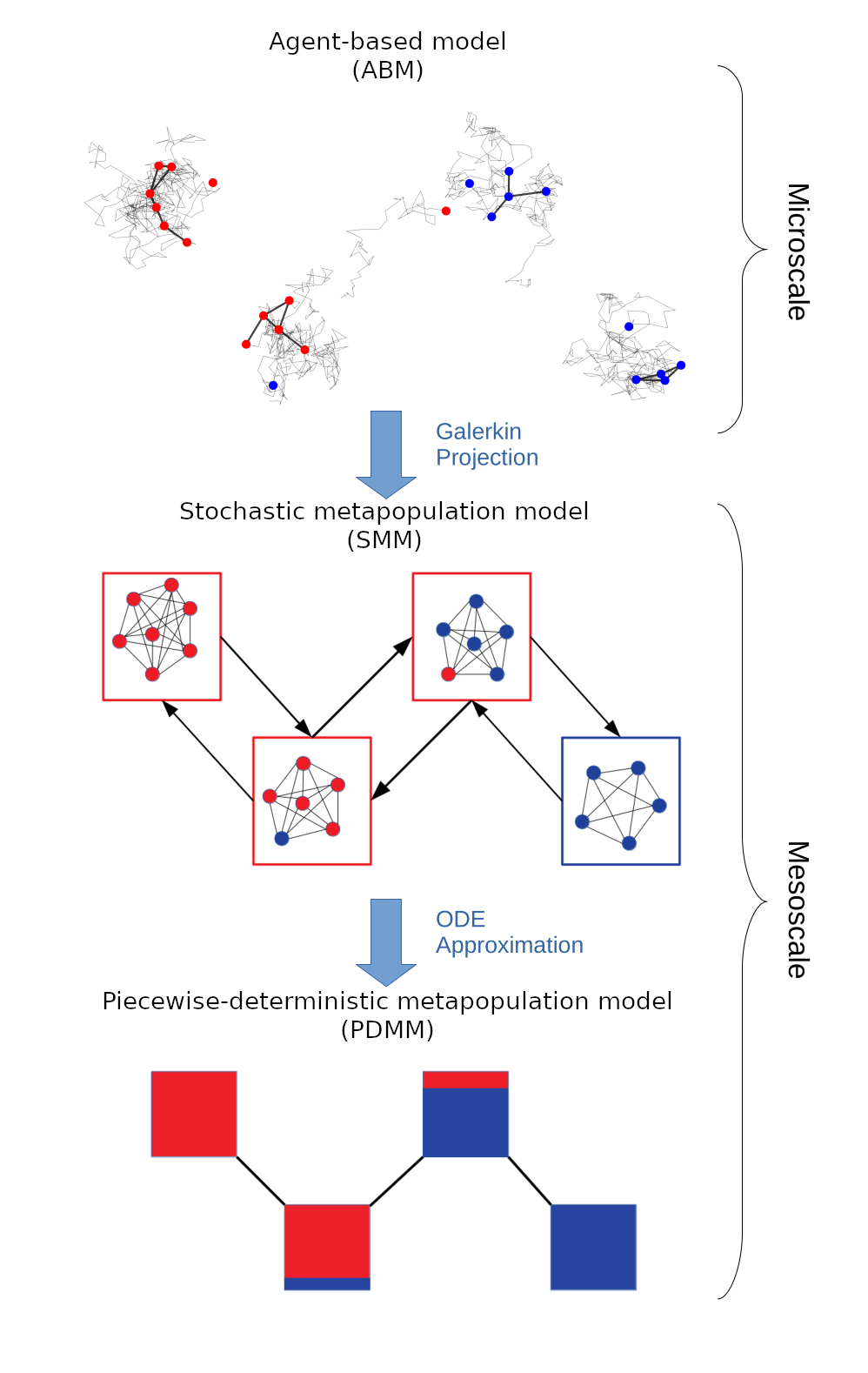}
    \caption{
    Illustration of hierarchy of modeling approaches considered herein and the related steps of the coarse-graining 
    procedure. Dots refer to agents that have a status (blue or red), boxes on the mesoscale represent subpopulations. Straight lines between dots indicate connections between agents that enable interaction. Grey lines on the microscale show movement trajectories of the agents, arrows and lines between boxes on the mesoscale represent possible transitions between subpopulations. For the piecewise-deterministic metapopulation model, the population state is depicted by the proportion of blue and red.
    }
    \label{Fig:Overview}
    \end{figure}

\paragraph{Outline.} We first introduce the mathematical formulations for each of the three modeling approaches and describe the details of the coarse-graining steps in Section \ref{sec:models}, see Figure \ref{Fig:Overview} for an overview. This is illustrated by a guiding example for which we also analyse the approximation quality and highlight the gain in computational efficiency.
Moreover, we derive explicit results for the relation between the \textit{microscopic} adoption rate functions, defining the propensity for an agent to change status in the ABM, and the \textit{macroscopic} adoption rate functions as an analogue in the SMM, see Section \ref{sec:MetaPop}.
Finally, we apply the PDMM framework in Section \ref{sec:studies} to model the spreading of the COVID-19 virus.

\section{Hierarchy of modeling approaches} \label{sec:models}

In this section, we present three basic modeling concepts for interacting agent dynamics and show how they are related to each other. We start in Section \ref{sec:ABM} with the most detailed \textit{agent-based model}, which tracks the movement and interactions of each individual agent present in the system. Assuming that within certain regions of the environment the spatial mixing is fast compared to the interaction dynamics, the dynamics may be approximated by a \textit{stochastic metapopulation model}  which is introduced in Section  \ref{sec:MetaPop}. A partial approximation of the stochastic dynamics by deterministic dynamics, which becomes reasonable in cases of large populations, leads to a \textit{piecewise-deterministic metapopulation process} which is considered in Section  \ref{sec:PDMM}.

\subsection{Agent-based model} \label{sec:ABM}
We consider a set of $n_a$ agents that are moving in a compact domain $\X\subset\R^d$ as independent realizations\footnote{Instead of many independent processes for each agent movement we consider one process that governs the development of all movements. The movements of specific agents are related to the respective marginal processes of $\XT$.} of a diffusion process $\XT$ on $\X^{n_a}$, where $\mathbb{T}$ is a continuous time interval, e.g. $\mathbb{T}=[0,\infty)$. 
In addition to the position $x_\alpha\in \X$, each agent $\alpha$ also has a status $s_\alpha \in\mathbb{S}:=\{1,...,n_s\}$. Let $X=(x_\alpha)_{\alpha=1,...,n_a}$ denote the vector of all agents positions, while $S=(s_\alpha)_{\alpha=1,...,n_a}$ is the vector of all agents' status and $Y=(X,S)$ is the system state. We define the space of all possible system states as $\mathbb{Y}:= \X^{n_a}\times\mathbb{S}^{n_a}$.

The dynamics in the agents' status are modelled by a continuous-time Markov jump process. Each agent can change its status by transitions of the form $ i \mapsto j$
for $i,j \in \mathbb{S}$, called \textit{status-change/adoption events}. We will distinguish between (i) \textit{first-order adoptions} given by status transitions  that require no interactions between agents and (ii) \textit{second-order adoptions} based on pairwise interactions with close neighbors. 
The rate for agent $\alpha$ to undergo such a transition may in general depend on the whole system state $Y=(X,S)$ (i.e. the positions and status of all other agents),\footnote{More generally, the rate function can also be time-dependent.}  
and is given by the adoption rate function $$\propensity_{i j}^{(\alpha)}:\mathbb{Y}\rightarrow [0,\infty). $$ 
Given that the agents are positioned according to $X$ and have a status according to $S$, 
it is $\propensity_{i j}^{(\alpha)}(X,S)$ the probability per unit of time for agent $\alpha$ to switch from status $i$ to status $j$. 
Note that the adoption rate function $\propensity_{i j}^{(\alpha)}$ does not explicitly depend on $\alpha$, i.e., it is not the case that every agent has its own propensity function. Dependence on $\alpha$ is only indirect,  with the propensity depending on the status and position of the agent. More concretely, as for first-order adoptions, we consider rate functions of the form
 \begin{equation} \label{f_first-order}
    f_{ij}^{(\alpha)}(X,S) = \delta_i(s_\alpha) \rate_{ij}(x_\alpha), 
 \end{equation}
where $\delta_i$ denotes the indicator function of a status, i.e. $\delta_i(s)=1$ if $s=i$ and $\delta_i(s)=0$ otherwise, and $\rate_{ij}:\X \to [0,\infty)$ gives the rate for a status transition from $i$ to $j$ depending on the spatial location of the acting agent. In particular, we set $\rate_{ii}(x)=0$ for all positions $x$.  

Equivalently, for second-order adoptions, we set 
\begin{equation} \label{f_second-order}
    f_{ij}^{(\alpha)}(X,S) = \delta_i(s_\alpha)\sum_{\substack{\beta=1\\ \beta\neq\alpha}}^{n_a}\delta_j(s_\beta) \rate_{ij}(x_\alpha,x_\beta),
\end{equation}
where $\rate_{ij}:\X^2 \to [0,\infty)$ is a function defining the rate for a status adoption from $i$ to $j$ depending on the positions of two interacting agents. 
In a more special setting,
this rate $\rate_{ij}(x_\alpha,x_\beta)$  depends only on the distance between the interacting agents, e.g., they need to be closer than some \textit{interaction radius} $r>0$ to interact, as in the Doi-model \cite{Doi} in the context of chemical reaction systems. The underlying idea is that interactions of agents (which may induce  adoptions of the status from other agents) require proximity of the agents in physical space. 
For this case, we set
\begin{equation} \label{Doi}
    \rate_{ij}(x_\alpha,x_\beta) := c_{ij}\cdot d_r(x_\alpha,x_\beta)
\end{equation}
 for a constant $c_{ij}\geq 0$, where  $d_r:\X^2 \to \{0,1\}$ for $r > 0$ is the distance indicator function:
 \begin{equation}\label{dr}
    d_r(x_\alpha,x_\beta) = \left\{
\begin{array}{ll}
1,\quad & \mathrm{if} \; |x_\alpha-x_\beta| \revision{\leq} r \\
0, & \mathrm{otherwise}.
\end{array}
\right. 
 \end{equation}
 For $i=j$ there are no status transitions and we set $\revision{c}_{ii}=0$.

The coupling of the diffusion process  and the  jump dynamics given by status-changes leads to a Markov process $\YT$ on the system state space $\mathbb{Y}$. 
Let $p(X,S,t)$ be the probability mass function for the process $\YT$ to be in the system state $(X,S)$ at time $t$, where the marginal with respect to $X$ is a continuous density function.

We define for each status $i$ 
an operator $L_{i}$ that describes the change of the probability mass function through the motion of a single agent  under the condition that he is currently in status $i$. Then we can write down an operator $L$ for the movement of the agents, 
\begin{equation}\label{def:L}
    L p(X,S,t):=\sum_{\alpha =1}^{n_a} L^{(\alpha)}_{s_\alpha}p(X,S,t),
\end{equation}
where $L_{s_\alpha}^{(\alpha)}$ is defined as $L_i$ for $i=s_\alpha$ acting on a function $p(X)$ with respect to the component $x_\alpha$ of $X$ (see Example \ref{ex:guiding} for details). Note that, consequently,   $L_{s_\alpha}^{(\alpha)}$ acts only on the position part of the probability mass function $p(X,S,t)$ in \eqref{def:L}.

As for the adoption dynamics, we define 
\begin{align}\begin{split}\label{def:G}
G p(X,S,t):= &-\sum_{i,j=1}^{n_s}\sum_{\alpha=1}^{n_a}\propensity_{i j}^{(\alpha)}(X,S)p(X,S,t)\\
				  &+\sum_{i,j=1}^{n_s}\sum_{\alpha=1}^{n_a}\propensity_{ij}^{(\alpha)}(X,S+i e_{\alpha}-j e_{\alpha})p(X,S+i e_{\alpha}-j e_{\alpha},t),\end{split}
\end{align}
where $e_{\alpha}$ denotes the $\alpha$th unit vector of $\R^{n_a}$.\footnote{\revision{Note that in the second line of \eqref{def:G} it might hold $(X,S+i e_{\alpha}-j e_{\alpha})\notin \mathbb{Y}$ for some $i,j, \alpha$. These terms, however, are multiplied by $\propensity_{ij}^{(\alpha)}(X,S+i e_{\alpha}-j e_{\alpha})=0 $. Still, for the sake of completeness, we extend the definition of $p$ and set $p(X,S,t):=0$ for $(X,S)\notin \mathbb{Y}$.}} 
The first term on the right-hand side refers to  the outflow from the current state through adoption events and the second term to the inflow through adoption events that would lead to the current state. 
The change of $p(X,S,t)$ including movement and status transitions is then given by the set of differential equations
\begin{equation}\label{eq:ABM}
    \partial_t p(X,S,t)= L p(X,S,t) + Gp(X,S,t)
\end{equation}

The strength of this general agent-based approach is that  different types of restrictions regarding the interaction dynamics can be included and quite complicated dynamics can be formulated. 
On the other hand, the coupled differential equations which describe the system are usually not analytically solvable. Instead, Monte Carlo (MC) simulations of the dynamics are required to sample the quantities we are interested in, and the simulation of such complex systems can be numerically very costly, especially if local neighborhoods have to be computed in every time step. We realize the simulations with an algorithm that is a combination of the stochastic simulation algorithm, where we draw the time for the next adoption event, and a numerical scheme for stochastic differential equations (e.g. Euler-Maruyama) to update the positions of agents between the adoption events \cite{conrad2018human}. The computational cost can be very high, especially for small time steps and large agent numbers. As soon as there are second-order interactions, the effort does not scale  linearly with the number of agents. This motivates \revision{us} to consider approximate modeling approaches which reduce the numerical effort. 

\begin{figure}[h!]
    \centering
    \begin{subfigure}{0.5\textwidth}
    \centering
        \includegraphics[width=\linewidth]{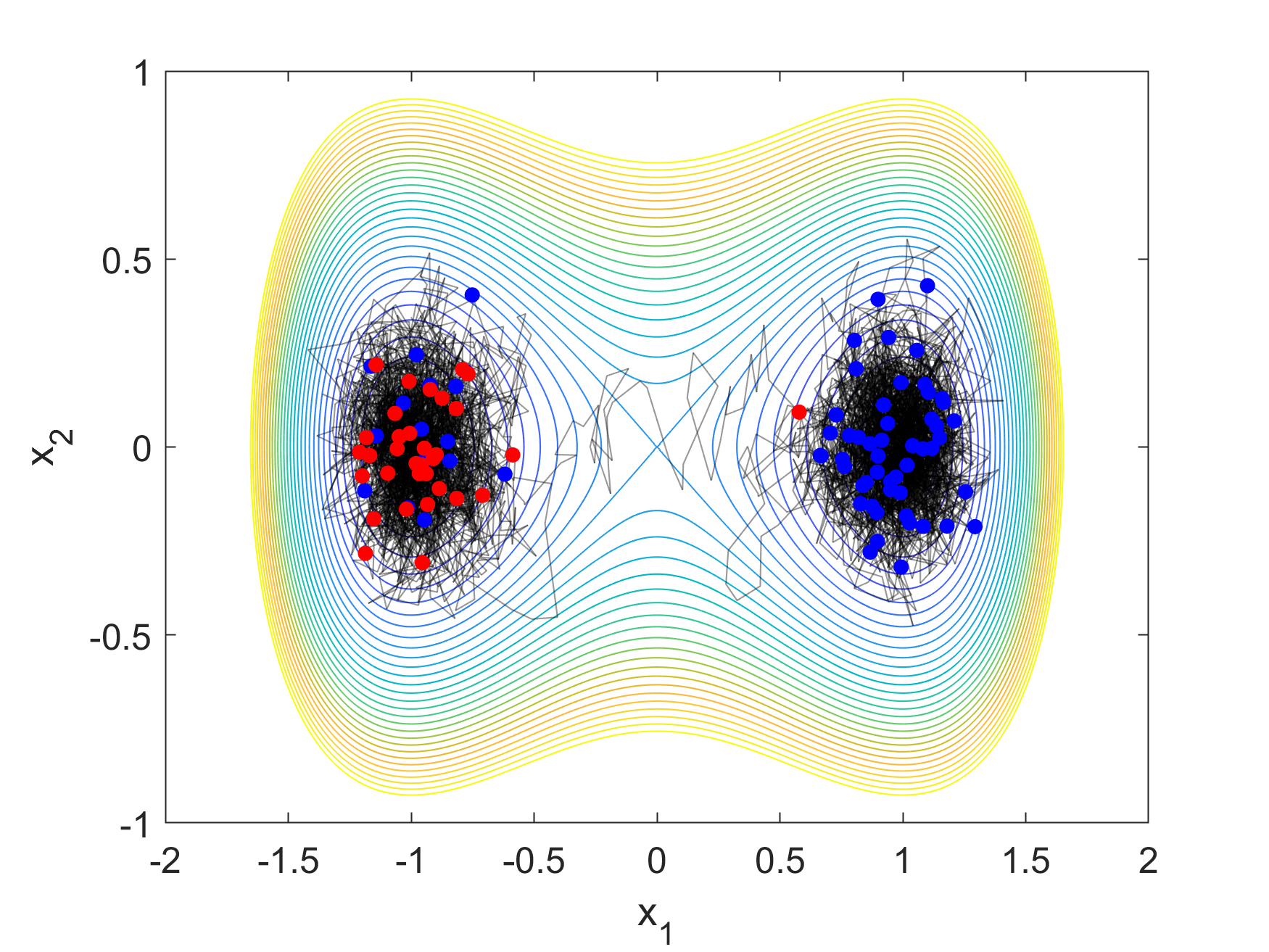}
        \caption{ABM Trajectory (snapshot)}
        \label{fig:pot}
    \end{subfigure}%
    \begin{subfigure}{0.5\textwidth}
    \centering
        \includegraphics[width=\linewidth]{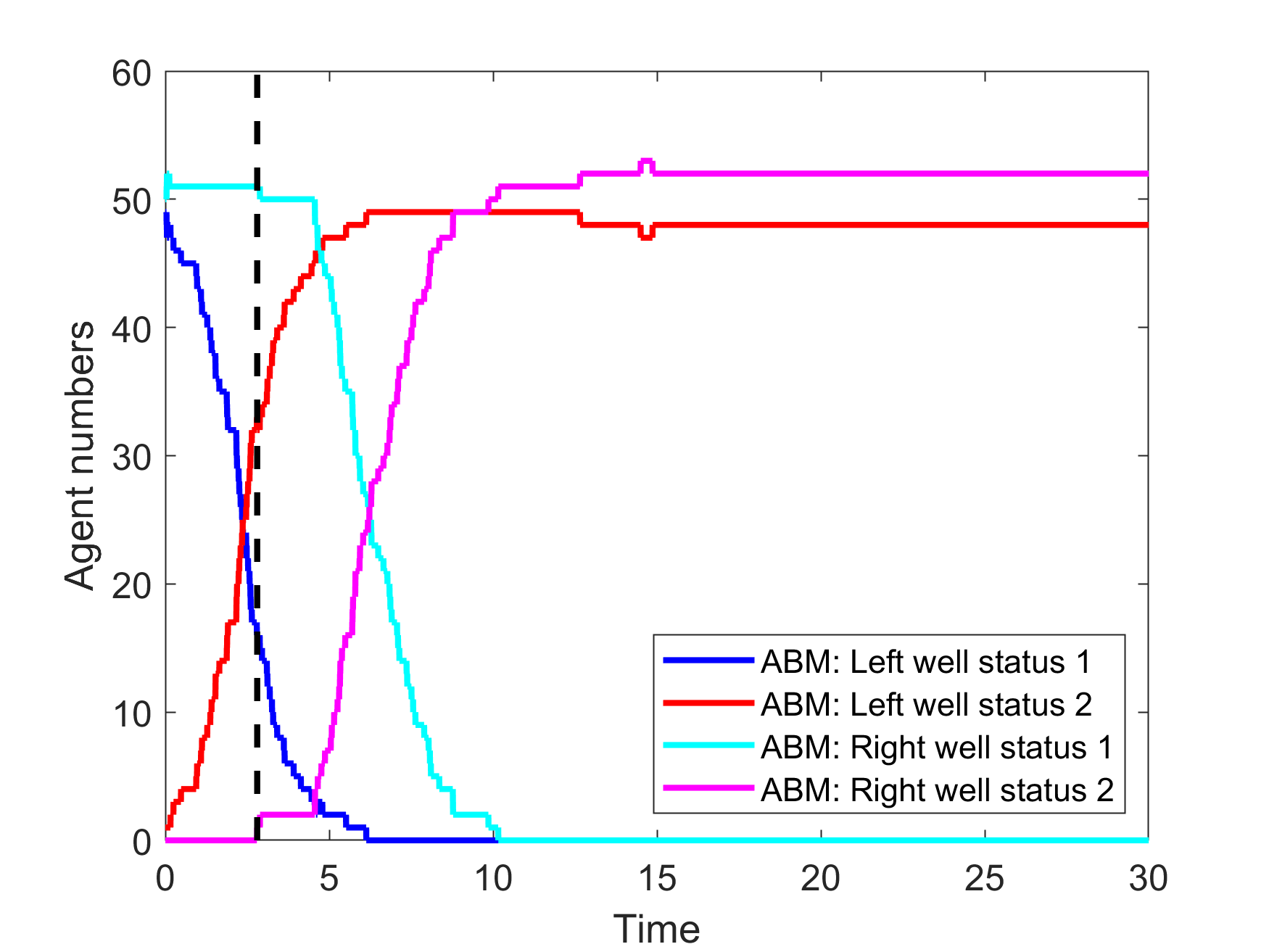}
        \caption{Projected ABM Trajectory}
        \label{fig:TrajMicro}
    \end{subfigure}
    \caption{Two-status dynamics in a double-well potential, see Example \ref{ex:guiding}. The adoption rate constant is chosen as $\revision{c}_{12}=0.1$ and the interaction radius as $r=0.15$ (a) Snapshot of an ABM trajectory for a simulation with $n_a=100$ agents with contour lines of the potential $U$ from Example \ref{ex:guiding} with diffusion constant $\sigma=1.2$. The blue dots refer to spatial positions of agents in status 1, the red dots represent agents in status 2. Grey lines show movement trajectories of the agents. (b) Projected trajectory of the same example, showing the temporal evolution of the total number of agents in each of the two wells (given by $(-\infty,0]\times \mathbb{R}$ and $(0,\infty)\times \mathbb{R}$) depending on the status. \revision{The dashed black line refers to the time of the snapshot.}
    }
    \label{fig:dw1}
\end{figure}

To illustrate the model reduction steps of the next section we introduce now the following guiding Example \ref{ex:guiding} for our ABM dynamics.

\begin{expl}[Two-status dynamics in a double-well potential]\label{ex:guiding}
We consider the continuous space for the agent movements to be $\X=\mathbb{R}^2$ and the discrete status space to be $\mathbb{S}=\{1,2\}$. The movement of a single agent is defined as a diffusion process and is given by the following stochastic differential equation 
\begin{equation*}
    d\boldsymbol{x}(t)=-\left(\frac{\sigma}{2}\right)^2\nabla U(\boldsymbol{x}(t))dt+\sigma d\boldsymbol{B}(t),
\end{equation*}
with $U(x_1,x_2)=(x_1^2-1)^2+7x_2^2$ (see Figure \ref{fig:dw1}),  diffusion constant $\sigma\in\R$ and $\boldsymbol{B}(t)$ a standard Brownian motion process in $\mathbb{R}^2$. The movement process of all $n_a=100$ agents is denoted by $\boldsymbol{X}(t)$. As for the status dynamics, we consider second-order adoptions with rate functions defined by \eqref{f_second-order} and \eqref{Doi}, assuming that only the changes into status $2$  are possible, namely with  rate constant $\revision{c}_{12}>0$, while transitions back to status $1$ are excluded by setting $\revision{c}_{21}=0$.
In the  initial state, all agents are in status 1 except for one agent in the left well, given by the subset $(-\infty,0]\times \mathbb{R} \subset \X$, which has status 2. The critical transition event that we are interested in is the first time that an agent with status 2 makes the transition from the left to the right well $(0,\infty]\times \mathbb{R}$. 
In Figure \ref{fig:TrajMicro} we see that this transition happens only after almost all agents in the left well have adopted status 2. This indicates that adoptions within one well happen faster than between two wells, which is due to the metastability of the diffusion process. In the following sections we will focus on this type of dynamics and derive reduced models that preserve its main properties. In particular, we will consider the critical transition event in order to compare the relevant statistics and quantify the approximation error of derived reduced models. \hfill $\diamondsuit$ 
\end{expl}

\subsection{Stochastic metapopulation model} \label{sec:MetaPop}
The diffusion process from Example \ref{ex:guiding} exhibits  metastable behavior: an agent remains for a comparatively long period of time within one well of the potential $U$ before it eventually jumps to the other well. For such metastable dynamics we consider a coarse-graining such that on the longest time-scale the dynamical properties of the original process are well approximated. \revision{More precisely, we assume that there exists a partition of the spatial domain $\X=\bigcup_{k=1}^m A_k$ into $m$ disjoint metastable sets $A_1 ,..., A_m\subset \X$, i.e., we assume rare transitions between the sets and comparatively fast mixing inside the sets.} Due to the fast mixing assumption we do not need to distinguish between the agents positions further than being in one of the metastable sets. Following these ideas, we will introduce here a stochastic metapopulation process, where the considered metastable sets (of the original dynamics) will represent the subpopulations and the rare transitions between the subpopulations will be reduced to jump dynamics. The dynamics within the subpopulations are affected by these transition events and by the adoption dynamics within subpopulations themselves. For more details, see \cite{winkelmann2016spatiotemporal} where such a coarsening has been proposed for chemical reaction-diffusion systems. 
Unlike in \cite{winkelmann2016spatiotemporal} where we focused on reactive particles, we will here consider the whole metapopulation state including several different statuses.
In the following, we will first describe the coarse-graining approach and then introduce details of the resulting stochastic metapopulation model.

We denote the system state of the stochastic metapopulation model at time $t$ by an \revision{$n_s\times m$} matrix $\boldsymbol{N}(t)=(\boldsymbol{N}_{i}^{(k)}(t))_{i=1,...,n_s,k=1,...,m}$, where $\boldsymbol{N}_{i}^{(k)}(t)$ refers to the number of members (agents) of subpopulation $k$ in status $i$ at time $t$. The set of all possible system states, given the total number $n_a$ of agents, is denoted by $\mathbb{M}_{n_a}$, 
\begin{equation}\label{def:M}
    \mathbb{M}_{n_a}:=\left\{N=(N_i^{(k)})_{i=1,...,n_s,k=1,...,m}\in\mathbb{N}_0^{n_s,m}: \sum_{i=1}^{n_s}\sum_{k=1}^{m}N_i^{(k)} = n_a\right\}.
\end{equation}
 A member of status $i$ in subpopulation $k$ transitioning to subpopulation $l$ at time $t$ leads to the immediate change $\boldsymbol{N}(t) \to \boldsymbol{N}(t)-E_{i}^{(k)}+E_{i}^{(l)}$, where $E_{i}^{(k)}$ is an \revision{$n_s\times m$} matrix with all entries zero except for the entry at index $(k,i)$ being one. Similarly, an adoption event from status $i$ to status $j$ in subpopulation $k$ leads to jump in the system state of the form $\boldsymbol{N}(t)\to\boldsymbol{N}(t)-E_{i}^{(k)}+E_{j}^{(k)}$. 

Let $P(N,t):=\prob(\boldsymbol{N}(t)=N)$ denote the probability to find the system in state $N\in \mathbb{M}_{n_a}$ at time $t$. In analogy to Eq. \eqref{eq:ABM} of the ABM we consider the equation  \revision{
\begin{equation}\label{SMM_equation}
 \frac{d}{dt} P(N,t) = \mathcal{L} P(N,t) + \mathcal{G} P(N,t)  
\end{equation}
for operators $\mathcal{L}, \mathcal{G}$ given by 
\begin{equation}\label{discreteOperators}
    \mathcal{L} P(M) := \sum_{N\in \mathbb{M}_{n_a}}\hcL_{NM}\cdot P(N), \quad  \mathcal{G} P(M) := \sum_{N\in \mathbb{M}_{n_a}}\hcG_{NM}\cdot P(N),
\end{equation}
}
where $\hcL_{NM}$ for $M\neq N$ is the rate to go from $N$ to $M$ by a spatial transition event between the subsets, 
while $\hcL_{NN}\revision{:=}-\sum_M\hcL_{NM} $. Analogously,  $\hcG_{NM}$ for $M\neq N$ is the rate to go from $N$ to $M$ by an adoption event, while $\hcG_{NN}=-\sum_M\hcG_{NM} $. 
In the following we specify the shape of $\hcL $ and $\hcG$ and show their connection to the ABM operators $L$ and $G$ by means of Galerkin projection methods.

\subsubsection*{Going from ABM to stochastic metapopulation}

Given the assumptions of well-mixed behavior, we will make use of the approach given in \cite{winkelmann2016spatiotemporal} and construct a Markov state model \revision{with respect to} the movements of the agents  by applying a Galerkin projection $Q$ (derived from the metastable partition) to the generator $L+G$ defined in \eqref{eq:ABM}, see also \cite{sarich2010approximation}. In particular, we will study how this projection acts on the transition rate function and on the rate functions for both first-order and second-order adoption events. For deriving the analytical results, we will focus on the case of a full partition of the state space, \revision{i.e., it holds $\X=\bigcup_{k=1}^m A_k$ for disjoint sets $A_k$}.

For any $N=(N_i^{(k)})\in \mathbb{M}_{n_a}$ we define the \revision{indicator} ansatz functions
\begin{equation}\label{def_Phi1}
    \Phi_N(X,S):=\revision{\prod_{k=1}^m\prod_{i=1}^{n_s}} \phi_{N_i^{(k)}}(X,S)
\end{equation}
with
\begin{equation}\label{def_Phi}
   \phi_{N_i^{(k)}}(X,S):= \delta_{N_i^{(k)}}\left(\revision{\sum_{\alpha=1}^{n_a}} \delta_{A_k}(x_\alpha)\delta_i(s_\alpha)\right), 
\end{equation}
where $\delta$ denotes Kronecker delta and set based indicator functions, as well. That is, $\Phi_N(X,S)$ has the value $1$ whenever there are, for each $i,k$, exactly $N_i^{(k)}$ agents $\alpha$ with position $x_\alpha \in A_k$ and status $s_\alpha=i$ , otherwise it is zero. \revision{These ansatz functions are thus non-negative and  fulfill $\sum_{N\in \mathbb{M}_{n_a}}\Phi_N(X,S)=1$ for all $(X,S)\in \mathbb{Y}$, thus they form a \textit{partition of unity} \cite{schutte2013metastability}. }

Next, we define the inner product of two functions $f,g: \mathbb{Y} \to \mathbb{R}$ as
\[
\langle f,g\rangle:=\frac{1}{(\mu(\X) n_s)^{n_a}}  \sum_{S\in \mathbb{S}^{n_a}}\int_{\X^{n_a}} f(X,S)g(X,S) \,dX,
\]
where $\mu$ denotes the Lebesgue measure. \revision{For the indicator ansatz functions $\Phi_N$ defined above we observe that it holds $\langle \Phi_M ,\Phi_N\rangle =0$ for $M\neq N$, while for $M=N$ we have $\langle \Phi_M,\Phi_N\rangle =\langle \Phi_N,\Phi_N\rangle=\langle \Phi_N,\mathbbm{1}\rangle $, where $\mathbbm{1}$ denotes the constant $1$-function on $\mathbb{Y}$.}

\revision{By means of this inner product we can consider  the full-partition  projection $Q:L^2(\mathbb{Y})\rightarrow D$ to the ansatz space $D = \mbox{span}\{\Phi_N,N \in \mathbb{M}_{n_a}\}$ given by \cite{schutte2013metastability}
\begin{equation}\label{def:Galerkin3}
      Qv = \sum_{N\in \mathbb{M}_{n_a}} \frac{\langle \Phi_N,v\rangle}{\langle \Phi_N,\mathbbm{1}\rangle} \Phi_N.
\end{equation}
Given any linear operator $H:L^2(\mathbb{Y})\to L^2(\mathbb{Y})$, a  Galerkin projection with $Q$ yields the \textit{projected operator} $QHQ:L^2(\mathbb{Y}) \to D$. 
The goal is now to find the matrix representations $\hcL =(\hcL_{NM})_{N,M\in \mathbb{M}_{n_a}}$ and $\hcG =(\hcG_{NM})_{N,M\in \mathbb{M}_{n_a}}$ (see Eq. \eqref{discreteOperators}) of the projected operators $QLQ$ and $QGQ$ for the operators $L$ and $G$ defined in \eqref{def:L} and \eqref{def:G}, respectively. }

\quad 

At first, we consider the spatial dynamics. Define 
 \begin{equation}\label{lambda}
    \lambda_i^{(kl)} := \frac{\langle \delta_{A_l}, L_i \delta_{A_k}\rangle_\X }{\langle \delta_{A_k}, \mathbbm{1}
   \rangle_\X} = \frac{\int_{\revision{\X}} \delta_{A_l}(x) (L_i \delta_{A_k})(x) dx}{\int_{\revision{\X}}{\delta_{A_k}(x)}  dx }
 \end{equation}
where $\langle \cdot , \cdot \rangle_\X$ refers to the standard scalar product for functions in $L^2(\mathbb{X})$ \revision{and $\mathbbm{1}$ denotes the constant $1$-function on $\X$} \cite{winkelmann2016spatiotemporal}.

\quad 

\begin{Thm}\label{Thm0}
The matrix representation of the projected generator $QLQ $ is given by $\hcL$ with
\begin{eqnarray*}
 \hcL_{NM} 
 & = & \left\{\begin{array}{ll}
  \lambda_i^{(kl)} N_i^{(k)}, & \mathrm{if} \;M=N+E_i^{(l)}-E_i^{(k)}\;,  k\neq l,\\
  -\sum_{i=1}^{n_s} \sum_{\substack{k,l=1;}{ l\neq k}}^m \lambda_i^{(kl)}N_i^{(k)},   & \mathrm{if} \;M=N, \\
0, & \mathrm{otherwise.}
 \end{array}
 \right.
\end{eqnarray*}

\end{Thm}

\noindent The proof can be found in the Appendix \ref{Appendix}. 

\quad

Now, we derive the matrix representations of the projected generator $QGQ$ for the interaction dynamics. Here, we consider the two fundamental cases of first- and second-order adoptions separately. 

For first-order interactions with adoption rate functions given in \eqref{f_first-order}, we define  the conditional expectation of $\rate_{ij}(x)$ given that $x\in A_k$:
\begin{equation}\label{gamma_k}
    \rate_{ij}^{(k)} :=\frac{\langle \rate_{ij},\delta_{A_k} \rangle_\X}{\langle \delta_{A_k}, \mathbbm{1}
\rangle_\X}=\frac{\int_{\revision{\X}} \rate_{ij}(x)\delta_{A_k}(x) dx}{\int_{\revision{\X}} \delta_{A_k}(x) dx}.
\end{equation}
Then, we obtain the following result.

\begin{Thm}\label{Thm1}
For first-order adoptions with an ABM rate function $f_{ij}^{(\alpha)}$ of the form \eqref{f_first-order}, the projected generator $QGQ$ has the matrix representation  $\hcG$ with
\[
\hcG_{NM} =\left\{\begin{array}{ll}
 \hat{\propensity}_{ij}^{(k)}(N),\quad & \mathrm{if} \; M = N+ E_{j}^{(k)} - E_{i}^{(k)} \; , i\neq j,\\
-\sum_{i,j=1}^{n_s}\sum_{k=1}^m \hat{\propensity}_{ij}^{(k)}(N),& \mathrm{if} \;M=N,\\
0, & \mathrm{otherwise},
\end{array}
\right.
\]
where 
\begin{equation*}
    \hat{\propensity}_{ij}^{(k)}(N):=\rate^{(k)}_{ij}N_i^{(k)}.
\end{equation*}
\end{Thm}

\noindent The proof can be found in the Appendix \ref{Appendix}. 

\quad

For the second-order status-changes, we consider the propensity function $f_{ij}^{(\alpha)}$ defined in \eqref{f_second-order}. 
Assuming the interaction distance $r$ to be small (compared to the size of the sets $A_k$),  second-order adoptions will mainly take place between agents of the same subpopulation as these agents are located relatively close to each other in space. 
Nevertheless, near the boundaries of the spatial subsets also agents of different subpopulations can be close enough to interact with each other, even if the interaction distance $r>0$ is chosen to be small.  However, in a metastable system, the sojourn of agents near the boundaries becomes unlikely and the probability for cross-over interactions (between different subpopulations) approaches zero. 
More precisely, let 
\begin{equation}\label{b_kl}
    b_{kl} :=
    \frac{\int_{\revision{\X^2}} d_r(x_1,x_2)\delta_{A_k}(x_1)\delta_{A_l}(x_2) dx_1dx_2}{\int_{\revision{\X^2}} \delta_{A_k}(x_1)\delta_{A_l}(x_2) dx_1dx_2}
\end{equation}
\revision{for $d_r$ given in \eqref{dr}} denote the conditional probability for two agents to  
be close enough to each other to interact,
given that they are located in the sets $A_k$ and $A_l$, respectively. 
Given the number state $N$ of the SMM system, we define
\begin{equation}\label{eps}
    \epsilon_{ij}^{(k)}(N) := c_{ij} \sum_{\substack{l=1\\ l\neq k}}^{m} b_{kl}  N_i^{(k)} N_j^{(l)}
\end{equation}
as the equilibrium propensity for an adoption event $i \to j$ to take place in subpopulation $k$ by means of a cross-over interaction with agents from a different subpopulation $l\neq k$, and
\begin{equation}\label{hatgamma}
    \hat{\rate}_{ij}^{(k)} := c_{ij} b_{kk} 
\end{equation}
as the macroscopic rate constant for adoptions 
within subpopulation $k$. Using these definitions, we obtain the following result.

\begin{Thm}\label{Thm2}
 For second-order adoptions with an ABM rate function $f_{ij}^{(\alpha)}$ given by \eqref{f_second-order} and  \eqref{Doi}, the  projected generator $QGQ$  has the matrix representation $\hcG$ with
 \small{
\begin{equation*} 
 \hcG_{NM} =\left\{\begin{array}{ll}
\hat{\propensity}_{ij}^{(k)}(N)+\epsilon_{ij}^{(k)}(N),\quad & \mathrm{if} \; M = N+ E_{j}^{(k)} - E_{i}^{(k)} , \; i\neq j,\\
-\sum_{i,j=1}^{n_s}\sum_{k=1}^m \left( \hat{\propensity}_{ij}^{(k)}(N)+\epsilon_{ij}^{(k)}(N)\right), & \mathrm{if} \;M=N,\\
0, & \mathrm{otherwise}.
\end{array}
\right.   
\end{equation*}}
where  
\begin{equation}\label{2ndOrderProp}
   \hat{\propensity}_{ij}^{(k)}(N) :=  \hat{\rate}^{(k)}_{ij}N_i^{(k)} N_j^{(k)} .
\end{equation}
\end{Thm}

\noindent   The proof is given in the Appendix \ref{Appendix}.

\quad

Theorem \ref{Thm2} shows how the projected rate functions of second-order adoptions can be decomposed into: (1) a part coming from adoptions that take place between agents of the same subpopulation with propensities given by $\hat{f}_{ij}^{(k)}$; and (2) a part coming from adoptions that take place between agents of different subpopulations. As discussed above, in a metastable system and for a good choice of sets $A_1,\ldots,A_m$, the probability of cross-over interactions is small and thus the value of $\epsilon_{ij}^{(k)}(N)$ will be negligibly small; in fact, as we will see below in Example ~\ref{ex:guide_cont}, it often is orders of magnitude smaller than the discretization error of the Galerkin projection.

In the SMM we only take the first type (1) of interactions, i.e., we  assume that there are no (cross-over) adoptions taking place between agents of different subpopulations. Using the results from Theorems \ref{Thm0}-\ref{Thm2}, the 
\revision{SMM equation given by  \eqref{SMM_equation}  and \eqref{discreteOperators} can be written as the following spatio-temporal master equation:}
\begin{align}\label{SMM_MasterEq}\begin{split}
	\frac{dP(N,t)}{dt}=&-\sum_{\substack{k,l=1\\ k\neq l}}^{m}\sum_{i=1}^{n_s}\lambda_{i}^{\revision{(kl)}}N_{i}^{\revision{(k)}}P(N,t)\\
	&+\sum_{\substack{k,l=1\\ k\neq l}}^{m}\sum_{i=1}^{n_s}\lambda_{i}^{(kl)}(N_{i}^{(k)}+1)P(N+E_{i}^{(k)}-E_{i}^{(l)},t)\\
	&-\sum_{i,j=1}^{n_s}\sum_{k=1}^{m}\hat{\propensity}_{ij}^{(k)}(N)P(N,t)\\
	&+\sum_{i,j=1}^{n_s}\sum_{k=1}^{m}\hat{\propensity}_{ij}^{(k)}(N+E_{i}^{(k)}-E_{j}^{(k)})P(N+E_{i}^{(k)}-E_{j}^{(k)},t),
	\end{split}
\end{align}
where the first two terms on the right-hand side refer to the change caused by the exchange between subpopulations \revision{(given by the operator $\mathcal{L}$)}, while the other two lines are describing the change through status adoptions inside the subpopulations \revision{(referring to operator $\mathcal{G}$)}.\footnote{\revision{In the second line of \eqref{SMM_MasterEq}, we need the rate to go from $M:=N+E_i^{(k)}-E_i^{(l)}$ to the given $N$. By Theorem \ref{Thm0} we know that this rate is given by $\lambda_i^{(kl)} M_i^{(k)} =\lambda_i^{(kl)} (N_i^{(k)}+1)$.}}

Using the definitions of the functions $\hat{f}_{ij}^{(k)}$ given in Theorems \ref{Thm1} and \ref{Thm2}, the interaction propensities agree with the standard law of mass-action from the chemical context \cite{bookC03}. This is due to the fact that we assume the agents to interact independently of each other (and of the overall system state)  -- which we do by choosing the ABM adoption rate functions according to Equations \eqref{f_first-order} and \eqref{f_second-order}. 
As for the spatial dynamics, 
state-of-the art metapopulation models  \cite{metapopulation2009} usually assume the commuting flow  between two subpopulations $k$  and $l$  to be  of the form \revision{ $\lambda_i^{(kl)} (N_k)^a (N_l)^b$}
for exponents $a,b\geq 0$ which tune the dependence
with respect to each subpopulation size.   
In our setting, we assume that the spatial movement of each agent is independent of the population sizes, which corresponds to setting $a=1$ and $b=0$.

\paragraph{Extension to core set approach}
The Galerkin projection does not have to be constructed on a full partition of the state space using indicator function as ansatz functions, but can also use ansatz functions that form a partition of unity. Appropriate ansatz functions of this form are defined via so-called \emph{core sets}, i.e., sets that do not form a partition of state space but
cover only the core areas of the metastable sets \cite{Djurdjevac2012,schutte2011markov,sarich2011projected}. For example, for a diffusive process in a potential energy landscape, the core sets are given by the vicinities around the energy's local minima (i.e., the valleys or wells of the landscape), while the transition regions around the energy's local maxima are not explicitly assigned to any core set. 
The advantage of considering a core set approach is a reduced approximation error in the estimated transition rates compared to the full partition case \cite{sarich2011projected,Djurdjevac2012}. The ansatz functions associated with the core set $C_k$ is its \textit{committor function} \cite{MetznerSchuetteEijnden09} $q_k$, with $q_k(x)$ denoting the probability that the agent visited the core set $C_k$ last, conditional on being in \revision{position} $x$, see \cite{schutte2013metastability}. \revision{Thus, for agents $x\in C_k$ we have that $q_k(x)=1$ and $q_l(x) = 0, l\neq k$. Agents that are in the transition region can be affiliated to several core sets with different probabilities, but such that it holds $\sum_{k=1}^m q_k(x) = 1$. }

For the core set approach, we have to redefine the conditional probabilities $b_{kl}$ given in \eqref{b_kl} according to 
\begin{equation}\label{b_kl2}
    b_{kl} :=
    \frac{\int_{\revision{\mathbb{X}^2}} d_r(x_1,x_2)q_{C_k}(x_1)q_{C_l}(x_2) dx_1dx_2}{\int_{\revision{\mathbb{X}^2}} q_{C_k}(x_1)q_{C_l}(x_2) dx_1dx_2}.
\end{equation}
Since in the core set approach we choose the sets $C_1,\ldots,C_m$ to be only the core areas of metastability, the probability of cross-over interactions of agents from different subpopulations is extremely small and thus, the value of $\epsilon_{ij}^{(k)}(N)$ is even smaller than in the case of a full partition. 

\begin{rem}[Approximation quality]
The step from the full-scale ABM to the SMM (\ref{SMM_MasterEq}) involves two approximations: the discretization error originating from the Galerkin projection, and the error resulting from neglecting the cross-over interactions between different spatial domains. While latter error can easily be monitored by estimating the neglected cross-over rates, controlling the discretization error is more difficult. For the cases where the spatial movement exhibits metastable sets, estimation of the discretization error is possible, see  \cite{schutte2011markov} for the core set approach, but requires sufficient ABM simulation data.
\end{rem}

In order to illustrate how the stochastic metapopulation model can approximate the ABM, we now return to our guiding example. 

\quad

\begin{figure}[h!] 
\centering
  \begin{subfigure}{0.5\textwidth}
    \centering
        \includegraphics[width=\linewidth]{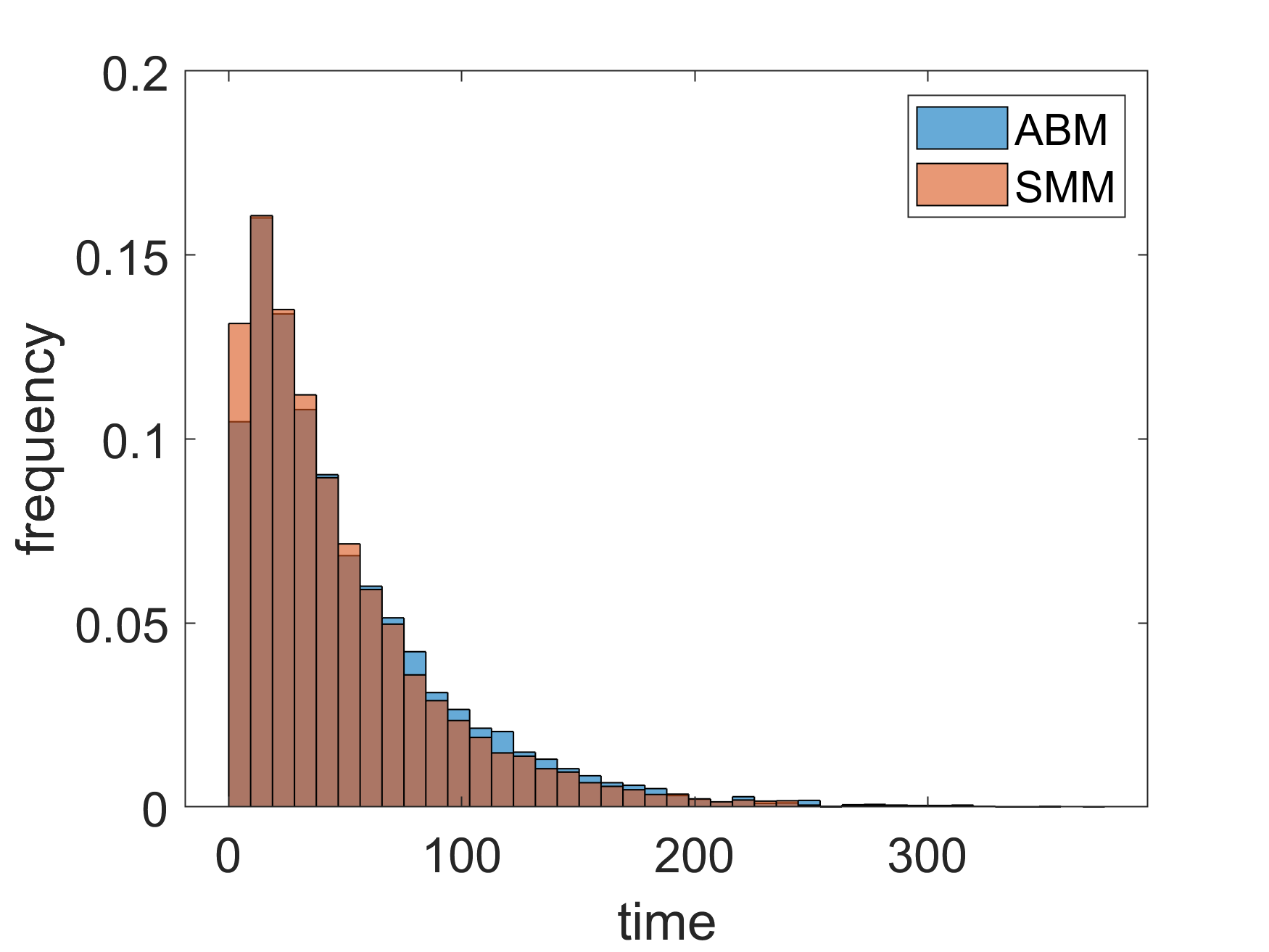}
        \caption{Critical transition time distribution \newline  for $\sigma=0.6$}
        
    \end{subfigure}%
    \begin{subfigure}{0.5\textwidth}
    \centering
        \includegraphics[width=\linewidth]{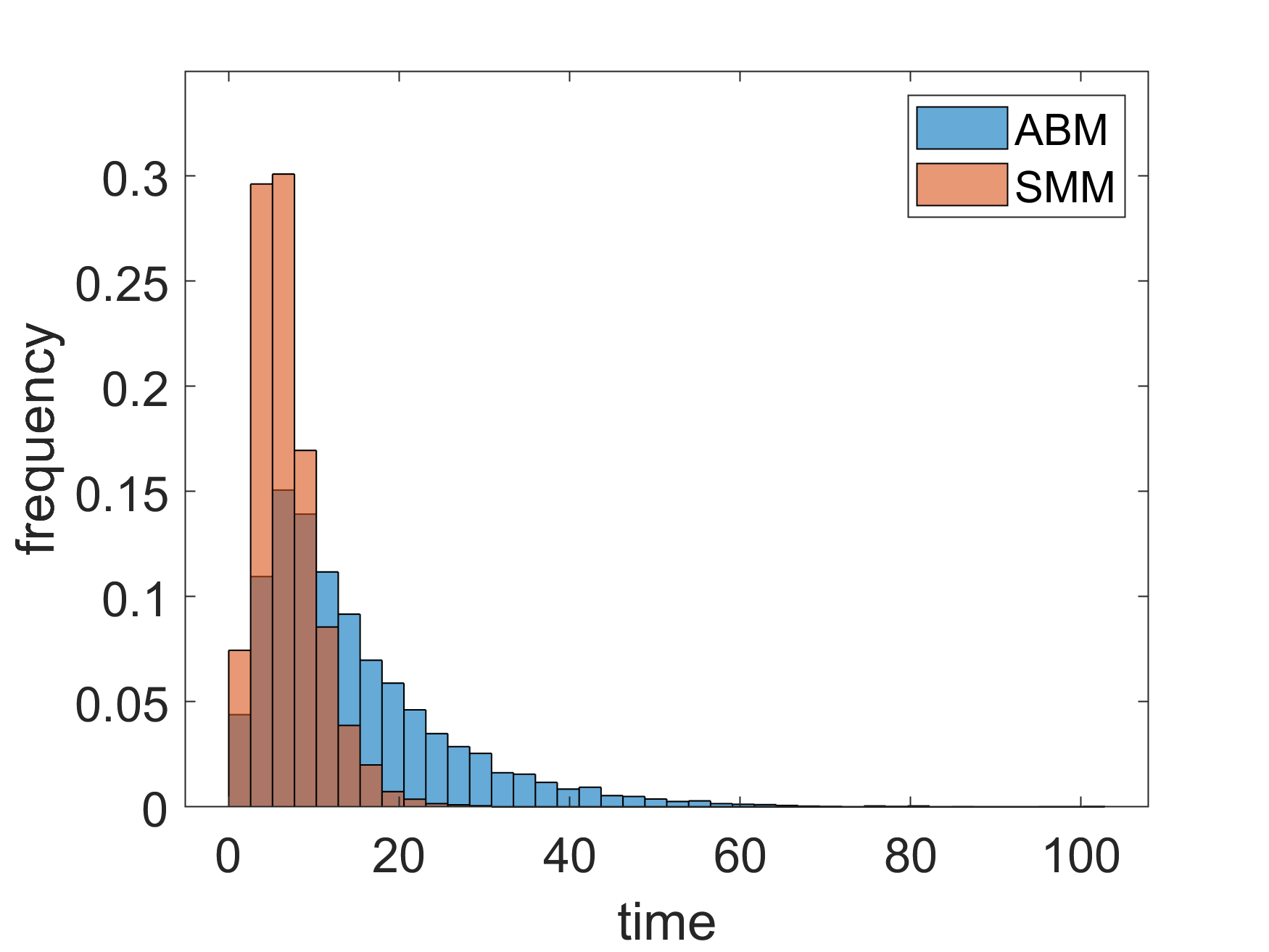}
        \caption{Critical transition time distribution \newline for $\sigma=1.2$}
        
    \end{subfigure}
   \caption{Distribution of the critical transition time for the ABM (blue) and the SMM (orange) given in Example \ref{ex:guiding}, sampled over $10000$ MC-simulations. \revision{The overlap of the two distributions is colored brown.} For small $\sigma$ the distribution is very well matched, while for larger $\sigma$ the critical transition happens faster in the SMM.}
   \label{fig:CritEvtDistMetaABM}
\end{figure}

\noindent \textbf{Example \ref{ex:guiding} (continued).}\label{ex:guide_cont}
Given the double-well potential shown in Figure \ref{fig:dw1} and using the \revision{Markov state model} approach, we partition the space $\X=\R^2$ of movement into the two core sets  $\revision{C}_1=(-\infty,-0.5)\times \mathbb{R}$ and $\revision{C}_2=(0.5,\infty)\times \mathbb{R}$ and the transition region $\X \setminus (\revision{C}_1 \cup \revision{C}_2)$. The jump rates $\lambda_i^{(12)}$ and $\lambda_i^{(21)}$ between the two subpopulations are 
the transition rates between $\revision{C}_1$ and $\revision{C}_2$. We analyze the quality of the SMM approximation for two different values of the diffusion constant $\sigma=0.6$ and $\sigma=1.2$, where the first case is more metastable than the other. We compare the SMM process to the projected ABM dynamics regarding the temporal distribution of a \textit{critical transition event} given by the first agent with status $2$ switching from one to the other subpopulation (meaning for the projected ABM that an agent of status $2$ who visited core set $\revision{C}_1$ last reaches core set $\revision{C}_2$ for the first time), see Figure \ref{fig:CritEvtDistMetaABM}. As this transition has a major impact on the overall dynamics, we consider the total approximation error to be small if the difference in the distributions of this critical transition event time is small. We observe that for smaller $\sigma$ the approximation is better, which is due to an increase in metastability of the dynamics. 

\begin{figure}[h!]
    \centering
    \begin{subfigure}{0.5\textwidth}
    \centering
        \includegraphics[width=\linewidth]{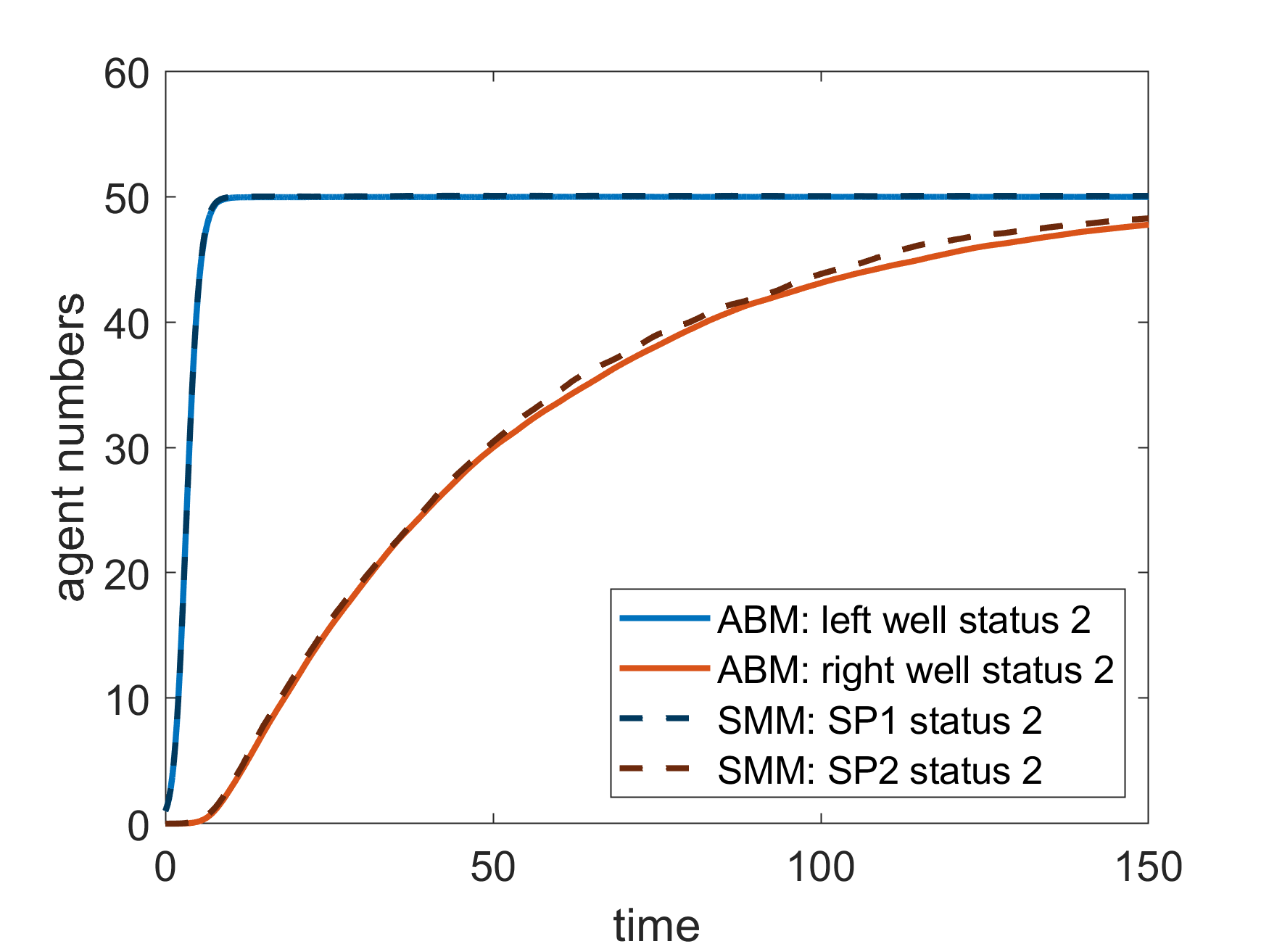}
        \caption{$\sigma=0.6$}
        \label{fig:abmcomp1}
    \end{subfigure}%
    \begin{subfigure}{0.5\textwidth}
    \centering
        \includegraphics[width=\linewidth]{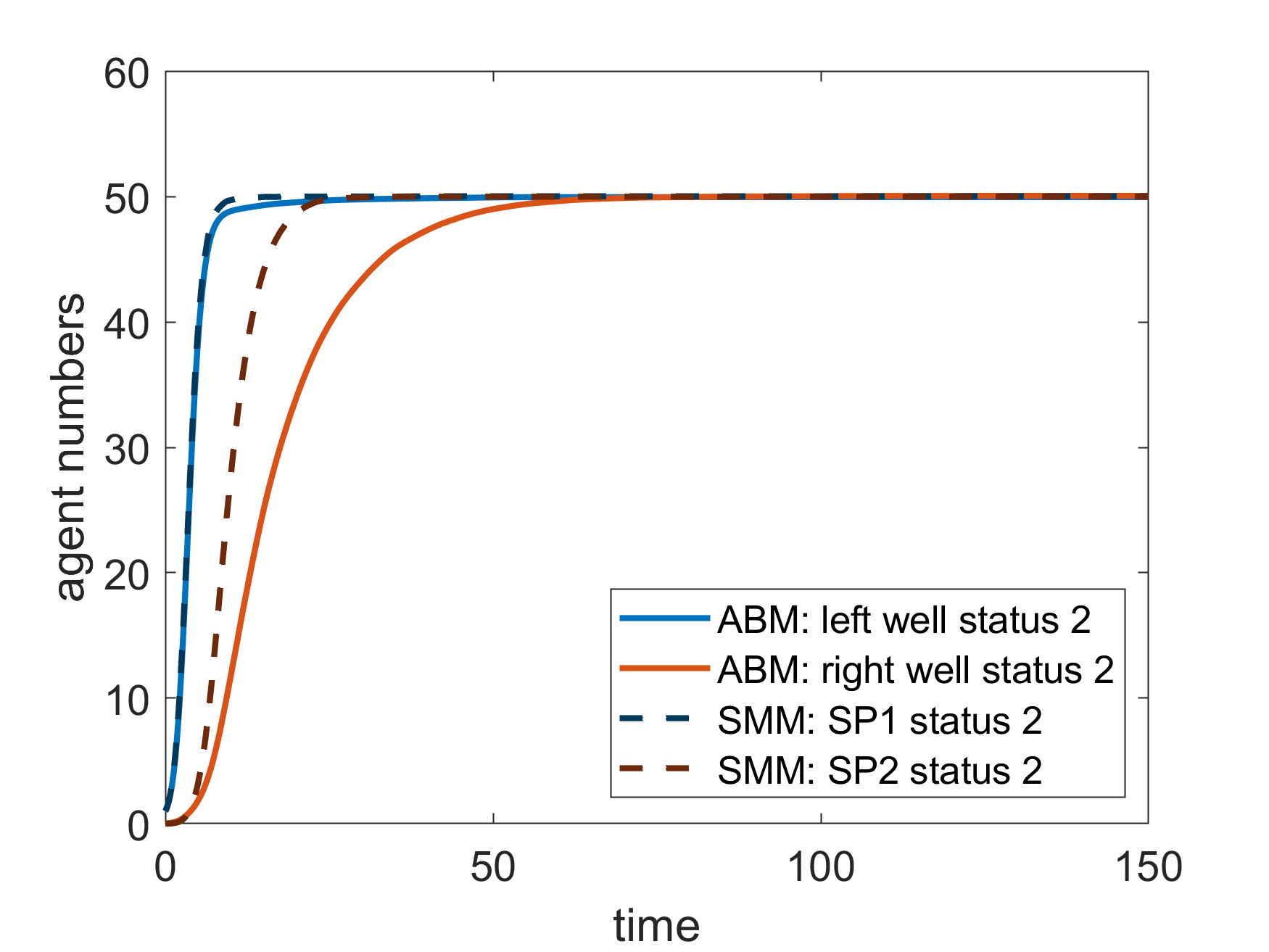}
        \caption{$\sigma=1.2$}
        \label{fig:abmcomp2}
    \end{subfigure}
  
    \caption{ Comparison of the time-dependent mean population size of agents with status $2$ over $10000$ MC-simulations for the projected ABM and the corresponding stochastic metapopulation model. The projection is based on a core set approach with $\revision{C}_1$ and $\revision{C}_2$ given  in Example \ref{ex:guiding} and defining two subpopulations denoted by SP1 and SP2. \revision{In scenario (a) the solid blue and dashed blue line are indistinguishable due to the good approximation quality.}}
    \label{fig:ABMcomp}
\end{figure}

This becomes even more clear when comparing  the temporal evolution of the average number of agents of status 2 in the two subpopulations, see Figure \ref{fig:ABMcomp}. For the smaller value $\sigma = 0.6$, 
the first-order moments agree very well (Figure \ref{fig:abmcomp1}). In contrast, there is  a significant difference between the model outcomes regarding these first-order moments for the larger diffusion constant $\sigma=1.2$ (Figure \ref{fig:abmcomp2}). This is due to the approximation quality of the Markov state model being worse because the diffusion process is less metastable and thus the first spatial transitions in the SMM happen too fast on average. The deviation of the spatial transition dynamics is also the main contribution to the approximation error of the critical transition event time.
\hfill $\diamondsuit$

\subsubsection*{Stochastic simulation}

For the stochastic metapopulation model we realize the sampling with the stochastic simulation algorithm \revision{, which produces statistically exact  realizations of the process, without any numerical approximation error \cite{gillespie1977exact}}. The effort is much lower than for ABM-simulations but still scales at least linearly with the number of agents. This means that for large population numbers the sampling of quantities of interest can still be infeasible due to high computational costs. This motivates \revision{us to consider} a further model approximation as discussed in the following subsection. 

\subsection{Piecewise-deterministic metapopulation model} \label{sec:PDMM}

For the case of a large population of interacting agents, the stochastic simulation of the stochastic metapopulation dynamics becomes computationally very expensive, because it tracks every single adoption event. \revision{Here, a further model reduction can be very useful. } 
Given that the number of agents is large in each subpopulation, we can apply standard convergence results and approximate the jump process which describes the internal adoption dynamics by a deterministic evolution equation \cite{kurtz1971limit}. \revision{Such approximations, which  are based on the law of large numbers, are well-known in the context of chemical reaction systems, where they are used to reduce the model complexity for systems with large molecular populations \cite{bookC03}. }
For the exchange process between the subpopulations, on the other hand, we assume that transitions are quite rare in time and occur sporadically at variable time points. In order to reduce the computational effort for simulations, while keeping the discrete, stochastic nature of the (inter-nodal) transition events, we approximate the overall dynamics by a piecewise-deterministic Markov process, see \cite{davis1984piecewise,alfonsi2005adaptive,franz2012piecewise,MONTAGNON2020} for details. This model will be called \textit{piecewise-deterministic metapopulation model (PDMM)}.

\subsubsection*{Going from stochastic metapopulation dynamics to PDMM}

The stochastic process $\NT$  given by \eqref{SMM_MasterEq}  can be rewritten in a pathwise notation of the form
\begin{align}\label{eq:Nt_pathwise} \begin{split}
    \boldsymbol{N}(t)=\boldsymbol{N}(0)&+\sum_{\substack{k,l=1\\ k\neq l}}^{m}\sum_{i=1}^{n_s}\mathcal{R}_i^{(kl)}\left(\int_0^t\lambda_{i}^{(kl)}\bN_i^{(k)}(s)ds\right)(E_{i}^{(l)}-E_{i}^{(k)})\\
    &+\sum_{i,j=1}^{n_s}\sum_{k=1}^{m}\mathcal{P}_{ij}^{(k)}\left(\int_0^t\hat{\propensity}_{ij}^{(k)}(\boldsymbol{N}(s))ds\right)(E_{j}^{(k)}-E_{i}^{(k)}) \end{split}
\end{align}
where $\mathcal{P}_{ij}^{(k)}$ and $\mathcal{R}_i^{(kl)}$ refer to independent, unit-rate Poisson processes \cite{alfonsi2005adaptive,menz2013hybrid}.

 Assuming that the jumps induced by the Poisson processes $\mathcal{P}_{ij}^{(k)}$, which refer to spatial transitions between metastable domains, occur with much less frequency than jumps induced by the Poisson processes $\mathcal{R}_i^{(kl)}$ referring to the adoption dynamics within the subpopulations, we can apply standard convergence results for Markov processes \cite{ethier2009markov} in order to approximate the stochastic dynamics given by the second line of  Eq. \eqref{eq:Nt_pathwise} by  deterministic dynamics and obtain the PDMM process $\hatnt$ given by the equation
\begin{align}\label{PDMM}\begin{split} 
    \hatn(t) =&\hatn(0)+ \sum_{k,l=1,k\neq l}^{m}\sum_{i=1}^{n_s}\mathcal{P}_i^{(kl)}\left(\int_0^t\lambda_{i}^{(kl)}\hatn_i^{(k)}(s)ds\right)(E_{i}^{(l)}-E_{i}^{(k)}) \\ & +\sum_{i,j=1}^{n_s}\sum_{k=1}^{m}\int_0^t\hat{\propensity}_{ij}^{(k)}(\hatn(s)) ds \, (E_{j}^{(k)}-E_{i}^{(k)}) .
    \end{split}
\end{align}
\revision{This means that the population dynamics on the local scale are modeled by a system of ODEs (adoption dynamics given in the second line), while the rare interactions between subpopulations are modeled as a stochastic jump process (first line), just as in \cite{MONTAGNON2020}. } 
It is well-known \cite{bookC03,kurtz1978strong} that the relative error produced by this approximation (relative with respect to the population size) decreases with an increasing number of agents with a Monte Carlo like rate. We will here consider a finite sized population which is large enough for the approximation to be reasonable.

\quad 

\begin{figure}[h!]
    \centering
    \begin{subfigure}{0.5\textwidth}
    \centering
        \includegraphics[width=\linewidth]{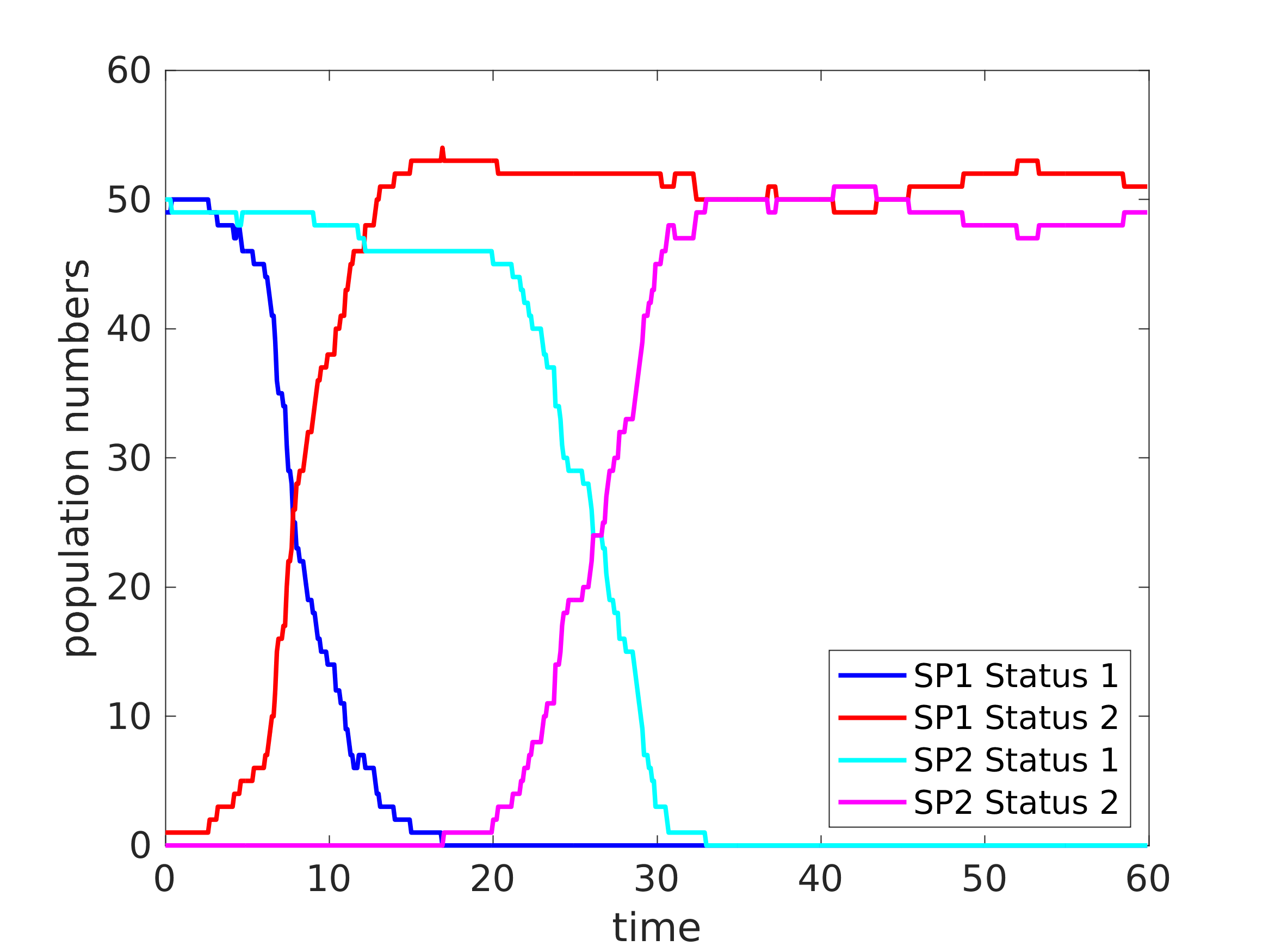}
        \caption{Stochastic metapopulation trajectory}
        \label{fig:mptr2}
    \end{subfigure}%
    \begin{subfigure}{0.5\textwidth}
    \centering
        \includegraphics[width=\linewidth]{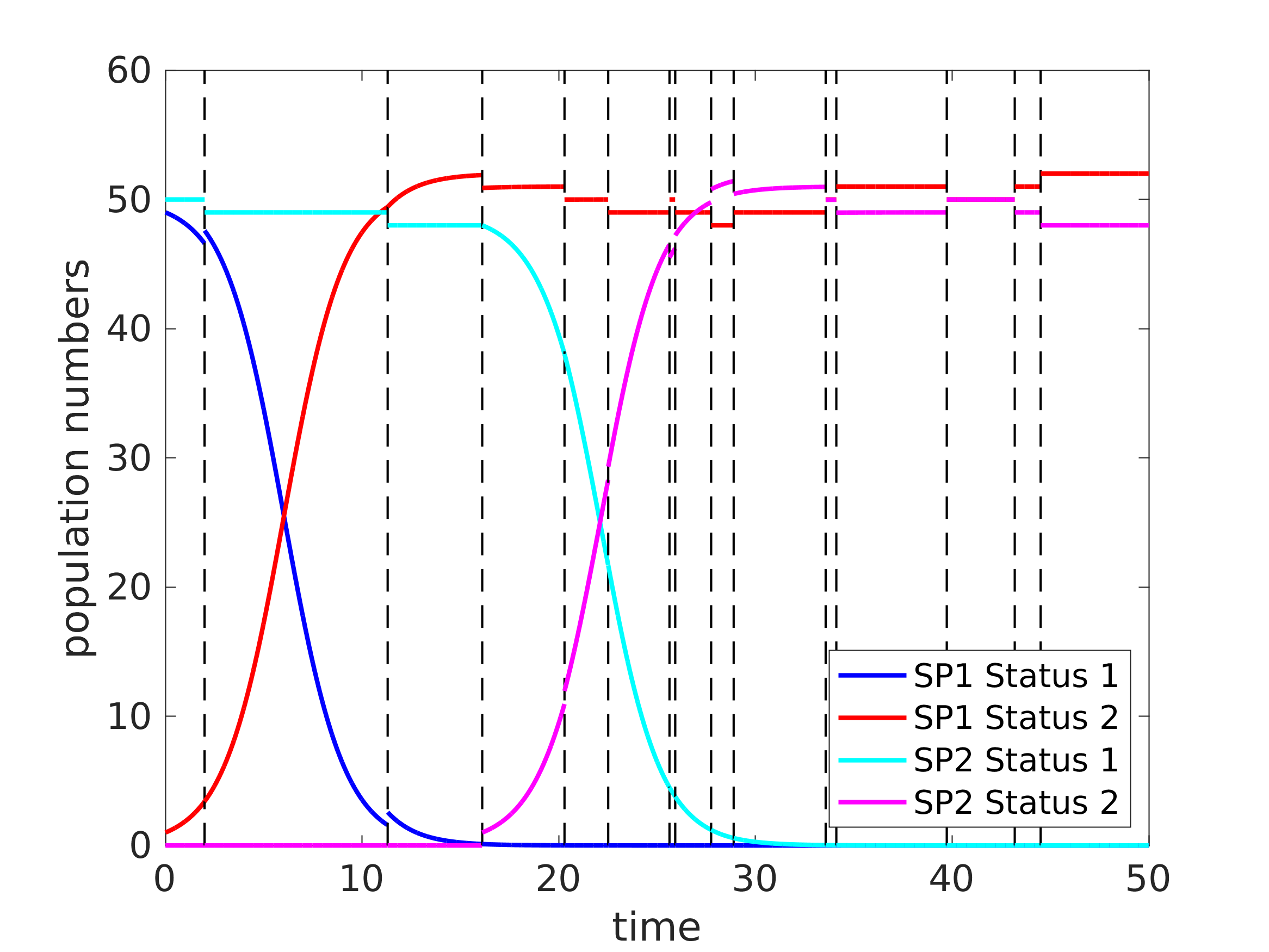}
        \caption{PDMM trajectory}
        \label{fig:PDMMtr}
    \end{subfigure}
    \caption{Comparison between (a) SMM and (b) PDMM trajectories for $\sigma=1.2$. In (b) the rare jump events are marked by vertical dotted lines.}
    \label{fig:TrajPDMM}
\end{figure}

\noindent \textbf{Example \ref{ex:guiding} (continued).}
For our example of two-status dynamics in a double-well potential, the deterministic status-adoption dynamics from the second line of \eqref{PDMM} are given by
\[    \hatn(t_0+\tau) = \hatn(t_0)+\sum_{k=1}^{m}\int_{t_0}^{t_0+\tau}\hat{\propensity}_{12}^{(k)}(\hatn(s))(E_{2}^{(k)}-E_{1}^{(k)}) ds \]
for $\tau < t_1-t_0$, 
where $t_0,t_1$ denote the time points of two subsequent stochastic transition events induced by the first line of \eqref{PDMM}. 
Using the definition \eqref{2ndOrderProp} of $\hat{\propensity}_{12}^{(k)}$, we get the following ODE for the number $\hatn^{(k)}_2$ of agents in subpopulation $k$ having status $2$:
\begin{equation} \label{ODE1}
 \frac{d\hatn^{(k)}_2(t)}{dt} = 
 \revision{\hat{\rate}_{12}^{(k)}}\cdot \hatn^{(k)}_1(t)\hatn^{(k)}_2(t) 
 \end{equation}
 for $t_0<t<t_1$.
Let $n^{(k)}_0:=\hatn_{1}^{(k)}(t_0)+\hatn_{2}^{(k)}(t_0)$ denote the total number of agents in subpopulation $k$ at time $t_0$. Between two transition events this number is constant, so we can substitute $\hatn_{1}^{(k)}(t)=n^{(k)}_0-\hatn_{2}^{(k)}(t)$ in Eq.  \eqref{ODE1} to arrive at 
$$
\frac{d\hatn_{2}^{(k)}(t)}{dt} =\revision{\hat{\rate}_{12}^{(k)}} \hatn_{2}^{(k)}(t)\left(n^{(k)}_0-\hatn_{2}^{(k)}(t)\right).
$$
The solution  is given by the logistic function, that is, we obtain an analytical solution
\begin{equation*}
        \hatn_2^{(k)}(t)=n^{(k)}_0\left(1+e^{-\revision{\hat{\rate}_{12}^{(k)}} n^{(k)}(t_0) t}\left(n^{(k)}_0-\hatn_{2}^{(k)}(t_0)\right)\right)^{-1}
\end{equation*}
for $t_0<t<t_1$. 
Treating the diffusive transitions between the subpopulations as stochastic events which induce jumps in the state $\hatn$ of the PDMM process, we obtain trajectories as depicted in Fig.  \ref{fig:TrajPDMM} (b). 

\begin{figure}[h!] 
\centering
  \begin{subfigure}{0.5\textwidth}
    \centering
        \includegraphics[width=\linewidth]{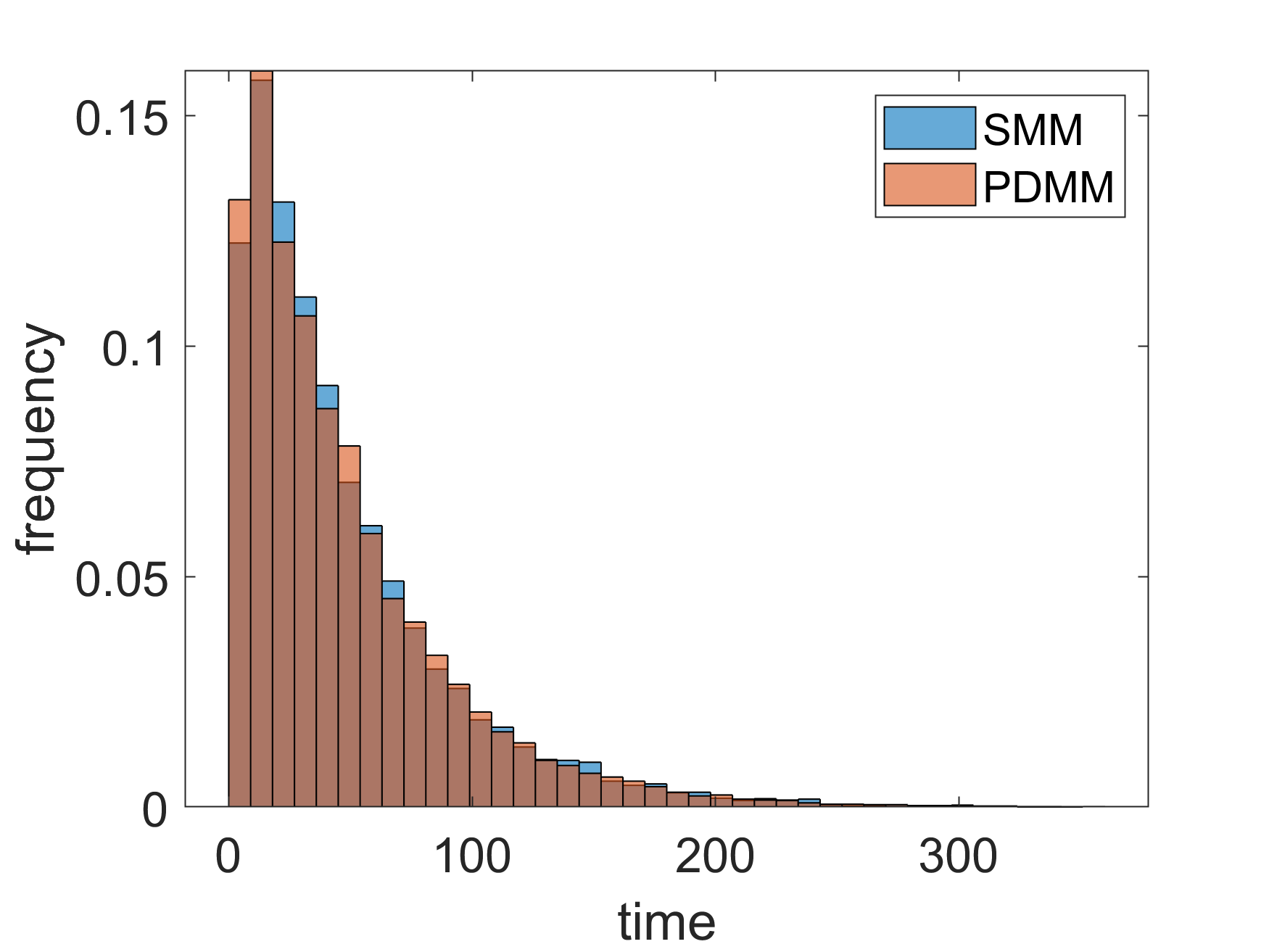}
        \caption{$\sigma=0.6$}
        
    \end{subfigure}%
    \begin{subfigure}{0.5\textwidth}
    \centering
        \includegraphics[width=\linewidth]{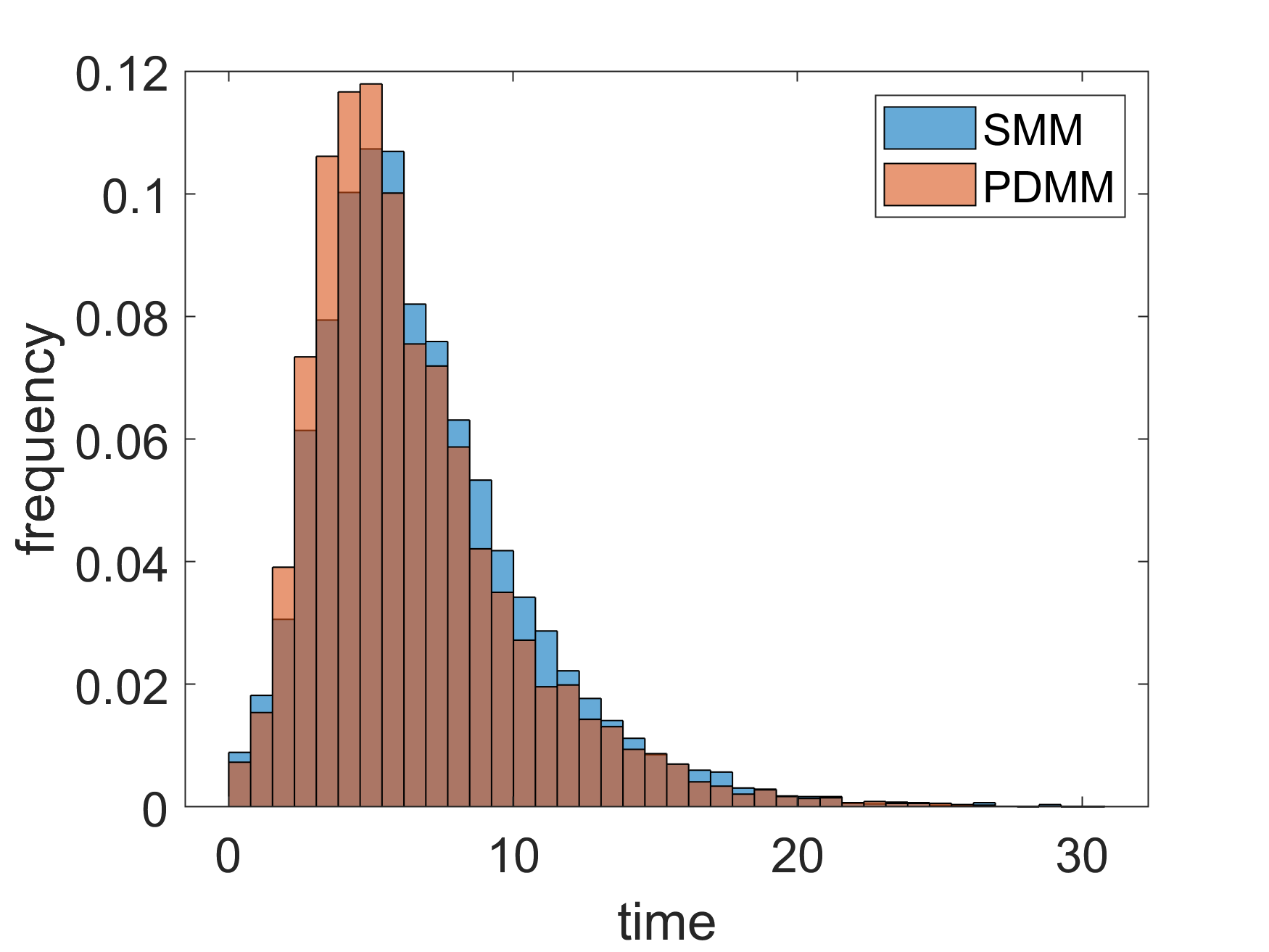}
        \caption{$\sigma=1.2$}
        
    \end{subfigure}
   \caption{Distribution of the critical transition time for SMM (blue) and PDMM (orange) sampled over 10000 MC-simulations. \revision{The overlap of the two distributions is colored brown.} For both values of $\sigma$ the distribution is very well matched. }
   \label{fig:CritEvtDistMetaPDMM}
\end{figure}

The main difference between the trajectories (see Figure \ref{fig:TrajPDMM}) is that in the stochastic metapopulation model we have many discontinuous jumps, while in the PDMM only a few remain.  The error for the estimation of the critical event time caused by the PDMM approximation (see Figure \ref{fig:CritEvtDistMetaPDMM}) is small compared to the error originating from the spatial discretization (see Figure \ref{fig:CritEvtDistMetaABM}). Even though our population number of 100 individuals is not very large, the critical transition time distribution of the SMM is already approximated well by the PDMM  for both choices of $\sigma$. \hfill $\diamondsuit$

\subsubsection*{Stochastic simulation}

\begin{figure}[h!]
    \centering
    \begin{subfigure}{0.45\textwidth}
    \centering
        \includegraphics[width=\linewidth]{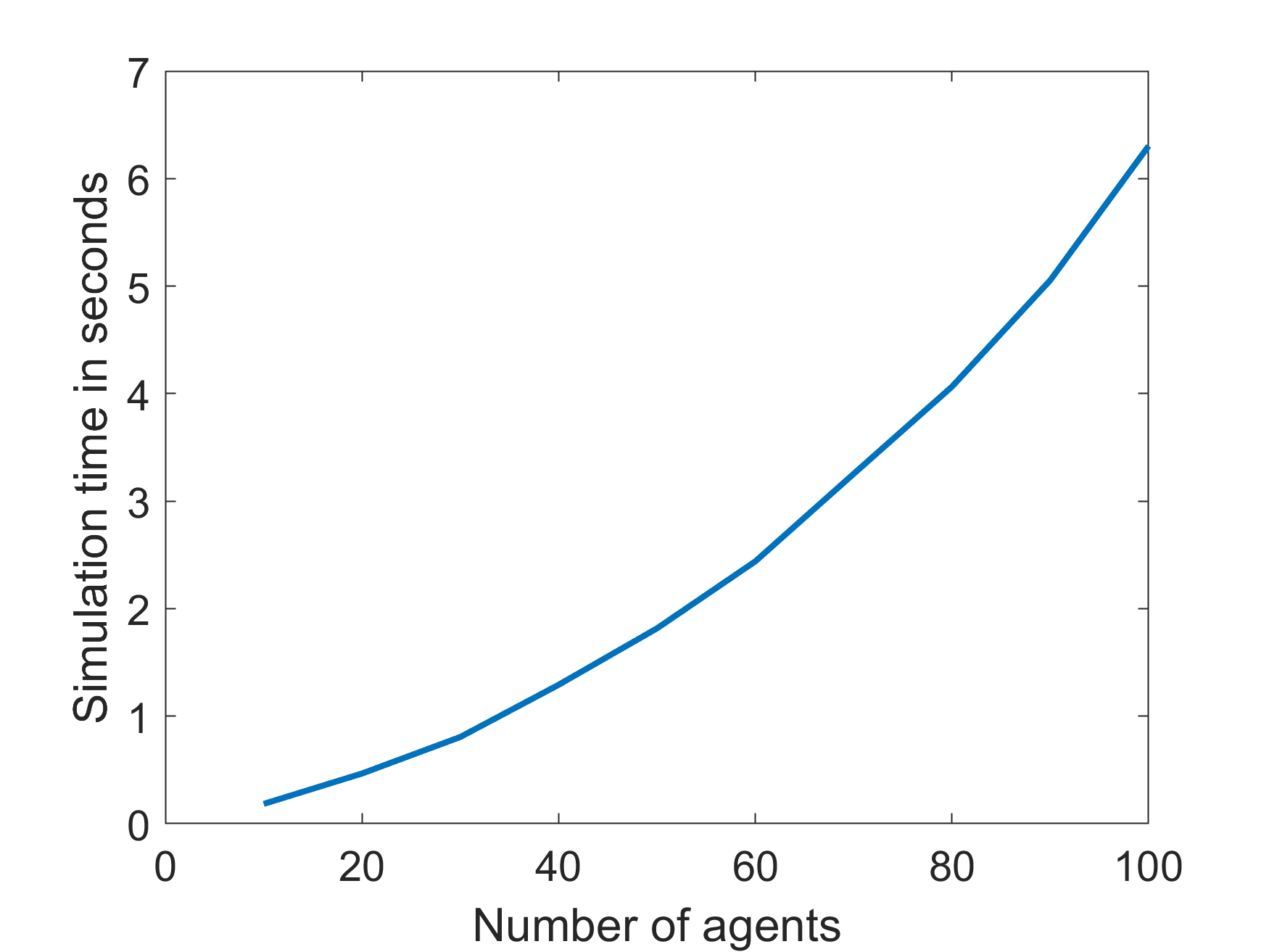}
        \caption{ABM effort}
        \label{fig:abmeff}
    \end{subfigure}
    \begin{subfigure}{0.45\textwidth}
    \centering
        \includegraphics[width=\linewidth]{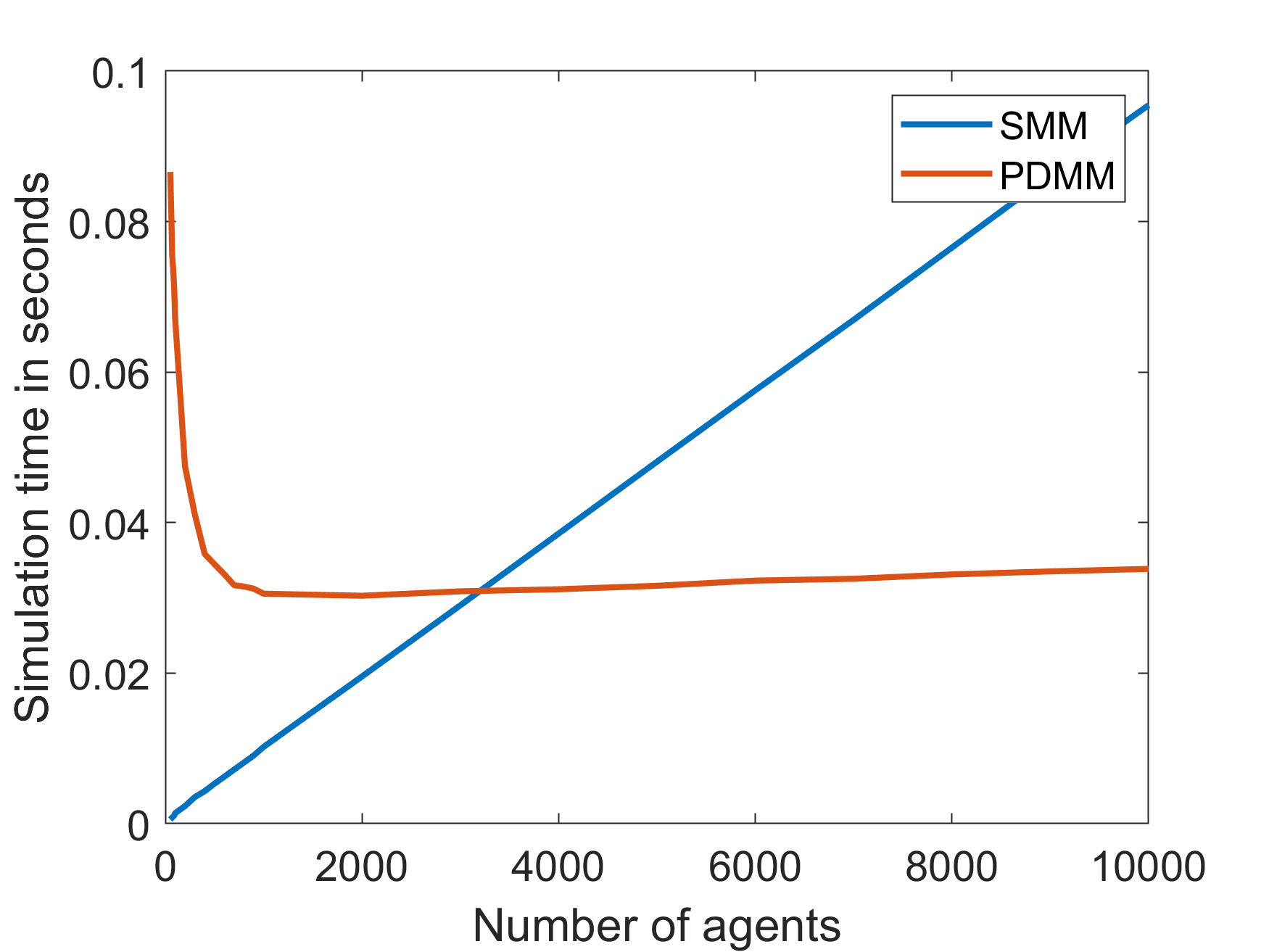}
        \caption{SMM and PDMM effort}
        \label{fig:discrteff}
    \end{subfigure}
    \caption{ Numerical effort for the simulations of the guiding example system of Ex. \ref{ex:guiding} for different choices of agent numbers $n_a$, depending on the number of agents. 
    Even for low agent numbers the approximate models are at least two orders of magnitude more efficient. \revision{For low agent numbers the PDMM effort is higher compared to the SMM as the critical transition event happens significantly later and thus more timesteps are needed in the PDMM computation.}}
    \label{fig:EffortCompare}
\end{figure}

In order to simulate a PDMM process given by \eqref{PDMM} one has to simultaneously integrate the deterministic flow of the ODE part  and the rate functions $\lambda_i^{(kl)}$ for the stochastic jumps, see \cite{menz2013hybrid}. Using the concept of the temporal Gillespie algorithm \cite{vestergaard2015temporal} one can thereby  determine the time point of the next stochastic  jump. After a jump event the initial conditions for the ODE are updated accordingly.\footnote{
\revision{More precisely, in our example simulations we compute the ODE solution with a standard forward Euler scheme. We adapt the time step size to restart the solver with new initial conditions once the parameters change by a jump event.}}
The numerical effort does not scale significantly with the number of agents, and thus, from the approaches considered in this paper, the PDMM is the most efficient choice for simulating systems with large agent numbers, see Figure \ref{fig:EffortCompare} for a comparison of the computational effort for the three modeling approaches. 

\section{Modeling COVID-19 epidemic spreading}
\label{sec:studies}
In this section, we discuss how 
the above presented approaches can be applied to the modeling of epidemic spreading. Agent-based models, as the most detailed models, have been extensively used for data-driven modeling with newly available large data-sets, e.g. on human mobility \cite{muller2020realistic,muller2020mobility}. However, most ABMs lack formalizations that would allow for more thorough mathematical analysis. Also, due to their large computational complexity, ABMs are usually focusing on epidemic spreading on smaller spatial scales, e.g. cities. Stochastic metapopulation models offer precise mathematical descriptions but, to the best of our knowledge, have not yet been investigated in the context of epidemic modeling. Piecewise-deterministic metapopulation models 
have been largely applied to model pathogen spreading \cite{abboud2018piecewise,soubeyrand2009spatiotemporal} due to  (1) their smaller computational complexity compared to ABMs; (2) better spatial resolution as opposed to ODE models, which are typically used. Despite the enormous recent advancement in adapting these models to realistic epidemic scenarios, this topic is still at the infancy of its development. In particular, mathematical model formalization with respect to data-based  descriptions on different levels needs to be explored in much more detail.\\
The results presented above show us how ABMs, SMMs and PDMMs are coupled. However, integrating all three types of models in one application would go beyond the scope of this article. Therefore, we will focus here on the most coarse-grained model, namely the piecewise-deterministic metapopulation model, and demonstrate how it can be applied to conceptually analyze COVID-19 spreading dynamics.

In the case of COVID-19, 
\revision{metapopulation models provide a good approximation of the original epidemic dynamics since the metastability assumption can be observed in mobility data \cite{SchlosserBrockmann2020}, e.g. in rare spatial transitions between different cities or countries caused by mobility restrictions. Also } 
large population sizes are realistic, such that an approximation by piecewise deterministic dynamics is justified and allows to drastically reduce the model complexity compared to the other two approaches. 
Moreover, the PDMM can more easily be calibrated to real world data than the underlying ABM, as we have good estimates for rates on the population scale, but not for individual interactions.

For simplicity, we will consider two subpopulations that have frequent local interactions and rare transitions in-between. 
The model will be calibrated to the parameters estimated in recent studies, for details see below. \revision{Note that in this paper, the model is not applied to analyze a particular real-world dataset, but representative results are used in order to show how the model can be applied to possible real-world scenarios.} Additionally, we will analyze the effectiveness of different containment measures taken within subpopulations and of global measures that impact the spreading between subpopulations. However, the main goal of this section is not to identify the optimal choice of containment measures, but to demonstrate the applicability of the PDMM on stochastic spreading processes and their efficiency on large real-world systems.

\subsection{The PDMM dynamics}
\label{subsec:LocalSystem}

As a first approach to formulate the piecewise-deterministic spreading dynamics, we will consider the rate constants $\rate$, which define the propensity of agents to change their status, to be independent of the evolution of the process. In Section \ref{sec:rate functions} we will generalize the dynamics by letting these constants depend on time and on the process' history. 

\subsubsection*{The system state}

On a local scale of each subpopulation, we use a compartmental  Susceptible-Exposed-Infected-Recovered-Deceased (SEIRD) model \cite{kermack1927contribution,weitz2015modeling} to describe the\\ COVID-19 dynamics. In this model, each individual can be in one of five possible statuses: susceptible (S), exposed (E), infected (I), recovered (R) and deceased (D). Susceptible individuals are the ones who have not yet been in contact with the virus and have no immunity against it. After being exposed to the virus, a susceptible individual is first in an asymptomatic status E that changes to a symptomatic status I after the end of an incubation period. Recovered and deceased individuals are considered to be immune for the time scale of our model and are no longer able to be infected and transmit the disease.
\begin{figure}
    \centering
    \includegraphics[scale=0.8]{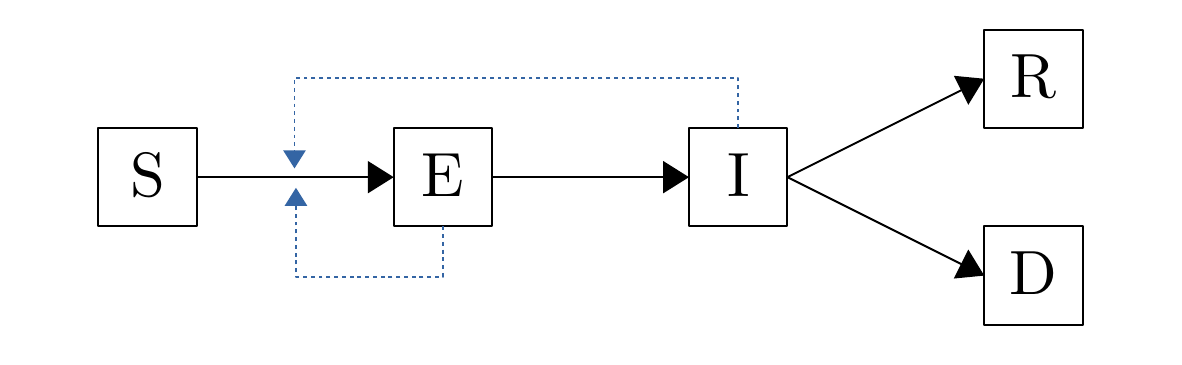}
    \caption{Visualization of the SEIRD model. Black arrows stand for possible status transitions, blue arrows indicate an impact by interaction.}
    \label{fig:seird}
\end{figure}

Given these statuses as well as a set of $m$ subpopulations, a possible state of the system has the form
\[ N = \left(N_S^{(k)},N_E^{(k)},N_I^{(k)},N_R^{(k)},N_D^{(k)}\right)_{k=1,...,m}\]
where $N^{(k)}_i$ denotes the number of individuals in status $i\in \{S,E,I,R,D\}$ within subpopulation $k$. 

\subsubsection*{Deterministic local interaction  dynamics}

We assume that exposed individuals are already able to transmit the virus \revision{\cite{SEIRKieran2020}}, so the status change from $S$ to $E$ can be caused by the second-order interactions of type either 
$S+E\rightarrow 2E$ or $S+I\rightarrow E+I$, with respective transition rate constants $\rate^{(k)}_{SE}> 0$ and  $\rate^{(k)}_{SI}>0$, see Figure \ref{fig:seird} for an illustration.  \revision{In general, it is possible to distinguish between infection rates from a contact with an exposed individual and from a contact with an infected individual. However, for simplicity, we assume here that infectiousness is constant from the moment of exposure until recovery or death, i.e. $\rate^{(k)}_{SE}=\rate^{(k)}_{SI}$.}\\
Given the system state $N$, the rate functions for the second-order status-change  
are given by 
\[\propensity_{SE}^{(k)}\left(N\right) = \rate^{(k)}_{SE} N_E^{(k)} N_S^{(k)}\]
and
\[\propensity_{SI}^{(k)}\left(N\right) = \rate^{(k)}_{SI} N_I^{(k)} N_S^{(k)}, \]
respectively, see \eqref{2ndOrderProp}.\footnote{ \revision{Note that in \eqref{2ndOrderProp} our notation includes $\hat{\cdot}$, e.g. $\hat{\propensity}_{ij}^{(k)}(N)$ which stands for a propensity in a reduced model. In this section we omit the notation with $\hat{\cdot}$ in order to make our calculation more clear.}} The remaining status transitions that we consider are given by first-order events of the form 
$E\rightarrow I, I\rightarrow R$ and $I\rightarrow D$ with respective rate constants $\rate^{(k)}_{EI},\rate^{(k)}_{IR},\rate^{(k)}_{ID}> 0$. For the transition $E\rightarrow I$ we accordingly obtain
\[ \propensity_{EI}^{(k)}(N) = \rate^{(k)}_{EI} \cdot N_E^{(k)},  \]
and equivalently for $I\rightarrow R$ and $I \rightarrow D$
\[ \propensity_{IR}^{(k)}(N) = \rate^{(k)}_{IR} \cdot N_I^{(k)},   
\quad  \propensity_{ID}^{(k)}(N) = \rate^{(k)}_{ID}\cdot N_I^{(k)}. \]
The resulting ODE-system describing the local interaction dynamics within a subpopulation $k$ is then given by 
\begin{equation}
\begin{aligned}\label{eq:SEIRD_constant}
    \frac{d}{dt}\boldsymbol{N}_S^{(k)}&=- \left( \rate^{(k)}_{SE} \boldsymbol{N}_E^{(k)}+ \rate^{(k)}_{SI}\boldsymbol{N}_I^{(k)}\right)\boldsymbol{N}_S^{(k)}\\
    \frac{d}{dt}\boldsymbol{N}_E^{(k)}&=
    \left( \rate^{(k)}_{SE} \boldsymbol{N}_E^{(k)}+ \rate^{(k)}_{SI}\boldsymbol{N}_I^{(k)}\right)\boldsymbol{N}_S^{(k)}
    -\rate^{(k)}_{EI}\boldsymbol{N}_E^{(k)}\\
    \frac{d}{dt}\boldsymbol{N}_I^{(k)}&=\rate^{(k)}_{EI}\boldsymbol{N}_E^{(k)}-\left(\rate^{(k)}_{IR}+\rate^{(k)}_{ID}\right)\boldsymbol{N}_I^{(k)}\\
    \frac{d}{dt}\boldsymbol{N}_R^{(k)}&=\rate^{(k)}_{IR} \boldsymbol{N}_I^{(k)}\\
    \frac{d}{dt}\boldsymbol{N}_D^{(k)}&=\rate^{(k)}_{ID} \boldsymbol{N}_I^{(k)}.
    \end{aligned}
\end{equation}

\subsubsection*{Stochastic dynamics for spatial exchange}

While the local interaction dynamics are given by deterministic evolution equations \eqref{eq:SEIRD_constant}, the spatial transitions between the subpopulations are described by sto\-chastic jump events. At any point in time $t\geq 0$, an agent of  subpopulation $k$ can switch  to another subpopulation $l\neq k$, which induces a discrete change in the system's state of the form
\[ \boldsymbol{N}(t) \mapsto \boldsymbol{N}(t) - E_i^{(k)} + E_i^{(l)} \]
depending on the agent's status  $i\in \{S,E,I,R\}$. (For $i=D$ we naturally assume that jumps cannot take place.) The PDMM process combines these discrete, stochastic jump events between the subpopulations with the ODE dynamics \eqref{eq:SEIRD_constant} for local status-changes, which 
in total, leads to a stochastic process $\boldsymbol{N}(t)_{t\geq 0}$,
\[ \boldsymbol{N}(t)=\left(\boldsymbol{N}_S^{(k)}(t),\boldsymbol{N}_E^{(k)}(t),\boldsymbol{N}_I^{(k)}(t),\boldsymbol{N}_R^{(k)}(t),\boldsymbol{N}_D^{(k)}(t)\right)_{k=1,...,m},  \]
described by an equation of the form \eqref{PDMM}. 
The terms in the first line of \eqref{PDMM} thereby correspond to an integrated version of the ODE \eqref{eq:SEIRD_constant} for each $k$, and the second line describes the spatial jumps for given rate constants $\lambda_{i}^{(kl)}$ between subpopulations $k$ and $l$. In the following subsection, both the rate constants $\rate$ for the local interactions and the rate constants $\lambda$ for the spatial transitions will depend on the  evolution of the overall stochastic process  $\boldsymbol{N}(t)_{t\geq 0}$.

\revision{\begin{rem} modeling choices presented above are made for simplicity and to demonstrate how the PDMM can be used to conceptually analyze COVID-19 spreading. Using rich data-sets and extensive literature on COVID-19 pandemic, our model can be easily extended to account for more realistic scenarios. For example, considering additional compartments such as symptomatic, asymptomatic, quarantined individuals \cite{Meyer_Hermann2021}; including more general infection rates with possible time dependency \cite{jo2020analysis}, adding demographics information \cite{bubar2021model}\cite{keeling2021predictions}, introducing vaccination effects \cite{bubar2021model} are only some of the extensions that would make this model more realistic. 
\end{rem}}

\subsection{Adaptive regulation of rate constants}
\label{sec:rate functions}

Until effective pharmaceutical treatment for \revision{COVID-19} is found, a lot of effort is taken to slow down the virus spreading by introducing measures for reducing social contacts. This is achieved, for example, by targeting the individual interactions (e.g. social distancing and wearing masks), by reducing the number of interactions (e.g. closing of schools, offices), but also by introducing measures on a global level, such as travel bans between countries and continents.\\
Having this in mind, the choice of transition rate constants $\rate$ that are independent of time and of the process' evolution (as assumed in Section \ref{subsec:LocalSystem}) appears to be unrealistic because containment measures are taken depending on the dynamics in order to influence the future evolution of the process. This is why we will in the following consider rate constants that are adapted in the course of time according to given rules.

The transition rates are depending on the local contacts within a population, which are changing over time with the installation of measures. Additionally, the dependence can also be on the interaction type, e.g. exposed individuals might be less contagious than infected ones or symptomatic cases could cause less infections than expected because they have already reduced their number of contacts. In order to include many types of possible dependencies, we define the transition rates in a quite general way as functions of the process' history $N_{\leq t}:=(N(s))_{s\leq t}$ and time $t$ 
\[ \rate^{(k)}_{SE}=\rate^{(k)}_{SE}\left(N_{\leq t},t\right), \quad  \rate^{(k)}_{SI}=\rate^{(k)}_{SI}\left(N_{\leq t},t\right).\] 
That is, these rates not only depend on the current state but also on the history of the process. By this we can define rules such as implementing a strict lock-down when case numbers are rising for the first time. In our model, the rate for developing symptoms $\rate_{EI}\geq 0$ is a constant, while the recovery and case fatality rates depend on the capacity of the health care system of a population, so they are defined as state-dependent rates $\rate^{(k)}_{IR}=\rate^{(k)}_{IR}(N)$ and $\rate^{(k)}_{ID}=\rate^{(k)}_{ID}(N)$.  

One goal of the introduced measures is to control the number of infections such that the limits of the health care system are not reached. As part of the global measures, the transitions between the subpopulations will be reduced. Thus, the spatial transition rates $\lambda_i^{(kl)}$ between subpopulations  are defined as functions of the entire process' history and time for $i\in \{S,E,I,R\}$:  \[ \lambda_i^{(kl)}=\lambda_i^{(kl)}\left(N_{\leq t},t\right) . \] 
\subsubsection*{Concrete choice of rate constants for status changes}
When modeling the implementation of measures for virus containment, we will assume that each measure is followed by a phase for which the infection rate remains constant. The transition between the phases can be triggered by deterministic as well as stochastic events, e.g. by the process crossing a threshold number of infections for the first time. In total, we will consider three different phases: 
\begin{enumerate}
\item \textbf{Initial phase:} In the beginning of the pandemic, the infection rates $\rate_{SE}^{(k)}$ and $\rate_{SI}^{(k)}$ have values $\delta_{SE}>0$ and $\delta_{SI}>0$, respectively, and the interaction dynamics start with an unmitigated spreading. 
\item \textbf{Strict measures phase:} The first measures to reduce the infection rates are taken in  subpopulation $k$ as soon as the number of infected individuals crosses a critical value $h_I^{(k)}>0$ for the first time. That is, the strict measures start at the random first hitting time
\[ t_1^{(k)}(N_{\leq t}) := \min\left\{ 0\leq s \leq t\left| \boldsymbol{N}_I^{(k)}(s) \geq h_I^{(k)}\right.\right\} \in [0,\infty],\]
with value $t_1^{(k)}(N_{\leq t}) = \infty$ in case of $N_I^{(k)}(s)< h_I^{(k)}$ for all $s\leq t$. In this phase, the infection rates $\rate_{SE}^{(k)}$ and $\rate_{SI}^{(k)}$ 
are reduced by a factor $\kappa_1^{(k)}\in (0,1)$. These strict measures are kept until the number of infected individuals falls below the critical value $\frac{h_I^{(k)}}{2}$, i.e. until the random time point 
\[ t_2^{(k)}(N_{\leq t}) := \min\left\{ t_1^{(k)}<s\leq t \left| \boldsymbol{N}_I^{(k)}(s) < \frac{h_I^{(k)}}{2}\right.\right\} \in [0,\infty].\] 

\item \textbf{Moderate measures phase:} After the strict measures are lifted, the interactions inside the population do not go back to normal, i.e. to the values from the initial phase. Instead, we introduce moderate measures where the infection rates are scaled with a factor $\kappa_2^{(k)}$, s.t. $\kappa_1^{(k)}<\kappa_2^{(k)}<1$, allowing for more contacts than in the previous phase. These measures are maintained for the remaining time of the model \revision{even if the number of infected individuals crosses again the value $h_I^{(k)}$}. 
\end{enumerate}

\quad

Taken all together, this means that the infection rate function is defined by 
\begin{equation}\label{eq:HdepInfect}
   \rate^{(k)}_{SE}\left(N_{\leq t},t\right):=
    \begin{cases}
    \delta_{SE} &\text{for } t\leq t_1^{(k)}\left(N_{\leq t}\right)\\
    \kappa_1^{(k)}\delta_{SE} &\text{for } t_1^{(k)}\left(N_{\leq t}\right)<t\leq t_2^{(k)}\left(N_{\leq t}\right)\\
    \kappa_2^{(k)}\delta_{SE} &\text{for } t_2^{(k)}\left(N_{\leq t}\right)<t
    \end{cases}
\end{equation}
and equivalently for $\rate_{SI}^{(k)}$, possibly with different reduction factors $\kappa_i^{(k)}$ $i=1,2$. More generally, also the rate constants $\delta_{SE}$ may depend on the subpopulation $k$, but here we omit the corresponding indices for the purpose of simplicity.

\quad
In order to make our model more realistic, we include in each subpopulation $k$ a limited health care capacity given by a threshold $h_R^{(k)}$. 
We assume that the case fatality rate $\rate_{ID}^{(k)}\geq 0$ increases from a given value $\delta_{ID}$ to another value  $\tilde{\delta}_{ID}> \delta_{ID}$  if the number of infected individuals exceeds this threshold $h_R^{(k)}$, giving
\begin{equation*}
   \rate^{(k)}_{ID}(N):=
    \begin{cases}
    \delta_{ID} &\text{for } N_I^{(k)}\leq h_R^{(k)}\\
    \tilde{\delta}_{ID} &\text{for } N_I^{(k)}>h_R^{(k)}.
    \end{cases}
\end{equation*}
Vice versa, the recovery rate $\rate_{IR}$ is reduced in case of an exhausted health care capacity, such that 
\begin{equation*}
    \rate^{(k)}_{IR}(N):=
    \begin{cases}
    \delta_{IR} &\text{for } N_I^{(k)}\leq h_R^{(k)}\\
    \tilde{\delta}_{IR} &\text{for } N_I^{(k)}>h_R^{(k)}
    \end{cases}
\end{equation*}
for constants $\delta_{IR} > \tilde{\delta}_{IR} \geq 0$. Additionally, within each subpopulation $k$ we consider $\rate^{(k)}_{ID}+\rate^{(k)}_{IR}$ to be constant, i.e. $\delta_{IR}+\delta_{ID}=\tilde{\delta}_{IR}+\tilde{\delta}_{ID}$. 

Finally, we assume that exposed individuals develop symptoms after an incubation period of average length $\tau_{EI}>0$ and set
\begin{equation*}
   \rate^{(k)}_{EI}=\frac{1}{\tau_{EI}}
\end{equation*}
for all $k$. 
\subsubsection*{Concrete choice of rate constants for spatial transitions}
The global spatial transition rate functions between the subpopulations $\lambda_i^{(kl)}$ will depend on the local phases within each of the subpopulations.
More precisely, we define
\begin{equation*}
    \tau_1\left(N_{\leq t}\right):=\min \left\{t_1^{(1)}\left(N_{\leq t}\right),t_1^{(2)}\left(N_{\leq t}\right)\right\}
\end{equation*}
to be the first time that one of the subpopulations initiates the lock-down phase and 
\begin{equation*}
    \tau_2\left(N_{\leq t}\right):=\max \left\{t_2^{(1)}\left(N_{\leq t}\right),t_2^{(2)}\left(N_{\leq t}\right)\right\}
\end{equation*}
to be the first time that both subpopulations have ended the lock-down phase. Whenever in one of the subpopulations the strict measures are applied, the spatial transition rates are reduced by a factor $\kappa_1^{(kl)}\in(0,1)$. After the strict measures have ended in both populations, the spatial transition rates are scaled by a factor $\kappa_2^{(kl)}$, where  $\kappa_1^{(kl)}<\kappa_2^{(kl)}<1$. Thus, the  spatial transition rate functions are defined as
\begin{equation}\label{eq:HdepSpatial}
   \lambda_i^{(kl)}\left(N_{\leq t},t\right):=
    \begin{cases}
    \delta^{(kl)}, &\text{for } t\leq \tau_1\left(N_{\leq t}\right)\\
    \kappa_1^{(kl)}\delta^{(kl)} &\text{for } \tau_1\left(N_{\leq t}\right)<t\leq \tau_2\left(N_{\leq t}\right)\\
    \kappa_2^{(kl)}\delta^{(kl)} &\text{for }  \tau_2\left(N_{\leq t}\right)<t
    \end{cases}
\end{equation}
for $i\in \{S,E,R\}$. We assume that people with symptoms do not travel, i.e. we set $\lambda_I^{(kl)}=\lambda_D^{(kl)}=0$ for all $k,l$ independently of time.
\subsection{PDMM-based simulations of COVID-19 spreading}\label{subsec:Scenarios}
We simulate the dynamics for  model scenarios with differences in the infection and spatial transition dynamics. In particular, we compare the following three scenarios:
\begin{itemize}
    \item[] \textbf{Scenario 1:} Choose constant infection and spatial transition rates, which can be interpreted that no measures are implemented. 
    \item[] \textbf{Scenario 2:} Let the infection rate depend on the process history as defined in  \eqref{eq:HdepInfect}, but assume constant spatial transition rates between subpopulations. This corresponds to introducing local measures to control the infection dynamics within subpopulations, but no additional travel restrictions in between. 
    \item[] \textbf{Scenario 3:} Combine the measures, i.e., let both infection and spatial transition rates change according to the epidemic dynamics following the rules defined in \eqref{eq:HdepInfect} and \eqref{eq:HdepSpatial}.
\end{itemize}
\subsubsection*{Parameter choices}
Recently, lots of research has been done on inferring the parameters of COVID-19 dynamics from available data. However, for most parameters there is a wide range of estimates and thus the choices for a conceptual model can seem arbitrary. Since parameter estimation is not at the core of this manuscript, but our modeling approach is, we will choose the parameters based on a few recent publications \cite{byrne2020inferred,wang2020updated,linton2020incubation,rki2020:covfs,wom2020:dr}.\\
\textbf{Parameters for  status change $E \rightarrow I$:}
The average incubation period is estimated to be $5-6$ days \cite{linton2020incubation,rki2020:covfs}, so we choose $\tau_{EI}=5.5$. \\
\textbf{Parameters for status change $I\to R$ and $I\to D$:} 
The average time for transition from I to either R or D will be set to $14$ days \cite{linton2020incubation,rki2020:covfs}. For the infection fatality rate of our model, we will use the estimate from \cite{wom2020:dr}, which leads to the choices of  $\delta_{IR}=\frac{1-0.014}{14}$ and $\delta_{ID}=\frac{0.014}{14}$ for the recovery and case fatality rates of the populations. \revision{We choose $\tilde{\delta}_{ID}=3\delta_{ID}$ and $\tilde{\delta}_{IR}$ accordingly such that $\delta_{IR}+\delta_{ID}=\tilde{\delta}_{IR}+\tilde{\delta}_{ID}$ is fulfilled.}\\
\textbf{Parameters for status change $S\to E$:} \revision{Here we consider two reactions that can lead to the status change $S\to E$, namely $S+E \to 2E$ and $S+I \to E+I$. As discussed before, for simplicity, we assume that $\rate^{(k)}_{SE}=\rate^{(k)}_{SI}$ and thus  $\delta_{SI}=\delta_{SE}$.} Estimates for the initial reproduction number $R_0:=\frac{\delta_{SE}}{\delta_{IR}+\delta_{ID}}$ vary depending on the region of choice \cite{fernandez2020estimating} as well as on the estimation method  \cite{alimohamadi2020estimate}. This leads to a wide range of possible parameter choices that have the highest impact on the model outcome. For our model, we use $R_0=4.1$ which corresponds to the estimate for the New York City in  \cite{fernandez2020estimating}, assuming our subpopulations to be well-mixed and an urban area like NYC to meet this assumption. This choice of $R_0$ leads to the infection rate $\delta_{SE}:=\frac{4.1}{14}$. The remaining parameters will be  subject to changes due to different containment measures. During a period of strict measures phase we will reduce the infection rate to $10\%$ of the original value by setting a scaling factor $\kappa_1^{(k)}=0.1$ for $k=1,2$. In the phase of moderate measures we assume more interactions which lead to an increase of the infection rate. For illustration purposes, we consider different choices of moderate measures in each subpopulation, such that we set values of the infection rates to be $30\%$ and $40\%$ of the original $\delta_{SE}$, i.e. $\kappa_2^{(1)}=0.3$ and $\kappa_2^{(2)}=0.4$. For both subpopulations, the infection threshold for the first hitting time event is chosen to be $2\%$ of the initial  total population number $n_a^{(k)}$ in subpopulation $k$, i.e. $h_I^{(k)}=0.02\cdot n_a$ and the threshold for the capacity of the health care system is reached when $10\%$ of $n_a$ are infected, i.e. $h_R^{(k)}=0.1\cdot n_a^{(k)}$. \\
\textbf{Parameters for spatial transitons:} The spatial transition rates are chosen to be $\delta^{(kl)}=\delta^{(lk)}=0.0003$, \revision{which corresponds in our example to 3 out of 10000 agents transitioning per time unit and fits to the assumption of} a metastable setting with slow transitions between subpopulations compared to the infection dynamics within subpopulations. When at least one subpopulation is in the strict measures phase, we introduce travel restrictions by reducing the spatial transition rate to $5\%$ of the original value, i.e. we set $\kappa_1^{kl}=\kappa_1^{lk}=0.05$. When both subpopulations are in the moderate measures phase, we moderately relax the travel restrictions by increasing the spatial transition rate to $50\%$ of the original value, i.e. we set $\kappa_2^{kl}=\kappa_2^{lk}=0.5$.

\subsubsection*{Simulation results}

\begin{figure}[h!]
    \centering
    \begin{subfigure}{0.5\textwidth}
    \centering
        \includegraphics[width=\linewidth]{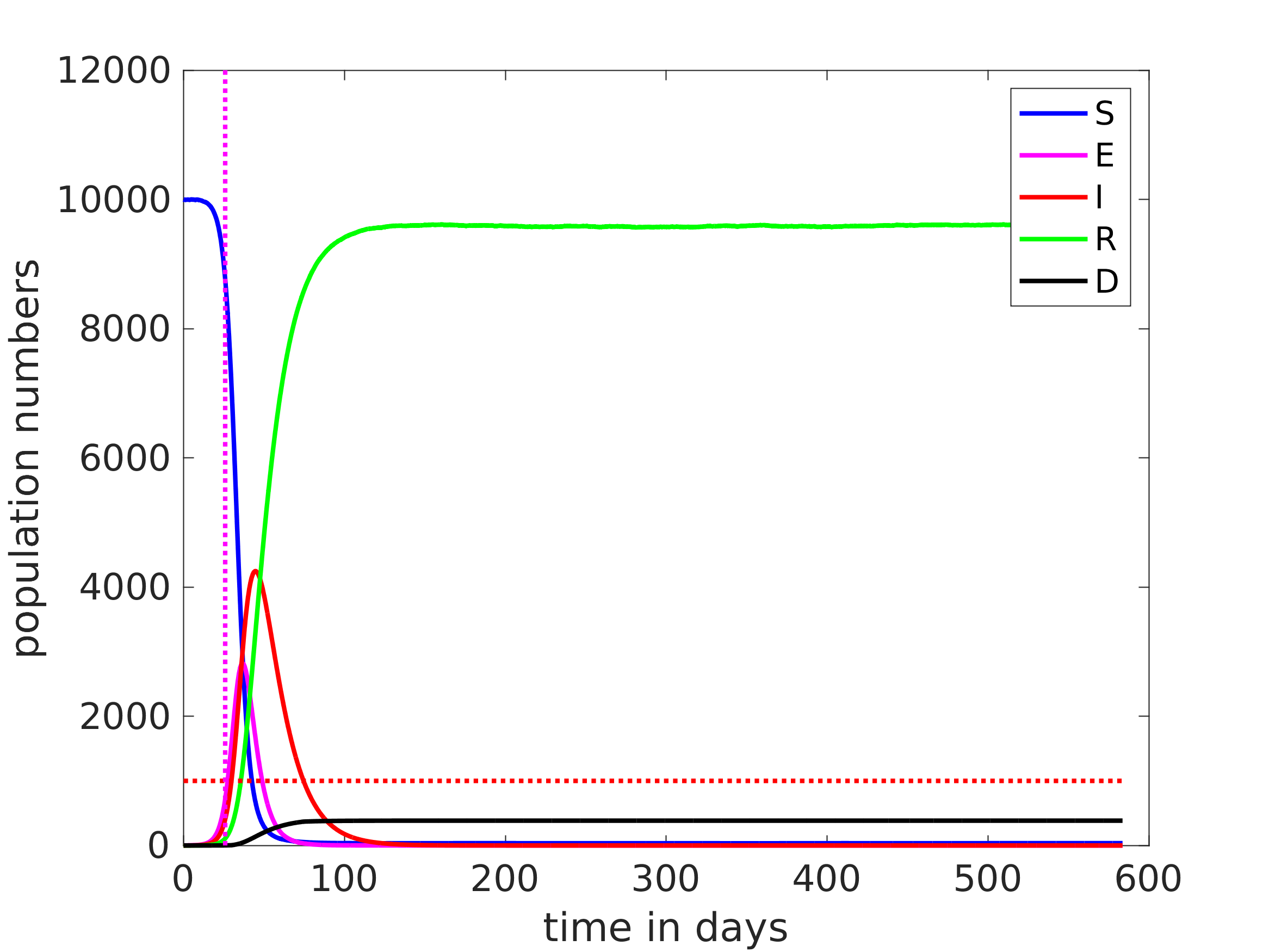}
        \caption{Subpopulation 1}
        
    \end{subfigure}%
    \begin{subfigure}{0.5\textwidth}
    \centering
        \includegraphics[width=\linewidth]{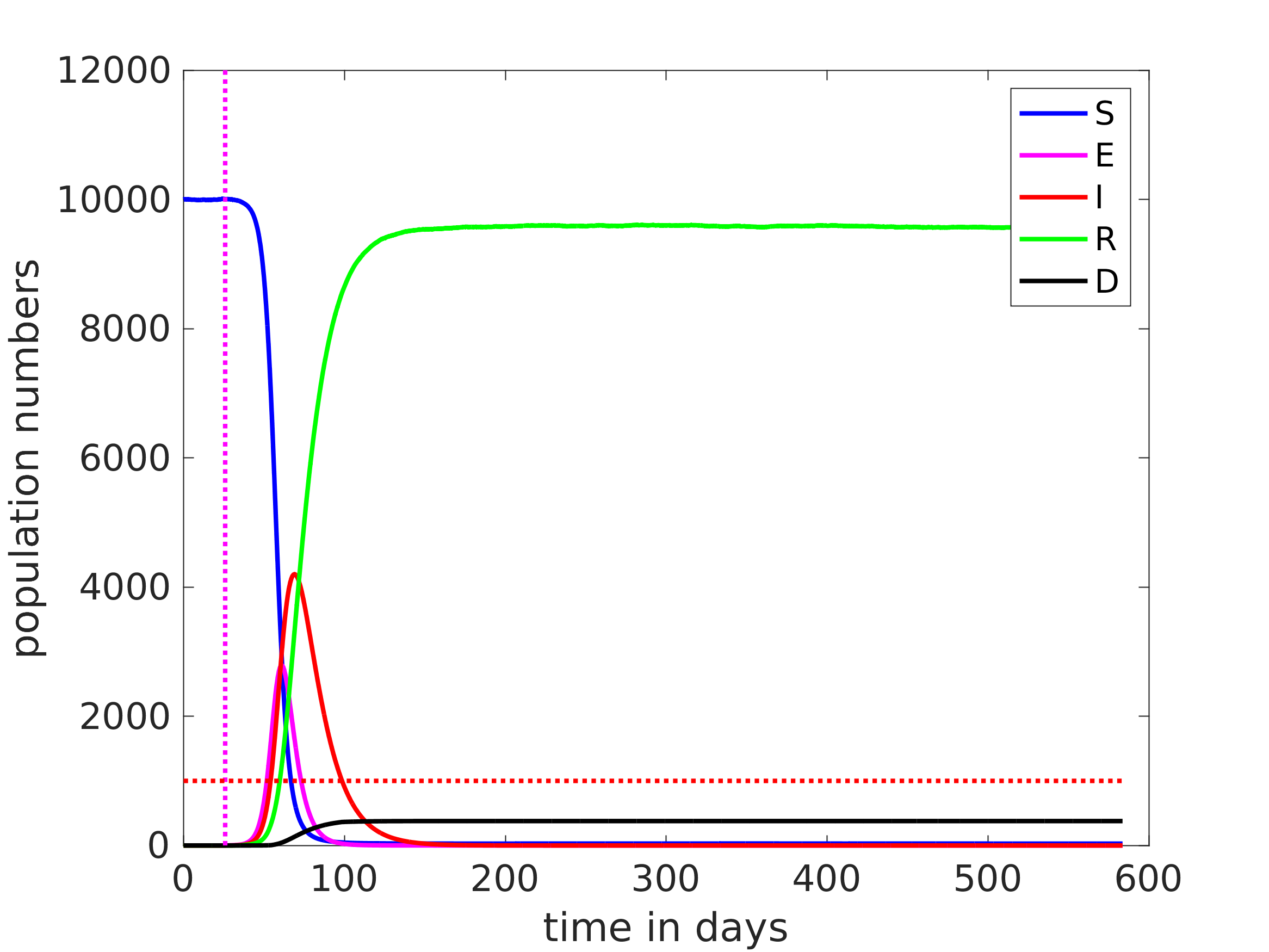}
        \caption{Subpopulation 2}
        
    \end{subfigure}
    
    \caption{Trajectory of a PDMM simulation for Scenario 1 (no measures). The dotted magenta line marks the critical transition event, the horizontal red line marks the threshold $h_R$ of health care capacity.}
    \label{fig:virusPDMM1}
\end{figure}

For the initial population sizes we choose the values $n_a^{(1)}=n_a^{(2)}= 10\,000$. We start with one member of the subpopulation 1 (SP1) being in status E and all other members of SP1 and subpopulation 2 (SP2) being in status S. The critical transition event is the first time when an individual with status E jumps from SP1 to SP2.
In Figure \ref{fig:virusPDMM1}, we see one outcome of Scenario 1, where no containment measures are taken. We observe one wave of infections in each subpopulation with the number of I cases quickly going above the threshold $h_R$ and staying there for a considerable amount of time. Almost all population members get infected. Due to the increased case fatality rate $\rate_{ID}$, in the end $3.8\%$ of the total population are in status D. This is a scenario that in reality should be avoided, for example by flattening the curve through the implementation of containment measures.

\begin{figure}[h!]
    \centering
    \begin{subfigure}{0.5\textwidth}
    \centering
        \includegraphics[width=\linewidth]{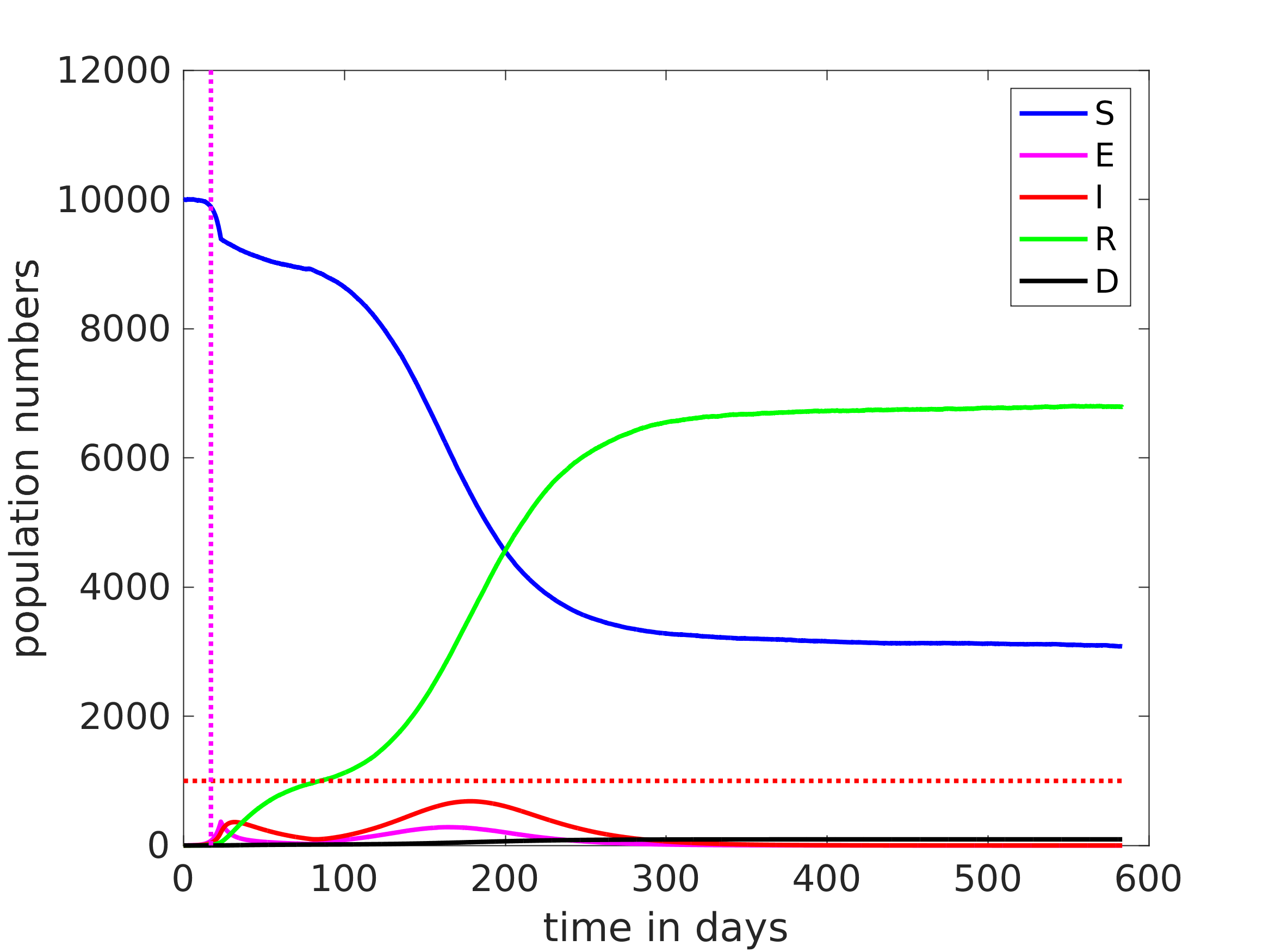}
        \caption{Subpopulation 1}
    
    \end{subfigure}%
    \begin{subfigure}{0.5\textwidth}
    \centering
        \includegraphics[width=\linewidth]{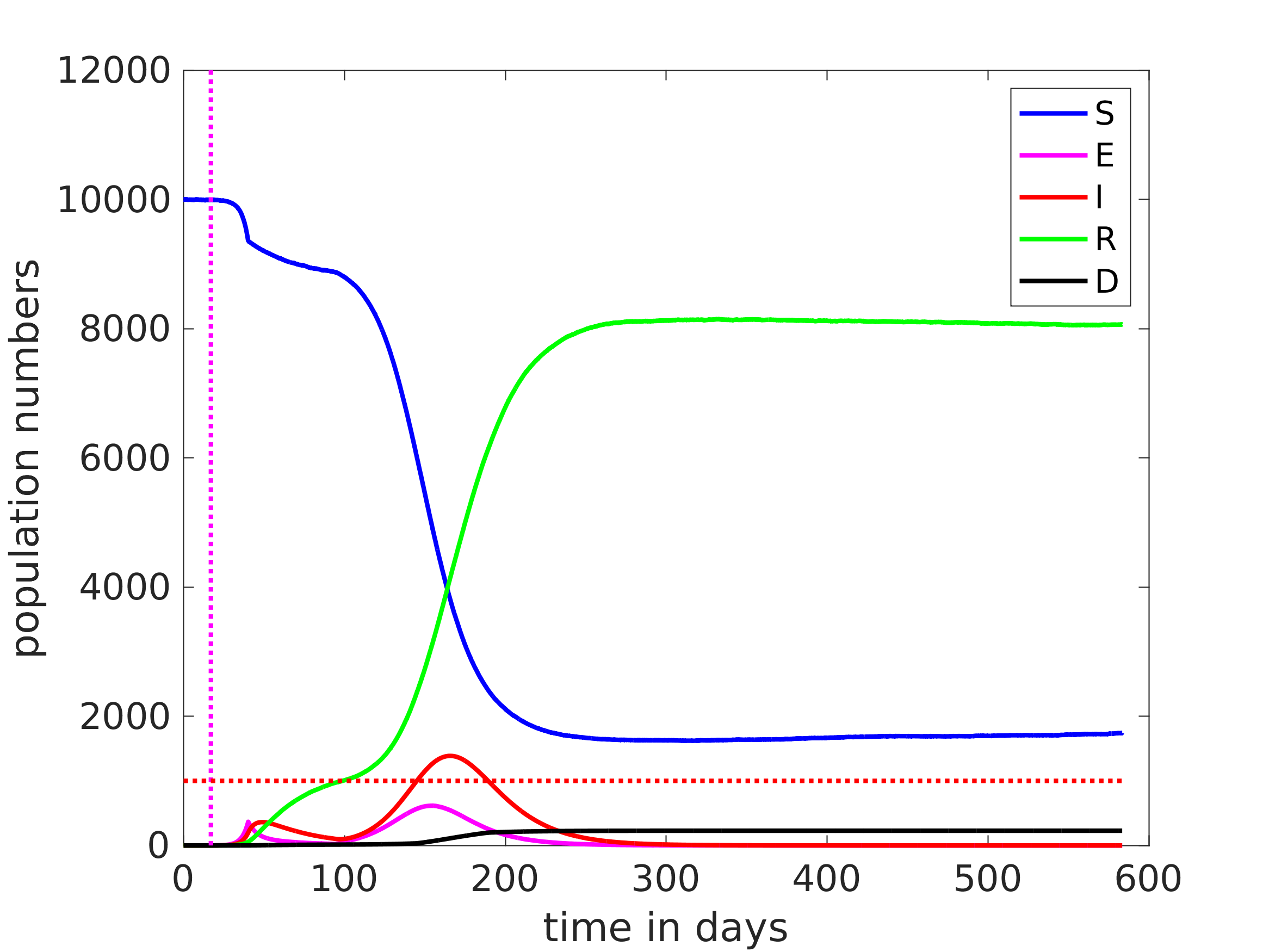}
        \caption{Subpopulation 2}
        
    \end{subfigure}
    
    \caption{Trajectory of a PDMM simulation for Scenario 2. 
    The dotted magenta line marks the critical transition event, the horizontal red line marks the threshold $h_R$.}
    \label{fig:virusPDMM3}
\end{figure}

As a result of the local measures that are present in Scenario 2, the curve of infections shows  two smaller waves instead of one large wave, see Figure \ref{fig:virusPDMM3}. In SP1 the number of infections stays below the threshold $h_R$ during the whole simulation period, while in SP2 the number of cases crosses $h_R$ during the peak of the second wave. This is due to a higher infection rate within SP2 in the phase of moderate measures, which leads to a higher total number of infections and more fatal cases in SP2 at the end of the simulation. Nevertheless, the outcome in both subpopulations is a much smaller number of I and D individuals than in Scenario 1. The same is true for Scenario 3 which has the same local measures. 

\begin{figure}[h!]
    \centering
    \begin{subfigure}{0.5\textwidth}
    \centering
        \includegraphics[width=\linewidth]{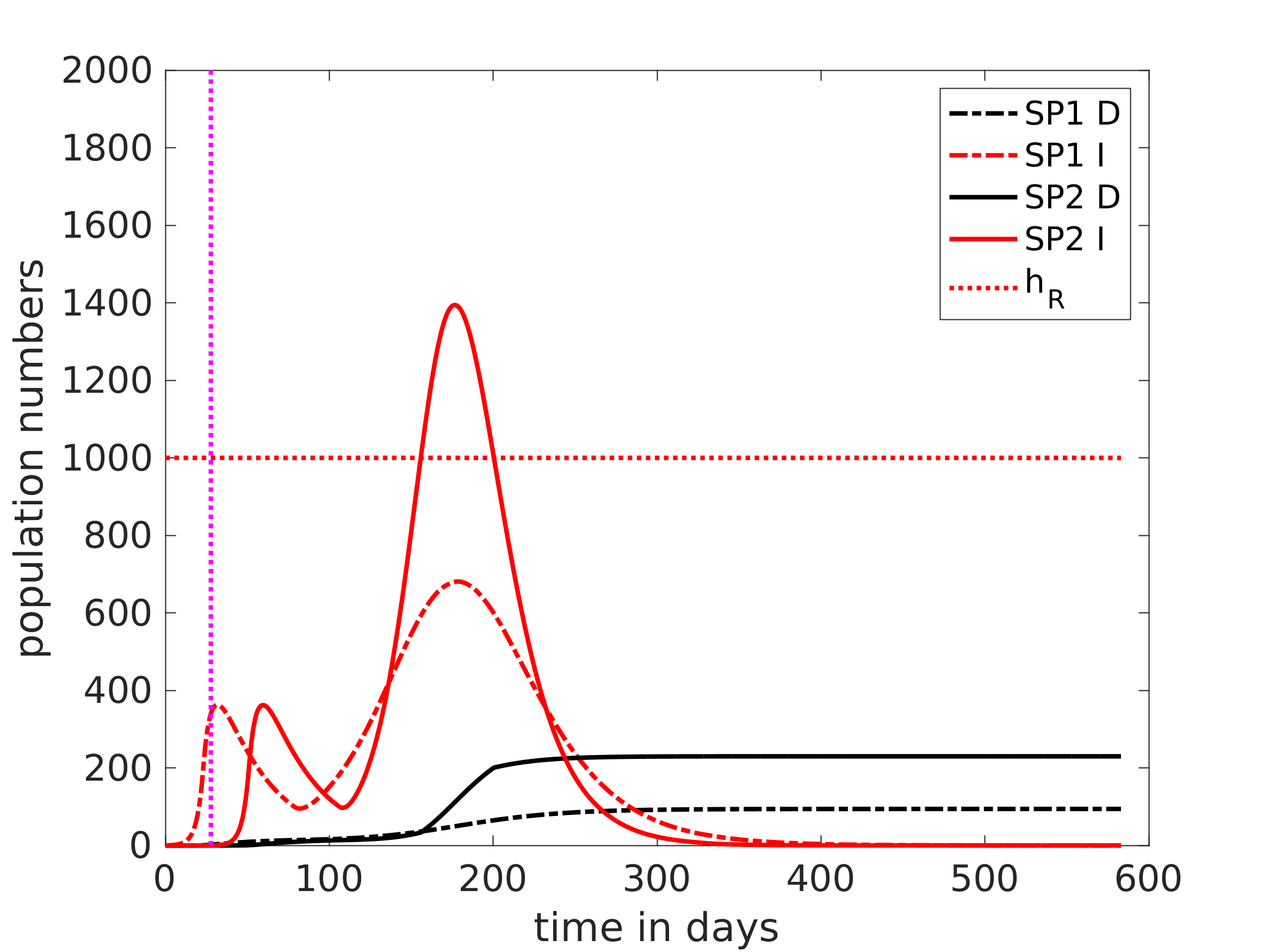}
        \caption{Scenario 2}
        
    \end{subfigure}%
    \begin{subfigure}{0.5\textwidth}
    \centering
        \includegraphics[width=\linewidth]{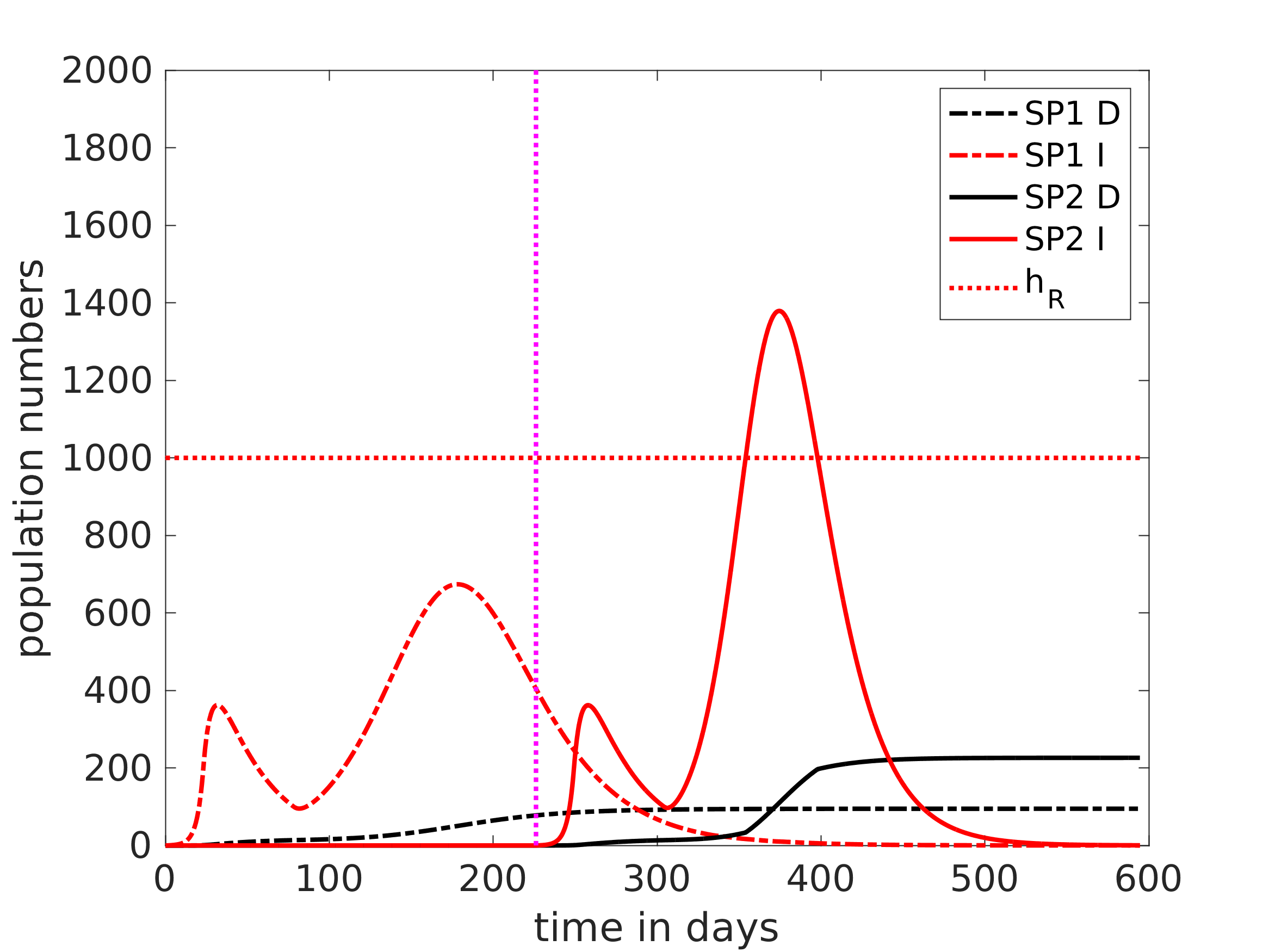}
        \caption{Scenario 3}
         \label{fig:virusPDMMcompb}
    \end{subfigure}
    
    \caption{Comparison between infection curves for single trajectories obtained in Scenarios 2 and 3. The vertical dotted magenta line marks the critical transition event, the horizontal red line marks the threshold $h_R$. The dashed lines refer to the development in the subpopulation 1 and solid lines to subpopluation 2.}
    \label{fig:virusPDMMcomp}
\end{figure}

Additionally, in Scenarios 2 and 3 the shape of the infection curves determined by the internal population dynamics is the same (see Figure \ref{fig:virusPDMMcomp}). However, the distribution of the critical transition time that starts the epidemic in SP2 is considerably different. Due to the introduction of travel restrictions in Scenario 3, we observe a later first infection in SP2 compared to the one from Scenario 2, see Figure \ref{fig:virusPDMMcompb}.

\begin{figure}[h!]
    \centering
    \begin{subfigure}{0.5\textwidth}
    \centering
        \includegraphics[width=\linewidth]{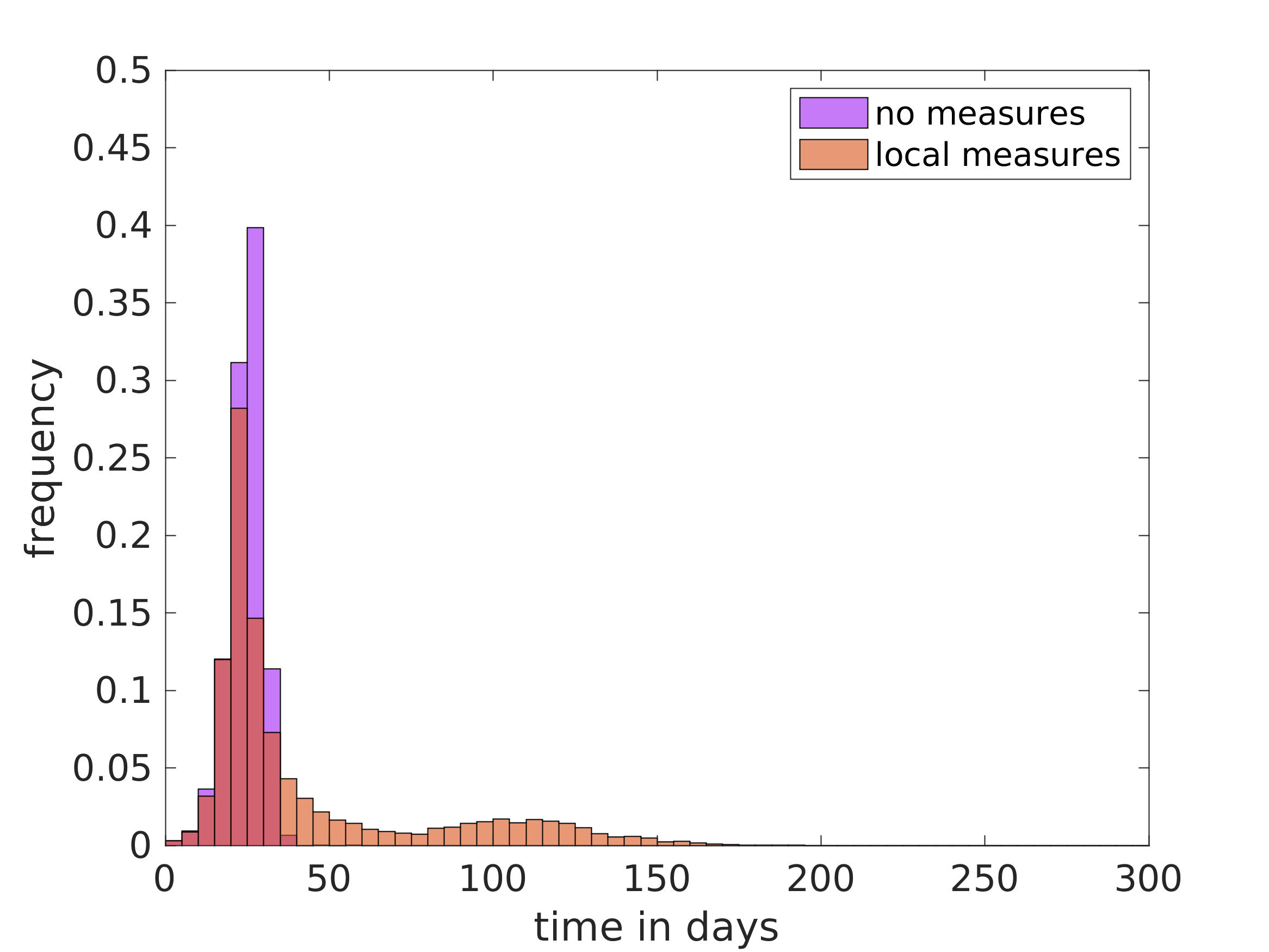}
        \caption{Scenario 1 and 2}
        
    \end{subfigure}%
    \begin{subfigure}{0.5\textwidth}
    \centering
        \includegraphics[width=\linewidth]{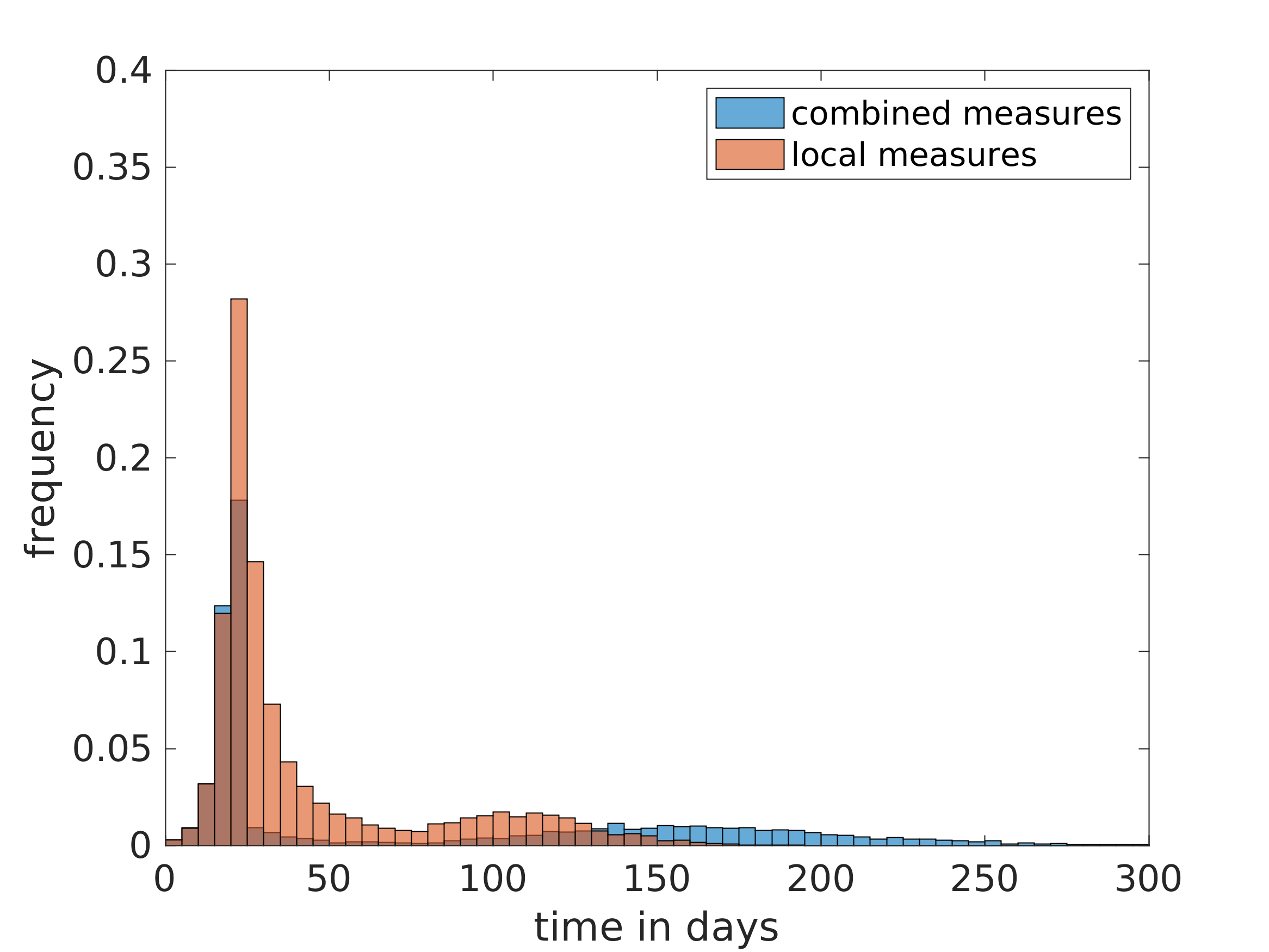}
        \caption{Scenario 2 and 3}
        
    \end{subfigure}
   
    \caption{Critical transition time distribution for different scenarios of the model. The number of MC-simulations for each scenario is $10000$.} 
    \label{fig:virusPDMMCTT}
\end{figure}

In order to compare the critical transition time distributions for different containment measures, we run $10000$ MC-simulations for each of the three scenarios, see Figure \ref{fig:virusPDMMCTT}.  For approximately one third of the simulations (no matter which scenario) the critical transition happens before the time $t_1^{(1)}$ when the measures are introduced in scenarios 2 and 3. After this time point we can observe the differences in the shapes of the critical transition time distributions. 
Namely, compared to Scenario 1, we observe for Scenario 2 a larger number of critical transitions  happening  later in time. 
This is due to the influence of the number of active cases in SP1, which is declining much faster in Scenario 2 compared to Scenario 1 due to the local measures. The mean time for the first infection  in SP2 is $24.5$ days for Scenario 1 and $43.9$ days for Scenario 2. Due to the reduced spatial transitions in Scenario 3, the probability for a critical transition is much smaller during the measures.  As a result, in $3989$ (out of $10000$) MC-simulations the critical transition event did not occur at all and the virus was successfully contained in SP1. Conditioned on the transition happening before the end of the simulation period, the mean critical transition time was $78.6$ days, which shows the benefit of the travel restrictions on the spreading dynamics.

In our conceptional model the outcome of the epidemic was improved by the implementation of measures. The spatial separation into multiple subpopulations has a high impact on the dynamics of the spreading process, especially when considering the scenario of combined measures where the spreading between populations could be delayed by a long time and sometimes even prevented.

\section{Conclusion}

In this work we introduced a hierarchy of modeling approaches for spatio-temporal dynamics of interacting agents. We showed how the stochastic meta\-population model can be derived from the underlying spatially continuous agent-based model by means of a Galerkin projection of the dynamics, considering  both  a full-partition approach and a core set approach. Especially, we specified the form of the projected generators and derived equations for the relation between the corresponding adoption rate functions for both first- and second-order status changes. In our guiding example, we analyze how the approximation error depends on the metastability of the dynamics.

Given the stochastic metapopulation model, we investigated a further approximation by piecewise-deterministic dynamics, where only the spatial transitions between metastable areas are modelled as stochastic jumps while the adoption dynamics within each subpopulation are approximated by deterministic dynamics. We analyzed the approximation quality as well as efficiency, finding that for large numbers of agents the PDMM delivers convincing results combined with enormous decrease in computational effort. 

Based on this insight, we formulated a PDMM for the spatio-temporal spreading dynamics of COVID-19 and analyzed the impact of different measures. We compared three basic scenarios with respect to critical transition times of the virus spreading between two subpopulations, finding out that the spatial separation into subpopulations can have a high impact on the epidemic spreading process,  e.g. by means of local measures and/or traveling restrictions.

Finally, in this paper we showed that, for large number of agents, stochastic metapopulation models and in particular PDMMs represent good approximations of ABMs which can be achieved with much less computational power, allowing for faster simulations of many different scenarios. Compared to other state of the art modeling approaches that were not discussed in detail in this manuscript, e.g. ODE-based models \cite{hethcote2000mathematics,alvarez2020simple,roda2020difficult,shao2020dynamic,chen2020time,ivorra2020mathematical}, in real world applications, SMMs and PDMMs appear to be better suited as they include both a certain level of stochasticity and a spatial resolution on a mesoscale and still have a reasonable computational cost.

\subsubsection*{Acknowledgements}
\revision{We would like to thank Francesco D'Amato for carefully reading the manuscript and giving helpful suggestions for revising the proofs. Furthermore, 
we wish to acknowledge the insightful comments of the two anonymous reviewers.
}
This research has been partially funded by the Deutsche Forschungsgemeinschaft (DFG, German Research Foundation) through grant CRC 1114/2 and under Germany’s Excellence Strategy -- The Berlin Mathematics Research Center MATH+ (EXC-2046/1 project ID:  390685689). 

\bibliography{bibliography}
\bibliographystyle{unsrt}
\newpage
\section{Appendix}\label{Appendix}

The following corollaries will be needed to show our main results. \revision{Again, we use the notation $e_\alpha$ for the $\alpha$th unit vector of $\R^{n_a}$, while $E_{i}^{(k)}$ denotes an \revision{$n_s\times m$} matrix with all entries zero except the entry at index $(k,i)$ which is one. }

\begin{cor}\label{CorAnsatzFunctions}
For any $N\in \mathbb{M}_{n_a}$ \revision{and given $i,j\in \mathbb{S}$}, it holds that
\begin{equation*}
\Phi_N(X,S+ie_\alpha-je_\alpha) = \Phi_{N+ E^{(k_\alpha)}_j-E^{(k_\alpha)}_i}(X,S)
\end{equation*} 
for each $\alpha\in\{1,...,n_a\}$ \revision{with $s_\alpha = j$ and  $N_i^{(k_\alpha)}>0$}.
\end{cor}

\revision{The condition $N_i^{(k_\alpha)}>0$ in Corollary \ref{CorAnsatzFunctions} is necessary to guarantee that it holds $N+ E^{(k_\alpha)}_j-E^{(k_\alpha)}_i\in \mathbb{M}_{n_a}$ such that  $\Phi_{N+ E^{(k_\alpha)}_j-E^{(k_\alpha)}_i}$ is actually defined.  In order to simplify the notation in all the following calculations, we extend the definition of $\Phi_N$ and set $\Phi_N(X,S):=0$ for $N\notin \mathbb{M}_{n_a}$. With this definition,  Corollary \ref{CorAnsatzFunctions} also works for  $N_i^{(k)}=0$ because both sides of the equation become zero.  }

\quad 

\noindent \textbf{Proof.} Choose $\ell\in \{1,...,n_s \}$ and $k\in \{1,...,m\}$. For $i\not=j$ we consider
\begin{eqnarray*}
\lefteqn{\delta_j(s_\alpha)\phi_{N_\ell^{(k)}}(X,S+ie_\alpha-je_\alpha)}\\
&\stackrel{\eqref{def_Phi}}{=}& \delta_j(s_\alpha)\cdot \delta_{N_\ell^{(k)}}\left(\sum_{\beta=1}^{n_a} \delta_{A_k}(x_\beta) \left\{\begin{array}{ll}
\delta_\ell(s_\beta) & \text{if } \beta\not=\alpha\\
\delta_\ell(i) & \text{if } \beta=\alpha\\
\end{array}\right.\right)\\
& = & \delta_j(s_\alpha)\cdot \delta_{N_\ell^{(k)}}\left( \left\{\begin{array}{ll}
\sum_{\beta=1}^{n_a} \delta_{A_k}(x_\beta)\delta_\ell(s_\beta) & \text{if } x_\alpha\not\in A_k\\
\sum_{\beta=1}^{n_a} \delta_{A_k}(x_\beta)\delta_\ell(s_\beta) +\delta_\ell(i)-\delta_\ell(j) & \text{if } x_\alpha\in A_k\\
\end{array}\right.
 \right).
\end{eqnarray*}
Thus, \revision{using definition \eqref{def_Phi}, we get for $k$ with $x_\alpha \notin A_k$:
\[ 
   \delta_j(s_\alpha)\phi_{N_\ell^{(k)}}(X,S+ie_\alpha-je_\alpha) = \delta_j(s_\alpha)\phi_{N_\ell^{(k)}}(X,S)
\]
for any $\ell$. For $k $ such that $x_\alpha \in A_k$, on the other hand, we distinguish between the following cases.
 For $\ell\not=i$ and $\ell\not=j$, it holds 
\begin{equation*}
   \delta_j(s_\alpha)\phi_{N_\ell^{(k)}}(X,S+ie_\alpha-je_\alpha) = \delta_j(s_\alpha)\phi_{N_\ell^{(k)}}(X,S), 
\end{equation*}
for $\ell=j$, we calculate
\begin{eqnarray*}
  \delta_j(s_\alpha)\phi_{N_j^{(k)}}(X,S+ie_\alpha-je_\alpha)
  & = &\delta_j(s_\alpha)\delta_{N_j^{(k)}}\left(\sum_{\beta=1}^{n_a}  \delta_{A_k}(x_\beta)\delta_j(s_\beta) -1\right)\\
    & = &\delta_j(s_\alpha)\delta_{N_j^{(k)}+1}\left(\sum_{\beta=1}^{n_a}  \delta_{A_k}(x_\beta)\delta_j(s_\beta) \right)\\
  & = &\delta_j(s_\alpha)\phi_{N_j^{(k)}+1}(X,S) 
  \end{eqnarray*}
and for $\ell=i$, we analogously get
\[ \delta_j(s_\alpha)\phi_{N_i^{(k)}}(X,S+ie_\alpha-je_\alpha) =  \delta_j(s_\alpha)\phi_{N_i^{(k)}-1}(X,S).\]
By definition of $E_i^{(k)}$, it is $N+ E^{(k)}_j-E^{(k)}_i$  the state where all  numbers $N_\ell^{(k)}$ stay the same except $N_i^{(k)}$, which is replaced be $N_i^{(k)}-1$, 
and $N_j^{(k)}$, which is replaced be $N_j^{(k)}+1$. 
Let $k_\alpha$ denote the index of the set $A_k$ for which $x_\alpha\in A_k$.
Then, combining the above calculations  and using the definition $\Phi_N=\prod_{k=1}^m \prod_{i=1}^{n_s}\phi_{N_i^{(k)}}$ of $\Phi_N$ given in \eqref{def_Phi1}, we obtain} 
\[ \Phi_N(X,S+ie_\alpha-je_\alpha) = \Phi_{N+ E^{(k_\alpha)}_j-E^{(k_\alpha)}_i}(X,S)\] 
for each $\alpha=1,...,n_a$ \revision{with $\delta_j(s_\alpha)=1$}. \hfill $\Box$

\quad

Now, using basic combinatorics we show that the following relation holds:

\begin{cor}\label{CorScalarProduct} For \revision{$N\in \mathbb{M}_{n_a}$ with $N_i^{(k)}>0$} and $M = N+ E^{(k)}_j-E^{(k)}_i$ it holds that
\begin{equation}\label{weight} 
\frac{\langle\Phi_M,\Phi_M\rangle}{\langle\Phi_N,\Phi_N\rangle} = \frac{N_i^{(k)}}{N_j^{(k)}+1}.
\end{equation}
\end{cor}

\noindent \textbf{Proof.} Using basic combinatorics, we obtain that
\begin{align}\label{eq:cor1}
\langle\Phi_N,\Phi_N\rangle  = & \frac{1}{(\mu(\X) n_s)^{n_a}} \int_{\X^{n_a}} \sum_{S\in \mathbb{S}^{n_a}} \Phi_N(X,S) \,dX  \nonumber \\
= & \frac{n_a!}{\prod_{\kappa,\ell}N_\ell^{(\kappa)}!} \prod_{\kappa,\ell} \Big(\frac{\mu(A_\kappa)}{n_s\mu(\X)}\Big)^{N_\ell^{(\kappa)}}.
\end{align}
This results from the multinomial distribution of $n_a$ particles into boxes $(\kappa,\ell)$, $\kappa=1,\ldots,m$, $\ell=1,\ldots,n_s$ with $N_\ell^{(\kappa)}$ particles each and the box probabilities $p_{\kappa,\ell}:=\frac{\mu(A_\kappa)}{n_s\mu(\X)}$. Then, for $M = N+ E^{(k)}_j-E^{(k)}_i$ by using equation \eqref{eq:cor1} we directly obtain \eqref{weight}.  \hfill $\Box$

\quad

\revision{Finally, given a linear operator $H:L^2(\mathbb{Y})\to L^2(\mathbb{Y})$, we derive the matrix representation of the  projected  operator $QHQ$ for the considered case of a full-partition projection. The proof is based on the analysis given in \cite{schutte2013metastability}.}

\revision{
\begin{cor}\label{CorGalerkin}
Given a linear operator $H:L^2(\mathbb{Y})\to L^2(\mathbb{Y})$, the Galerkin projection $QHQ$ with $Q$ defined in \eqref{def:Galerkin3} has the matrix representation $\hat{\mathcal{H}}$, where
\begin{equation}\label{A-matrix}
\hat{\mathcal{H}}_{NM}=\frac{\langle \Phi_M,H \Phi_N\rangle}{\langle \Phi_N,\mathbbm{1}\rangle}.
\end{equation}
\end{cor}
\noindent \textbf{Proof.} We construct a matrix representation for $QHQ$ with respect to the basis $\{\hat{\Phi}_N: N\in \mathbb{M}_{n_a}\}$ for probability densities given by 
\[ \hat{\Phi}_N := \frac{\Phi_N}{\langle \Phi_N,\mathbbm{1}\rangle}.\]
For a function $u = u(X,S,t)$ let 
\[ (Qu)(X,S,t) = \sum_{N\in\mathbb{M}_{n_a}} c_N(t) \hat{\Phi}_N(X,S),\]
for coefficients $c_N(t)\in \mathbb{R}$, 
be the basis representation of $Qu$. 
Then, 
\[ QHQu = \sum_{N\in\mathbb{M}_{n_a}} c_NQH\hat{\Phi}_N,\]
and with the formula for $Q$ given in \eqref{def:Galerkin3},
\[
     QH \hat{\Phi}_N = \sum_{M\in\mathbb{M}_{n_a}} \frac{\langle \Phi_M,H \hat{\Phi}_N\rangle}{\langle \Phi_M,\mathbbm{1} \rangle} \Phi_M.
\]
Putting this together, we obtain
\begin{align*}
   QHQu = & \sum_{M,N\in\mathbb{M}_{n_a}} c_N\frac{\langle \Phi_M,H \hat{\Phi}_N\rangle}{\langle \Phi_M,\mathbbm{1} \rangle} \Phi_M \\
    =  & \sum_{M,N\in\mathbb{M}_{n_a}} c_N\langle \Phi_M,H \hat{\Phi}_N\rangle \hat{\Phi}_M \\
    = & \sum_{M,N\in\mathbb{M}_{n_a}} c_N\frac{\langle \Phi_M,H \Phi_N\rangle}{\langle \Phi_N,\mathbbm{1}\rangle} \hat{\Phi}_M. 
\end{align*}
This means that the matrix $\hat{\mathcal{H}}$ is given by $\hat{\mathcal{H}}_{NM}=\frac{\langle \Phi_M,H \Phi_N\rangle}{\langle \Phi_N,\mathbbm{1}\rangle}$.  
Note that the coefficient $c_N(t)$ gives the probability for the state $N$ at time $t$, i.e., it refers to $P(N,t)$ in \eqref{discreteOperators}.\hfill $\Box$ \newline 
}

\revision{\begin{rem} In case of a non-full partition, i.e., for more general basis functions $\Phi_N$ that form a partition of unity but are not necessarily indicator functions, the result from Corollary \ref{CorGalerkin} can be extended by using Theorem 5.6 from  \cite{schutte2013metastability}. There, the adjoint operator $H^{adj}$ is considered, fulfilling $\langle H^{adj}\Phi_M,\Phi_N \rangle = \langle \Phi_M,H\Phi_N \rangle $. Translating their result to our generator $H$, the matrix representation of the projected generator $QHQ$ with $Q$ given by
\begin{equation}\label{def:Galerkin2}
    Qv = \sum_{M,N\in \mathbb{M}_{n_a}} (\mathcal{M}^{-1})_{MN} \langle \Phi_M,v\rangle \Phi_N,\quad  \mathcal{M}_{MN} := \langle \Phi_M,\Phi_N\rangle
\end{equation}
has the form $\hat{\mathcal{H}}\hat{\mathcal{M}}^{-1}$ with $\hat{\mathcal{H}}$ given in \eqref{A-matrix} and 
\[ \hat{\mathcal{M}}_{MN} := \frac{\langle \Phi_N,\Phi_M\rangle}{\langle \Phi_M,\mathbbm{1} \rangle }. \]
\end{rem}
}

\quad

\noindent  \revision{In the following, we will use these corollaries  for proving the main Theorems \ref{Thm0}-\ref{Thm2}. }

\subsection*{Proof of Theorem \ref{Thm0}}
\revision{
We first observe that for fixed $\alpha\in \{1,...,n_a\}$ and $i\in \mathbb{S}$ it holds  
\begin{equation}\label{split_alpha}
\delta_i(s_\alpha)\Phi_N(X,S) = \sum_{k=1}^{m}  \delta_k(x_\alpha)\delta_i(s_\alpha)\Phi_{N-E_i^{(k)}}(X^{\neg \alpha},S^{\neg \alpha}),
\end{equation}
where $S^{\neg \alpha}\in \mathbb{S}^{n_a-1}$ denotes the vector resulting from $S\in \mathbb{S}^{n_a}$ by skipping the entry $s_\alpha$, and   $X^{\neg \alpha}\in \mathbb{X}^{n_a}$ results from $X\in\mathbb{X}^{n_a}$ by skipping the entry $x_\alpha$.\footnote{\revision{Again, we use the extended definition $\Phi_N:=0 $ for $N\notin \mathbb{M}_{n_a}$, such that $\Phi_{N-E_i^{(k)}} $ is also defined in case of $N_i^{(k)}=0$.}}  }

\revision{Set $\mu_0:=\frac{1}{(\mu(\X) n_s)^{n_a}}$. Using the definition of $L$ given in \eqref{def:L}, we calculate \small
\begin{eqnarray*}
&&\langle \Phi_M,L\Phi_N\rangle \\
&=& \mu_0\sum_{S\in \mathbb{S}^{n_a}} \int_{\X^{n_a}} \Phi_M(X,S) \sum_{\alpha=1}^{n_a}\sum_{i=1}^{n_s} \delta_i(s_\alpha)L^{(\alpha)}_i \Phi_N(X,S)\, dX \\
&=&  \mu_0 \sum_{S\in \mathbb{S}^{n_a}}\int_{\X^{n_a}} \Phi_M(X,S) \sum_{\alpha=1}^{n_a}\sum_{i=1}^{n_s}\sum_{l=1}^m \delta_{A_l}(x_\alpha)\delta_i(s_\alpha)L^{(\alpha)}_i \Phi_N(X,S)\, dX \\
&\stackrel{\eqref{split_alpha}}{ =} & \mu_0\sum_{\alpha=1}^{n_a}  \sum_{i=1}^{n_s}\sum_{k,l=1}^m\sum_{S\in \mathbb{S}^{n_a}} \int_{\X^{n_a}} 
 \Phi_{M-E_i^{(l)}}(X^{\neg  \alpha },S^{\neg \alpha})\delta_{A_l}(x_\alpha)  \delta_i(s_\alpha) (L_{i}\delta_{A_k})(x_\alpha)\Phi_{N-E_i^{(k)}}(X^{\neg \alpha},S^{\neg \alpha})\, dX \\
& =&  \mu_0\sum_{\alpha=1}^{n_a} \sum_{i=1}^{n_s}\sum_{k,l=1}^m\sum_{S\in \mathbb{S}^{n_a}} \int_{\X^{n_a-1}} 
 \Phi_{M-E_i^{(l)}}(X^{\neg  \alpha },S^{\neg \alpha})\delta_i(s_\alpha)\int_\mathbb{X} \delta_{A_l}(x_\alpha)   (L_{i}\delta_{A_k})(x_\alpha) dx_\alpha \Phi_{N-E_i^{(k)}}(X^{\neg \alpha},S^{\neg \alpha})\, dX^{\neg \alpha} \\
& \stackrel{\eqref{lambda}}{=}& \mu_0\sum_{\alpha=1}^{n_a}  \sum_{i=1}^{n_s}\sum_{k,l=1}^m \sum_{S\in \mathbb{S}^{n_a}}\int_{\X^{n_a-1}}
 \Phi_{M-E_i^{(l)}}(X^{\neg  \alpha },S^{\neg \alpha}) \delta_i(s_\alpha)  \lambda_i^{(kl)}\int_\mathbb{X}\delta_{A_k}(x_\alpha)   dx_\alpha \Phi_{N-E_i^{(k)}}(X^{\neg \alpha},S^{\neg \alpha})\, dX^{\neg \alpha} \\
& =&\mu_0  \sum_{\alpha=1}^{n_a} \sum_{i=1}^{n_s}\sum_{k,l=1}^m \lambda_i^{(kl)}\sum_{S\in \mathbb{S}^{n_a}}\int_{\X^{n_a}}
 \Phi_{M-E_i^{(l)}+E_i^{(k)}}(X,S)\delta_i(s_\alpha)  \delta_{A_k}(x_\alpha)   \Phi_{N}(X,S)\, dX \\
 & =&\mu_0   \sum_{i=1}^{n_s}\sum_{k,l=1}^m \lambda_i^{(kl)}\sum_{S\in \mathbb{S}^{n_a}}\int_{\X^{n_a}}
 \Phi_{M-E_i^{(l)}+E_i^{(k)}}(X,S)\sum_{\alpha=1}^{n_a}\delta_i(s_\alpha)  \delta_{A_k}(x_\alpha)   \Phi_{N}(X,S)\, dX \\
  & \stackrel{(*)}{=}&\mu_0   \sum_{i=1}^{n_s}\sum_{k,l=1}^m \lambda_i^{(kl)}\sum_{S\in \mathbb{S}^{n_a}}\int_{\X^{n_a}}
 \Phi_{M-E_i^{(l)}+E_i^{(k)}}(X,S)N_i^{(k)}   \Phi_{N}(X,S)\, dX \\
  & =&  \sum_{i=1}^{n_s}\sum_{k,l=1}^m \lambda_i^{(kl)}N_i^{(k)}\left\langle  \Phi_{M-E_i^{(l)}+E_i^{(k)}},   \Phi_{N} \right \rangle 
\end{eqnarray*}
\normalsize
where $(*)$ results from the fact that it holds $\sum_{\alpha=1}^{n_a}\delta_i(s_\alpha)  \delta_{A_k}(x_\alpha)=N_i^{(k)}$ for all $(X,S)$ with $\Phi_N(X,S)\neq 0$.  
Assume that it holds $M=N+E_i^{(l)}-E_i^{(k)}$ for certain $k,l,i$, $k\neq l$. Then all summands  are zero except one summand, and we obtain
\begin{eqnarray*} 
\langle \Phi_M,L\Phi_N\rangle 
& =&  \lambda_i^{(kl)} N_i^{(k)} \langle  \Phi_{N},   \Phi_{N} \rangle .
\end{eqnarray*}
For the case $M=N$ we need $k=l$ such that $\Phi_{M-E_i^{(l)}+E_i^{(k)}}=\Phi_M=\Phi_N$, and obtain 
\[\langle \Phi_M,L\Phi_N\rangle = \sum_{i=1}^{n_s} \sum_{k=1}^m \lambda_i^{(kk)} N_i^{(k)}  \langle \Phi_{N}, \Phi_{N}\rangle. \]
For other combinations of $M,N$ the overall sum is zero. } In total, we get 
\begin{eqnarray*}
	&& \langle \Phi_M,L\Phi_N\rangle \\
	& = & \left\{\begin{array}{ll}
		N_i^{(k)} \lambda_i^{(kl)}  \langle \Phi_{N}, \Phi_{N}\rangle, & \mathrm{if} \;M=N+E_i^{(l)}-E_i^{(k)},\; k\neq l\\
		\sum_{i=1}^{n_s} \sum_{k=1}^m N_i^{(k)} \lambda_i^{(kk)} \langle \Phi_{N}, \Phi_{N}\rangle, & \mathrm{if} \;M=N, \\
		0, & \mathrm{otherwise,}
	\end{array}
	\right.\\
	& = & \left\{\begin{array}{ll}
		N_i^{(k)} \lambda_i^{(kl)} \langle \Phi_{N}, \Phi_{N}\rangle, & \mathrm{if} \;M=N+E_i^{(l)}-E_i^{(k)},\;  k\neq l\\
		-\sum_{i=1}^{n_s} \sum_{\substack{k,l=1}{ l\neq k}}^m N_i^{(k)} \lambda_i^{(kl)} \langle \Phi_{N}, \Phi_{N}\rangle, & \mathrm{if} \;M=N, \\
		0, & \mathrm{otherwise.}
	\end{array}
	\right.
\end{eqnarray*}
where for the case $M=N$ we used $\lambda_i^{(kk)}=-\sum_{l\neq k}\lambda_i^{(kl)}$. \revision{By means of Corollary \ref{CorGalerkin} we can have to divide by $\langle \Phi_{N}, \mathbbm{1}\rangle=\langle \Phi_{N}, \Phi_{N}\rangle$ and obtain 
\begin{eqnarray*}
\hcL_{NM}	& = & \left\{\begin{array}{ll}
		N_i^{(k)} \lambda_i^{(kl)} , & \mathrm{if} \;M=N+E_i^{(l)}-E_i^{(k)},\;  k\neq l\\
		-\sum_{i=1}^{n_s} \sum_{\substack{k,l=1}{ l\neq k}}^m N_i^{(k)} \lambda_i^{(kl)} , & \mathrm{if} \;M=N, \\
		0, & \mathrm{otherwise.}
	\end{array}
	\right.
\end{eqnarray*}
for the entries of $\hcL$.} \hfill $\Box$

\subsection*{Proof of Theorem \ref{Thm1}}

On the basis of Corollary \ref{CorAnsatzFunctions} we see that the action of the ABM generator $G$ on an individual indicator ansatz function can be written as
\begin{eqnarray*}
\lefteqn{G\Phi_N(X,S)}\\
& = & \sum_{i,j=1}^{n_s}\sum_{\alpha=1}^{n_a}\left(- f_{ij}^{(\alpha)}(X,S)\Phi_N(X,S) + f_{ij}^{(\alpha)}(X,S+ie_\alpha-je_\alpha)\Phi_{N+ E^{(k_\alpha)}_j-E^{(k_\alpha)}_i}(X,S)\right),
\end{eqnarray*}
with the consequence that
\begin{equation}\label{decomposition}
\hG_{MN}:=\langle\Phi_M,G\Phi_N\rangle
 =  -\hG_{1,M,N} + \hG_{2,M,N},
\end{equation}
where 
\begin{eqnarray*}
\hG_{1,M,N} 
& := & \left\langle\Phi_M,\sum_{i,j=1}^{n_s}\sum_{\alpha=1}^{n_a} f_{ij}^{(\alpha)}\Phi_N\right\rangle\\
& = & \revision{\mu_0\sum_{S\in\mathbb{S}^{n_a}} \int_{\X^{n_a}} \Phi_M(X,S) \sum_{i,j=1}^{n_s}\sum_{\alpha=1}^{n_a} f_{ij}^{(\alpha)}(X,S) \Phi_{N}(X,S)dX }\\
\hG_{2,M,N} 
& := & \revision{\left\langle\Phi_M,\sum_{i,j=1}^{n_s}\sum_{\alpha=1}^{n_a}f_{ij}^{(\alpha)}(X,S+ie_\alpha-je_\alpha) \Phi_{N + E_j^{(k_\alpha)} - E_i^{(k_\alpha)}}\right\rangle }\\
& = & \mu_0\sum_{S\in\mathbb{S}^{n_a}}\int_{\X^{n_a}}  \Phi_M(X,S) \sum_{i,j=1}^{n_s}\sum_{\alpha=1}^{n_a} f_{ij}^{(\alpha)}(X,S+ie_\alpha-je_\alpha) \Phi_{N+ E^{(k_\alpha)}_j-E^{(k_\alpha)}_i}(X,S)dX \end{eqnarray*}
and again $\mu_0:=\frac{1}{(\mu(\X) n_s)^{n_a}}$. 
\revision{We compute \small
\begin{eqnarray*}
\lefteqn{\hG_{1,M,N} } \\
& = & \mu_0\sum_{S\in\mathbb{S}^{n_a}} \int_{\X^{n_a}} \Phi_M(X,S) \sum_{i,j=1}^{n_s}\sum_{\alpha=1}^{n_a} f_{ij}^{(\alpha)}(X,S) \Phi_{N}(X,S)dX \\
& = & \mu_0\sum_{S\in\mathbb{S}^{n_a}} \int_{\X^{n_a}} \Phi_M(X,S) \sum_{i,j=1}^{n_s}\sum_{\alpha=1}^{n_a} \delta_i(s_\alpha)\gamma_{ij}(x_\alpha) \Phi_N(X,S)\, dX \\
& = &\mu_0\sum_{\alpha=1}^{n_a}\sum_{i,j=1}^{n_s}\sum_{k=1}^m \sum_{S\in\mathbb{S}^{n_a}} \int_{\X^{n_a}}  \Phi_M(X,S)  \delta_i(s_\alpha)\delta_{A_k}(x_\alpha)\gamma_{ij}(x_\alpha) \Phi_N(X,S)\, dX \\
& = &\mu_0\sum_{\alpha=1}^{n_a}\sum_{i,j=1}^{n_s}\sum_{k=1}^m \sum_{S\in\mathbb{S}^{n_a}} \int_{\X^{n_a-1}}  \Phi_{M-E_i^{(k)}}(X^{\neg \alpha},S^{\neg \alpha}) \delta_i(s_\alpha) \int_X\delta_{A_k}(x_\alpha)\gamma_{ij}(x_\alpha) dx_\alpha \Phi_{N-E_i^{(k)}}(X^{\neg \alpha},S^{\neg \alpha})\, dX^{\neg \alpha} \\
& \stackrel{\eqref{gamma_k}}{=} &\mu_0\sum_{\alpha=1}^{n_a}\sum_{i,j=1}^{n_s}\sum_{k=1}^m \sum_{S\in\mathbb{S}^{n_a}} \int_{\X^{n_a-1}}  \Phi_{M-E_i^{(k)}}(X^{\neg \alpha},S^{\neg \alpha}) \delta_i(s_\alpha) \gamma_{ij}^{(k)} \int_X\delta_{A_k}(x_\alpha) dx_\alpha \Phi_{N-E_i^{(k)}}(X^{\neg \alpha},S^{\neg \alpha})\, dX^{\neg \alpha} \\
& = &\mu_0\sum_{\alpha=1}^{n_a}\sum_{i,j=1}^{n_s}\sum_{k=1}^m \gamma_{ij}^{(k)}\sum_{S\in\mathbb{S}^{n_a}} \int_{\X^{n_a}}  \Phi_M(X,S)  \delta_i(s_\alpha)\delta_{A_k}(x_\alpha) \Phi_N(X,S)\, dX \\
& = &\mu_0\sum_{i,j=1}^{n_s}\sum_{k=1}^m \gamma_{ij}^{(k)}\sum_{S\in\mathbb{S}^{n_a}} \int_{\X^{n_a}}  \Phi_M(X,S) \sum_{\alpha=1}^{n_a} \delta_i(s_\alpha)\delta_{A_k}(x_\alpha) \Phi_N(X,S)\, dX \\
& \stackrel{(*)}{=} &\mu_0 \sum_{i,j=1}^{n_s}\sum_{k=1}^m\gamma_{ij}^{(k)} \sum_{S\in\mathbb{S}^{n_a}}\int_{\X^{n_a}} \Phi_M(X,S)  N_i^{(k)} \Phi_N(X,S)\, dX \\
& = & \sum_{i,j=1}^{n_s}\sum_{k=1}^m \gamma_{ij}^{(k)} N_i^{(k)} \langle \Phi_M , \Phi_N \rangle\\
& = & \left\{
\begin{array}{ll}
\sum_{k=1}^m\sum_{i,j=1}^{n_s}\gamma_{ij}^{(k)} N_i^{(k)}\langle \Phi_N , \Phi_N \rangle,\quad & \mathrm{if} \; M = N, \\
0, & \mathrm{otherwise}.
\end{array}
\right.
\end{eqnarray*}
\normalsize
In $(*)$ we used that it holds
\revision{$\sum_{\alpha=1}^{n_a} \delta_i(s_\alpha)\delta_{A_k}(x_\alpha)= N_i^{(k)} $}
for all $(X,S)$ with $\Phi_N(X,S)\neq 0$. 
Analogously, we calculate the non-diagonal entries, setting $\revision{\gamma}_{ii}(x)=0$, such that we can sum over all $i,j$:  \small
\begin{eqnarray*}
\lefteqn{ \hG_{2,M,N}} \\
& = &\mu_0\sum_{i,j=1}^{n_s} \sum_{S\in\mathbb{S}^{n_a}}\int_{\X^{n_a}}  \Phi_M(X,S) \sum_{\alpha=1}^{n_a} \delta_j(s_\alpha) \gamma_{ij}(x_\alpha) \Phi_{N + E_j^{(k_\alpha)} - E_i^{(k_\alpha)}}(X,S)\, dX \\
& = &\mu_0\sum_{\alpha=1}^{n_a} \sum_{i,j=1}^{n_s}\sum_{k=1}^m \sum_{S\in\mathbb{S}^{n_a}}\int_{\X^{n_a}}  \Phi_M(X,S) \delta_j(s_\alpha)\delta_{A_k}(x_\alpha) \gamma_{ij}(x_\alpha) \Phi_{N+ E_j^{(k)} - E_i^{(k)} }(X,S)\, dX \\
& = &\mu_0\sum_{\alpha=1}^{n_a} \sum_{i,j=1}^{n_s}\sum_{k=1}^m \sum_{S\in\mathbb{S}^{n_a}}\int_{\X^{n_a-1}}  \Phi_{M-E_j^{(k)}}(X^{\neg \alpha},S^{\neg \alpha}) \delta_j(s_\alpha)\int_\X\delta_{A_k}(x_\alpha) \gamma_{ij}(x_\alpha) dx_\alpha \Phi_{N - E_i^{(k)} }(X^{\neg \alpha},S^{\neg  \alpha})\, dX^{\neg \alpha} \\
& = &\mu_0\sum_{\alpha=1}^{n_a} \sum_{i,j=1}^{n_s}\sum_{k=1}^m \sum_{S\in\mathbb{S}^{n_a}}\int_{\X^{n_a-1}}  \Phi_{M-E_j^{(k)}}(X^{\neg \alpha},S^{\neg \alpha}) \delta_j(s_\alpha)\gamma_{ij}^{(k)}\int_\X\delta_{A_k}(x_\alpha)  dx_\alpha \Phi_{N - E_i^{(k)} }(X^{\neg \alpha},S^{\neg  \alpha})\, dX^{\neg \alpha} \\
& = &\mu_0 \sum_{i,j=1}^{n_s}\sum_{k=1}^m\gamma_{ij}^{(k)} \sum_{S\in\mathbb{S}^{n_a}}\int_{\X^{n_a-1}}  \Phi_{M}(X,S) \sum_{\alpha=1}^{n_a}\delta_j(s_\alpha)\delta_{A_k}(x_\alpha)   \Phi_{N+E_j^{(k)} - E_i^{(k)} }(X,S)\, dX \\
& \stackrel{(**)}{=} &\mu_0 \sum_{i,j=1}^{n_s}\sum_{k=1}^m\gamma_{ij}^{(k)} \sum_{S\in\mathbb{S}^{n_a}}\int_{\X^{n_a-1}}  \Phi_{M}(X,S) (N_j^{(k)}+1) \Phi_{N+E_j^{(k)} - E_i^{(k)} }(X,S)\, dX \\
& = & \sum_{i,j=1}^{n_s}\sum_{k=1}^m \gamma_{ij}^{(k)} (N_j^{(k)}+1) \left\langle \Phi_M , \Phi_{N+ E_j^{(k)} - E_i^{(k)} } \right\rangle\\
& = & \left\{\begin{array}{ll}
\gamma_{ij}^{(k)} (N_j^{(k)}+1)\langle \Phi_M,\Phi_M\rangle, & \mathrm{if} \;M=N+ E_{j}^{(k)}-E_{i}^{(k)} ,\;  i\neq j,\\
0, & \mathrm{otherwise,}
 \end{array}
 \right. \\
 & \stackrel{\eqref{weight} }{=} & \left\{\begin{array}{ll}
\gamma_{ij}^{(k)} N_i^{(k)}\langle \Phi_N,\Phi_N\rangle, & \mathrm{if} \;M=N+ E_{j}^{(k)}-E_{i}^{(k)} ,\;  i\neq j,\\
0, & \mathrm{otherwise.}
 \end{array}
 \right.
\end{eqnarray*} \normalsize
In $(**)$ we used that it holds $\sum_{\alpha=1}^{n_a} \delta_j(s_\alpha)\delta_{A_k}(x_\alpha)=N_j^{(k)}+1$ 
for all $(X,S)$ with $\Phi_{N + E_j^{(k)}- E_i^{(k)} }(X,S)\neq 0$.}

\revision{Combining the diagonal and non-diagonal part  and using Corollary \ref{CorGalerkin}, we obtain
\[
\hcG_{NM} =\left\{\begin{array}{ll}
\gamma_{ij}^{(k)}N_i^{(k)},\quad & \mathrm{if} \; M = N+ E_{j}^{(k)} - E_{i}^{(k)} ,\;  i\neq j,\\
-\sum_{i,j=1}^{n_s}\sum_{k=1}^m \gamma_{ij}^{(k)} N_i^{(k)}, & \mathrm{if} \;M=N,\\
0, & \mathrm{otherwise}
\end{array}
\right.
\]
for the entries of the matrix $\hcG$.}\hfill $\Box$

\subsection*{Proof of Theorem \ref{Thm2}}

\revision{At first, we observe that for each $i,j\in \mathbb{S}$, $S\in \mathbb{S}^{n_a}$, $\alpha\in \{1,...,n_a\}$, and $M\in \mathbb{M}_{n_a}$ it holds
\begin{align}\label{use_b_kl}\begin{split}
& \delta_i(s_\alpha)\delta_j(s_\beta)\int_{\X^{n_a}}   \delta_{A_k}(x_\alpha)\delta_{A_l}(x_\beta) d_r(x_\alpha,x_\beta)\Phi_{M}(X,S)\, dX \\
 = \; & b_{kl}\delta_i(s_\alpha)\delta_j(s_\beta)\int_{\X^{n_a}}   \delta_{A_k}(x_\alpha)\delta_{A_l}(x_\beta) \Phi_{M}(X,S)\, dX \end{split}
\end{align}
This can be seen by the following calculation: \small
\begin{eqnarray*}
\lefteqn{\delta_i(s_\alpha)\delta_j(s_\beta)\int_{\X^{n_a}}   \delta_{A_k}(x_\alpha)\delta_{A_l}(x_\beta) d_r(x_\alpha,x_\beta)\Phi_{M}(X,S)\, dX } \\
& = & \delta_i(s_\alpha)\delta_j(s_\beta)\int_{\X^{n_a-2}}\int_{\X^2}   \delta_{A_k}(x_\alpha)\delta_{A_l}(x_\beta) d_r(x_\alpha,x_\beta)dx_1dx_2\Phi_{M-E_i^{(k)}-E_j^{(l)}}(X,S)\, dX \\
& \stackrel{\eqref{b_kl}}{=} & \delta_i(s_\alpha)\delta_j(s_\beta)\int_{\X^{n_a-2}}b_{kl}\int_{\X^2}   \delta_{A_k}(x_\alpha)\delta_{A_l}(x_\beta) x_1dx_2\Phi_{M-E_i^{(k)}-E_j^{(l)}}(X,S)\, dX \\
& = & b_{kl}\delta_i(s_\alpha)\delta_j(s_\beta)\int_{\X^{n_a}}   \delta_{A_k}(x_\alpha)\delta_{A_l}(x_\beta) x_1dx_2\Phi_{M}(X,S)\, dX. 
\end{eqnarray*}
\normalsize
}

\revision{We use the same decomposition \eqref{decomposition} as in the proof before. 
Let the rate function $\rate_{ij}$ be given by $\rate_{ij}(x_1,x_2)=c_{ij}\cdot d_r(x_1,x_2)$ for  $c_{ij}\geq 0$, see definition \eqref{Doi}, where $c_{ii}=0$. Note that by setting $c_{ii}=0$ we can take the sum over all $\alpha,\beta$ without the condition $\beta \neq \alpha$, because it holds $\delta_i(s_\alpha)\delta_j(s_\beta)c_{ij}=0$ for all $i,j$ in case of $\alpha=\beta$. We compute
\begin{eqnarray*}
\lefteqn{ \hG_{1,M,N}} \\
& = & \mu_0\sum_{S\in\mathbb{S}^{n_a}} \int_{\X^{n_a}} \Phi_M(X,S) \sum_{i,j=1}^{n_s}\sum_{\alpha=1}^{n_a} f_{ij}^{(\alpha)}(X,S) \Phi_{N}(X,S)dX \\
& = &\mu_0\sum_{\alpha,\beta=1}^{n_a}\sum_{i,j=1}^{n_s}c_{ij} \sum_{S\in\mathbb{S}^{n_a}}\int_{\X^{n_a}}  \Phi_M(X,S)  \delta_i(s_\alpha) \delta_j(s_\beta) d_r(x_\alpha,x_\beta)\Phi_N(X,S)\, dX. 
\end{eqnarray*}
As we have $\Phi_M\cdot \Phi_N=0$ for $M\neq N$, we can follow $\hG_{1,M,N}=0$ for $M\neq N$. For $M=N$, on the other hand, we have $\Phi_M\cdot \Phi_N = \Phi_N$, such that we can skip $\Phi_M$ and get \small
\begin{eqnarray*}
\lefteqn{ \hG_{1,M,N}} \\
& = &\mu_0\sum_{\alpha,\beta=1}^{n_a} \sum_{i,j=1}^{n_s}c_{ij} \sum_{S\in\mathbb{S}^{n_a}}\int_{\X^{n_a}}  \delta_i(s_\alpha) \delta_j(s_\beta) d_r(x_\alpha,x_\beta)\Phi_N(X,S)\, dX \\
& = &\mu_0\sum_{\alpha,\beta=1}^{n_a} \sum_{i,j=1}^{n_s}\sum_{k,l=1}^m c_{ij} \sum_{S\in\mathbb{S}^{n_a}}\int_{\X^{n_a}} \delta_i(s_\alpha) \delta_j(s_\beta)\delta_{A_k}(x_\alpha)\delta_{A_l}(x_\beta) d_r(x_\alpha,x_\beta)\Phi_N(X,S)\, dX \\
& \stackrel{\eqref{use_b_kl}}{=} &\mu_0\sum_{\alpha,\beta=1}^{n_a} \sum_{i,j=1}^{n_s}\sum_{k,l=1}^m c_{ij}b_{kl} \sum_{S\in\mathbb{S}^{n_a}}\int_{\X^{n_a}} \delta_i(s_\alpha) \delta_j(s_\beta)\delta_{A_k}(x_\alpha)\delta_{A_l}(x_\beta) \Phi_N(X,S)\, dX \\
& = &\mu_0 \sum_{i,j=1}^{n_s}\sum_{k,l=1}^m c_{ij}b_{kl} \sum_{S\in\mathbb{S}^{n_a}}\int_{\X^{n_a}}\sum_{\alpha,\beta=1}^{n_a} \delta_i(s_\alpha) \delta_j(s_\beta)\delta_{A_k}(x_\alpha)\delta_{A_l}(x_\beta) \Phi_N(X,S)\, dX \\
& \stackrel{(*)}{=} &\mu_0 \sum_{i,j=1}^{n_s}\sum_{k,l=1}^m c_{ij}b_{kl} \sum_{S\in\mathbb{S}^{n_a}}\int_{\X^{n_a}} N_i^{(k)}N_j^{(l)} \Phi_N(X,S)\, dX \\
& = &\sum_{i,j=1}^{n_s}\sum_{k,l=1}^m c_{ij}b_{kl} N_i^{(k)}N_j^{(l)} \langle\Phi_N,\Phi_N \rangle, 
\end{eqnarray*} \normalsize
where  $(*)$ is true because of  $\sum_{\alpha=1}^{n_a} \delta_i(s_\alpha) \delta_{A_k}(x_\alpha)=N_i^{(k)} $ and $\sum_{\beta=1}^{n_a} \delta_j(s_\beta) \delta_{A_l}(x_\beta)=N_j^{(l)} $ for all $(X,S)$ with $\Phi_N(X,S) \neq0$. }

\revision{
 Using 
\[ f_{ij}^{(\alpha)}(X,S+ie_\alpha - je_\alpha) = c_{ij} \delta_j(s_\alpha) \sum_{\beta \neq \alpha} d_r(x_\alpha,x_\beta)\delta_j(s_\beta)\] 
with $c_{ij}=0$ for $i=j$,
we analogously  get for the non-diagonal entries:
\small
\begin{eqnarray*}
\lefteqn{\hG_{2,M,N}} \\
& = &\mu_0\sum_{i,j=1}^{n_s}c_{ij} \sum_{S\in\mathbb{S}^{n_a}}\int_{\X^{n_a}}  \Phi_M(X,S) \sum_{\alpha=1}^{n_a} \delta_j(s_\alpha) \sum_{\substack{\beta=1\\ \beta\neq\alpha}}^{n_a}\delta_j(s_\beta) d_r(x_\alpha,x_\beta)\Phi_{N+E_j^{(k_\alpha)}-E_i^{(k_\alpha)}}(X,S)\, dX \\
& = &\mu_0\sum_{\alpha=1}^{n_a} \sum_{\substack{\beta=1\\ \beta\neq\alpha}}^{n_a}\sum_{i,j=1}^{n_s}\sum_{k,l=1}^{m}c_{ij} \sum_{S\in\mathbb{S}^{n_a}}\int_{\X^{n_a}}  \Phi_M(X,S)  \delta_j(s_\alpha)\delta_j(s_\beta)\delta_{A_k}(x_\alpha)\delta_{A_l}(x_\beta) d_r(x_\alpha,x_\beta)\Phi_{N+E_j^{(k)}-E_i^{(k)}}(X,S)\, dX .
\end{eqnarray*}
\normalsize
In case of  $M\neq N+E_j^{(k)}-E_i^{(k)}$ for all $i,j =1,...,n_s$ and all $k=1,...,m$  we get $\hG_{2,M,N}=0$ because it holds $\Phi_M  \cdot \Phi_{N+E_j^{(k)}-E_i^{(k)}}=0$ for all $i,j,k$. For $M=N$, we have $\hG_{2,M,N}=0$ because of $c_{ii}=0$. If, on the other hand $M=N+E_j^{(k)}-E_i^{(k)}$ for some $i\neq j$ and $k\in \{1,...,m\}$, we have  $\Phi_M  \cdot \Phi_{N+E_j^{(k)}-E_i^{(k)}}=\Phi_{N+E_j^{(k)}-E_i^{(k)}}$ and 
\footnotesize
\begin{eqnarray*}
\lefteqn{ \hG_{2,M,N}} \\
& = &\mu_0\sum_{\alpha=1}^{n_a}  \sum_{\substack{\beta=1\\ \beta\neq\alpha}}^{n_a}\sum_{l=1}^m c_{ij} \sum_{S\in\mathbb{S}^{n_a}}\int_{\X^{n_a}}   \delta_j(s_\alpha)\delta_j(s_\beta)\delta_{A_k}(x_\alpha)\delta_{A_l}(x_\beta) d_r(x_\alpha,x_\beta)\Phi_{N+E_j^{(k)}-E_i^{(k)}}(X,S)\, dX\\
& \stackrel{\eqref{use_b_kl}}{=} &\mu_0\sum_{\alpha=1}^{n_a}  \sum_{\substack{\beta=1\\ \beta\neq\alpha}}^{n_a}\sum_{l=1}^m c_{ij}b_{kl} \sum_{S\in\mathbb{S}^{n_a}}\int_{\X^{n_a}}   \delta_j(s_\alpha)\delta_j(s_\beta) \delta_{A_k}(x_\alpha)\delta_{A_l}(x_\beta)\Phi_{N+E_j^{(k)}-E_i^{(k)}}(X,S)\, dX\\
& = &\mu_0\sum_{l=1}^m c_{ij}b_{kl} \sum_{S\in\mathbb{S}^{n_a}}\int_{\X^{n_a}}  \sum_{\alpha=1}^{n_a}  \sum_{\substack{\beta=1\\ \beta\neq\alpha}}^{n_a} \delta_j(s_\alpha)\delta_j(s_\beta) \delta_{A_k}(x_\alpha)\delta_{A_l}(x_\beta)\Phi_{N+E_j^{(k)}-E_i^{(k)}}(X,S)\, dX\\
& \stackrel{(**)}{=} &\mu_0\sum_{l=1}^m c_{ij}b_{kl} \sum_{S\in\mathbb{S}^{n_a}}\int_{\X^{n_a}} (N_j^{(k)}+1)N_j^{(l)}\Phi_{N+E_j^{(k)}-E_i^{(k)}}(X,S)\, dX\\
&  = &\sum_{l=1}^m c_{ij}b_{kl}(N_j^{(k)}+1) N_j^{(l)} \langle\Phi_{N+E_j^{(k)}-E_i^{(k)}}, \Phi_{N+E_j^{(k)}-E_i^{(k)}} \rangle \\
&  \stackrel{\eqref{weight}}{=} &\sum_{l=1}^mc_{ij} N_i^{(k)} N_j^{(l)} b_{kl}\langle\Phi_{N}, \Phi_N \rangle.
\end{eqnarray*} \normalsize
In $(**)$  we used that 
\[ \sum_{\alpha=1}^{n_a}\sum_{\substack{\beta=1\\ \beta\neq\alpha}}^{n_a} \delta_j(s_\alpha) \delta_{A_k}(x_\alpha)\delta_j(s_\beta) \delta_{A_l}(x_\beta)=(N_j^{(k)}+1)N_j^{(l)} \]
holds for all $(X,S)$ with $\Phi_{N+E_j^{(k)}-E_i^{(k)}}(X,S)\neq 0$.
For $k\neq l$ this is clear because for these $(X,S) $ there must be $N_j^{(k)}+1$ agents of status $j$ located in subset $A_k$ and  $N_j^{(l)}$  agents  of status $j$ in subset $A_l\neq A_k$. For $l=k$ we obtain $(N_j^{(k)}+1)N_j^{(k)} $ which is the number of pairs of different agents 
both of status $j$ and being located in $A_k$. 
}

\revision{
Dividing by $\langle \Phi_N,\mathbbm{1}\rangle=\langle \Phi_N,\Phi_N\rangle$ (see again Corollary \ref{CorGalerkin}) and combining diagonal and non-diagonal entries, we find that the matrix $\hcG$ has the entries
\[
\hcG_{NM} =\left\{\begin{array}{ll}
\sum_{l=1}^mc_{ij}b_{kl}N_i^{(k)}N_j^{(l)},\quad & \mathrm{if} \; M = N+ E_{j}^{(k)} - E_{i}^{(k)} ,\;  i\neq j,\\
-\sum_{i,j=1}^{n_s}\sum_{k,l=1}^mc_{ij}b_{kl} N_i^{(k)}N_j^{(l)}, & \mathrm{if} \;M=N,\\
0, & \mathrm{otherwise}.
\end{array}
\right.
\]
By using definitions \eqref{eps}, \eqref{hatgamma} and \eqref{2ndOrderProp} we complete the proof.} \hfill $\Box$

\end{document}